\documentclass[review]{elsarticle}

\usepackage{graphicx}
\usepackage{amsfonts, amsmath, amssymb}
\usepackage{bm}
\usepackage{cancel}
\usepackage{color}
\usepackage{enumerate}
\usepackage{float}
\usepackage{hyperref}
\usepackage[english]{babel}
\usepackage[margin=1in]{geometry}
\usepackage{mathtools}
\usepackage{setspace}
\usepackage{subfigure}
\usepackage{lineno}

\renewcommand \d [2]{\frac{{\rm d} #1}{{\rm d} #2}}
\newcommand \D [2]{\frac{\partial #1}{\partial #2}}

\renewcommand{\vec}[1]{\bm{\mathrm{#1}}}
\newcommand{\V}[1]{\bm{\mathrm{#1}}}

\def \div{\nabla \cdot \mbox{}}
\def \grad{\nabla}
\def \lap{\nabla^2}

\def \x{\vec{x}}
\def \q{\vec{q}}
\def \n{\vec{n}}
\def \r{\vec{r}}
\def \u{\vec{u}}
\def \uS{\vec{u}_S}
\def \e{\vec{e}}

\def \I{\vec{I}}
\def \F{\vec{F}}
\def \N{\vec{N}}
\def \T{\vec{T}}
\def \U{\vec{U}}

\def \cM{\vec{\mathcal{M}}}

\def \Sb{S_b}
\def \Sbt{S_b(t)}

\def  \Scv{S_\text{CV}}
\def  \Scvt{S_\text{CV}(t)}
\def \Vt{V(t)}
\def \Vmt{V_m(t)}
\def \Vbt{V_b(t)}
\def \Vb{V_b}
\def  \Vcv{V_\text{CV}}
\def  \Vcvt{V_\text{CV}(t)}

\def \cF{\vec{\mathcal{F}}}

\def \F{\vec{F}}

\def \I{\vec{I}}
\def \Ib{\I_{\text{b}}}

\def \Mb{M_{\text{b}}}
\def \N{\vec{N}}

\def \Omegat{\Omega(t)}

\def \cP{{\mathcal{P}}}
\def \vcP{\vec{\mathcal{P}}}
\def \Pb{\vec{\mathcal{P}}_\text{b}}
\def \R{\vec{R}}
\def \U{\vec{U}}
\def \Ub{\U_{\text{b}}}

\def \Uk{\U_{\text{k}}}

\def \Ur{\U_{\text{r}}}
\def \W{\vec{W}}

\def \Wr{\W_{\text{r}}}

\def \X{\vec{X}}

\def \cB{\mathcal{B}}

\def \cG{\mathcal{G}}

\def \cI{\vec{\mathcal{I}}}

\def \cR{{\mathcal{R}}}
\def \cL{{\mathcal{L}}}
\def \vcL{\vec{\mathcal{L}}}
\def \cS{\vec{\mathcal{S}}}
\def \cJ{\vec{\mathcal{J}}}
\def \cT{\mathcal{T}}

\def \e{\vec{e}}

\def \f{\vec{f}}

\def \half{\frac{1}{2}}
\def \3half{\frac{3}{2}}

\def \n{\vec{n}}

\def \nref{n_{\text{ref}}}

\def \s{\vec{s}}
\def \S{\vec{S}}
\def \u{\vec{u}}

\def \vomega{\vec{\omega}}

\def \x{\vec{x}}
\def \xL{x_L}
\def \yL{y_L}
\def \xU{x_U}
\def \yU{y_U}

\def \Re{\text{Re}}
\def \div{\nabla \cdot \mbox{}}
\def \grad{\nabla}
\def \lap{\nabla^2}

\def \dt{\Delta t}
\def \dx{\Delta x}
\def \dy{\Delta y}
\def \dz{\Delta z}
\def \delV{\Delta V}

\def \Ds{{\mathrm d}\s}
\def \Dels{\Delta \vec{s}}
\def \dS{\,\mathrm{dS}}
\def \dV{\,\mathrm{dV}}
\def \Dx{{\mathrm d}\x}

\def \dt{\Delta t}
\def \dx{\Delta x}

\def \ndot{\n \cdot}
\def \Ndot{\N \cdot}
\def \rcross{\r \wedge}
\def \Rcross{\R \wedge}

\newcommand{\REVIEW}[1]{{#1}}

\begin{document}

\begin{frontmatter}
	
\title{A moving control volume approach to computing hydrodynamic forces and torques on immersed bodies}

\author[Northwestern1]{Nishant Nangia}
\author[LBNL]{Hans Johansen}
\author[Northwestern1,Northwestern2]{Neelesh A. Patankar}
\author[LBNL]{Amneet Pal Singh Bhalla\corref{mycorrespondingauthor}}
\ead{apbhalla@lbl.gov}

\address[Northwestern1]{Department of Engineering Sciences and Applied Mathematics, Northwestern University, Evanston, IL}
\address[Northwestern2]{Department of Mechanical Engineering, Northwestern University, Evanston, IL}
\address[LBNL]{Applied Numerical Algorithms Group, Lawrence Berkeley National Laboratory, Berkeley, CA}
\cortext[mycorrespondingauthor]{Corresponding author}

\begin{abstract}
We present a moving control volume (CV) approach to computing hydrodynamic forces 
and torques on complex geometries. The method requires surface and volumetric 
integrals over a simple and regular Cartesian box that moves with an arbitrary 
velocity to enclose the body at all times. The moving box is aligned with Cartesian 
grid faces, which makes the integral evaluation straightforward in an  
immersed boundary (IB) framework. Discontinuous and noisy derivatives of velocity 
and pressure at the fluid-structure interface are avoided and far-field (smooth) velocity and 
pressure information is used. We re-visit the approach to compute 
hydrodynamic forces and torques through force/torque balance equation 
in a Lagrangian frame that some of us took in a prior work (Bhalla et al., J Comp Phys, 2013). 
We prove the equivalence of the two approaches for IB methods, thanks to 
the use of Peskin's delta functions. Both approaches are able to suppress 
spurious force oscillations and are in excellent agreement, as expected 
theoretically. Test cases ranging from Stokes to high Reynolds number regimes 
are considered. We discuss regridding issues for the moving CV method in an adaptive 
mesh refinement (AMR) context. The proposed moving CV method is not 
limited to a specific IB method and can also be used, for example, with embedded 
boundary methods.
\end{abstract}


\end{frontmatter}


\section{Introduction}

Fluid-structure interaction (FSI) problems involving moving bodies
is a challenging area in the computational fluid dynamics 
field that has vested the interest of researchers for several decades.
FSI modeling has traditionally been carried out in two ways: the 
\emph{body-fitted} mesh approach using unstructured grids~\cite{Hu01,Kern06} 
and the Cartesian grid approach based on the fictitious domain 
method~\cite{RGlowinski99,NAPatankar00}. Although the
body-fitted mesh approach to FSI resolves the fluid-structure interface 
sharply, it requires complex mesh management infrastructure along 
with high computational costs for solving linear equations. The fictitious 
domain method on the other hand extends the fluid equations inside 
the structure along with some additional non-zero body forcing term. 
As a result, regular Cartesian grids and fast linear solvers such as fast Fourier transform
can be used to solve the common momentum equation of the continua. 
However, fictitious domain methods tend to smear the fluid-structure
interface and hence reduce the solution accuracy near the interface.
 
One widely used fictitious domain approach to FSI is the 
immersed boundary (IB) method~\cite{Peskin02} which was 
originally proposed by Peskin in the context of cardiac flows~\cite{Peskin72b}. 
The main advantage to IB methods is that they do not require 
a body-fitted mesh to model the structure. The immersed body 
is allowed to freely cut the background Cartesian mesh, 
making the IB method easy to implement within an existing 
incompressible flow solver. A Lagrangian force density is 
computed on structure nodes, which is then transferred to the 
background grid via regularized delta functions. The use of 
regularized delta functions diffuses the interface and smears 
it over a number of grid cells proportional to the width of the 
delta function. This often leads to discontinuous
or noisy derivatives of the velocity and pressure field required 
to compute surface traction.

The original IB method uses an
explicit time stepping scheme and fiber elasticity 
to compute the additional body forcing inside 
the region occupied by the structure. The method works fairly well 
for soft elastic structures, but incurs severe time step restrictions 
if the stiffness of the material is substantially increased to model 
rigid bodies. Specialized versions of the original IB method have 
been developed to model rigid and stiff bodies in an efficient manner 
like the implicit IB method~\cite{NewrenEtAl07,Newren08,RDGuy12,Bhalla16,Mori08}, 
the direct forcing method~\cite{Uhlmann05,Bhalla13,Bhalla14}, and the 
fully constrained IB method~\cite{Taira07,Kallemov16,Usabiaga17}. The 
implicit IB method requires special solvers like algebraic~\cite{Ceniceros09} 
and geometric multigrid~\cite{NewrenEtAl07,Bhalla16} to treat the stiff 
elastic forces implicitly. The direct forcing and fully constrained IB 
methods that are primarily used for modeling rigid bodies impose 
rigidity constraint through Lagrange multipliers. The direct forcing 
method approximates the Lagrange multiplier with suitable 
penalty term and solves the fluid-structure equations in a fractional 
time stepping scheme. This is in contrast to fully constrained IB 
method where the fluid velocity, pressure, and Lagrange 
multipliers are solved together. Apart from
the obvious advantage of imposing the rigidity constraint exactly rather
approximately, the fully constrained IB method can be used for Stokes 
flow where fractional time stepping schemes do not work. However,
special preconditioners are required to solve for the Lagrange multipliers
exactly~\cite{Kallemov16,Usabiaga17}, which makes solving the system 
costly for big three dimensional volumetric bodies. 

A common feature of IB methods that are based on Peskin's IB 
approach~\cite{Peskin02} is that they do not require rich geometric 
information like surface elements and normals to compute the 
Lagrangian force density. Structure node position 
(and possibly the node connectivity) information suffices 
to compute Lagrangian forces. However, there 
are many other versions of the IB method and Cartesian grid 
based methods that require additional geometric information to 
sharply resolve the fluid-structure interface.
Examples include the immersed finite element method~\cite{LZhang04,WKLiu06,Heltai12,Griffith17}, 
the immersed interface method~\cite{LaiLi01,LiLai01}, 
the ghost-fluid method~\cite{YHTseng03}, and the cut-cell embedded 
boundary method~\cite{RMittal08,Uday01,Treb15}. For such approaches, 
the availability of surface elements and normals along with sharp 
interface resolution makes the surface traction computation easier 
and smooth (for at least smooth problems). 

Often times, the net hydrodynamic forces and torques on immersed 
bodies are desired rather than point-wise traction values. Force/torque 
balance equations can be used to compute the hydrodynamic 
force/torque contribution instead of directly integrating surface 
traction. In the context of ``Peskin-like" IB methods, this implies that 
one can essentially eliminate surface mesh generation as a 
post-processing step if only net hydrodynamic forces and torques 
are desired. This has another added advantage of not using noisy 
derivatives of velocity and pressure at the interface for evaluating 
hydrodynamic forces and torques. 

In this work we analyze and compare two approaches to computing net hydrodynamic forces and
 torques on an immersed body. In the first approach we use the Reynolds transport theorem 
(RTT) to convert the traction integral over an irregular body surface to a traction integral 
over a regular and simple Cartesian box (that is aligned with grid faces). The RTT is proposed for 
a moving control volume that translates with an arbitrary velocity to enclose the
immersed body at all times. We refer to this approach as the moving 
CV approach. In the context of locally refined grids, the moving 
control volume can span a hierarchy of grid levels. 
In the second approach, hydrodynamic forces and torques are computed using 
inertia and Lagrange multipliers (approximate or exact) defined in the body region. 
We refer to this approach as the LM approach and has been used before 
in~\cite{Bhalla13}. For IB methods, we show that both approaches are 
equivalent. This is due to a special property of Peskin's delta functions  
that makes Lagrangian and Eulerian force density equivalent~\cite{Peskin02}.
We show that both these approaches give smooth forces
and suppress spurious oscillations that arise by directly integrating spatial 
pressure and velocity gradients over the immersed body as reported in the
literature~\cite{Lee11}. 

Although application of the RTT on a stationary control volume to evaluate 
hydrodynamic force is a well known result~\cite{Kundu14,Pozrikidis11}; 
its extension to moving control volumes was first proposed by Flavio
Noca in 1997~\cite{Noca97}. They were motivated by the task of evaluating 
net hydrodynamic force on a moving bluff body using DPIV data 
from experiments. Flavio has also proposed force expressions 
that eliminate the pressure variable; a quantity not available 
in DPIV experiments~\cite{Noca97,Noca99}. We do not analyze 
such expressions in this work, however. In the context of IB method, 
Bergmann and co-workers~\cite{Bergmann11,Bergmann16} have 
used Flavio's moving control volume force expressions (involving both 
velocity and pressure) to compute hydrodynamic forces. They  
observed spurious force oscillations with a moving control volume approach~\cite{Bergmann11}. 
We show that by manipulating time derivatives in the original expressions, one 
can eliminate such spurious oscillations. We also present strategies to 
mitigate jumps in velocity derivatives in an AMR framework. Such jumps 
arise when the Cartesian grid hierarchy is regridded and velocity in the new grid hierarchy 
is reconstructed from the old hierarchy~\cite{Griffith07}. Lai and Peskin~\cite{Lai00} used stationary
control volume analysis on a uniform grid to compute steady state hydrodynamic forces on a stationary 
cylinder for Reynolds numbers between $100-200$ using the IB method. They did not 
consider time derivative terms in their analysis and temporal hydrodynamic force profiles 
were reported at steady state. In this work we include time derivative terms for 
finite Reynolds number flows (but not for steady Stokes flow) in our (moving) control volume 
analysis. 

If point-wise traction values are desired for IB-like methods, there are several 
recommendations proposed in the literature for smoothing them. Here we 
list a few of them. Verma et al.~\cite{Verma17} have recommended 
using a ``lifted'' surface: a surface two grid cell distance away from the 
actual interface to avoid choppy velocity gradients. This recommendation, 
based upon their empirical tests using a Brinkman penalization method, can 
change depending on the smoothness of the problem and the discrete delta 
function used in IB methods. Goza et al.~\cite{Goza16} obtain smooth point-wise force measurements
by using a force filtering post-processing step that penalizes inaccurate high frequency stress components.
Martins et al.~\cite{Martins17} enforce a continuity 
constraint in the velocity interpolation stencil to reconstruct a
second-order velocity field at IB surface points that is discretely divergence-free. 
They also impose a normal gradient constraint in the pressure 
interpolation stencil to reconstruct a second-order accurate pressure field 
at IB points. They have successfully eliminated spurious force 
oscillations using \emph{constrained} least-squares stencils. The idea 
of \emph{unconstrained} moving least-squares velocity interpolation and force  
spreading in the context of direct forcing IB method was first proposed by 
Vanella and Balaras~\cite{Vanella09}. Lee et al.~\cite{Lee11} attribute sources 
of spurious force oscillations to spatial pressure discontinuities across the
fluid-solid interface and temporal velocity discontinuities for moving bodies.
They recommend using fine grid resolutions to alleviate spurious oscillations.
\REVIEW{Their analyses and tests~\cite{Lee11} show that grid spacing has 
a more pronounced effect on spurious force oscillations than computational time 
step size.}

For a range of test cases varying from free-swimming at high Reynolds 
number to steady Stokes flow, we show that both LM and moving CV methods 
are in excellent agreement and are able to suppress spurious force 
oscillations. They do not require any additional treatment such as least-squares 
stencil (velocity and pressure) interpolation or force filtering beyond simple 
integration of force balance laws. For moderate to high Reynolds number 
test cases, we use a direct forcing IB method to estimate Lagrange multipliers, 
and for Stokes flow we use a fully constrained IB method 
to compute Lagrange multipliers exactly. 

\section{Equations of motion}

\subsection{Immersed boundary method}
The immersed boundary (IB) formulation uses an Eulerian description for the 
momentum equation and divergence-free condition for both 
the fluid and the structure. A Lagrangian description is employed for the structural 
position and forces. Let $\x = (x_1, \ldots, x_d) \in \Omega$ denote fixed 
Cartesian coordinates, in which $\Omega \subset \mathbb{R}^d$ is the fixed 
domain occupied by the entire fluid-structure system in $d$ spatial dimensions.
Let $\s = (s_1, \ldots s_d) \in U$ denote the fixed material coordinate system 
attached to the structure, in which $U \subset \mathbb{R}^d$ is the Lagrangian curvilinear 
coordinate domain. The position of the immersed structure occupying a volumetric 
region $\Vbt \subset \Omega$ at time $t$ is denoted by $\X (\s,t)$. We consider 
only neutrally buoyant bodies to simplify the implementation; this assumption implies 
that the fluid and structure share the same uniform mass density $\rho$. The deviatoric 
stress tensor of fluid, characterized by dynamic viscosity $\mu$, is extended inside 
the structure to make the momentum equation of both media appear similar. 
The combined equations of motion for fluid-structure system are~\cite{Peskin02}

\begin{align}
\rho\left(\D{\u}{t}(\x,t) + \u(\x,t) \cdot \grad \u(\x,t) \right) &= -\grad p(\x,t) + \mu \lap \u(\x,t) + \f(\x,t), \label{eqn_momentum}\\
  \div \u(\x,t) &= 0, \label{eqn_continuity} \\
\f(\x,t)  &= \int_{U} \F(\s,t) \, \delta(\x - \X(\s,t)) \, \Ds, \label{eqn_F_f} \\
  \U(\s,t) &= \int_{\Omega} \u(\x,t) \, \delta(\x - \X(\s,t)) \, \Dx, \label{eqn_u_interpolation} \\
   \D{\X}{t} (\s,t) &= \U(\s,t). \label{eqn_body_motion} 
\end{align}

Eqs.~\eqref{eqn_momentum} and \eqref{eqn_continuity} are the incompressible 
Navier-Stokes equations written in Eulerian form, in which $\u(\x,t)$ is the velocity, 
$p(\x,t)$ is the pressure, and $\f(\x,t)$ is the Eulerian force density, which is non-zero 
only in the structure region. Interactions between Lagrangian and Eulerian quantities 
in Eqs.~\eqref{eqn_F_f} and~\eqref{eqn_u_interpolation} are mediated by integral 
equations with Dirac delta function kernels, in which the $d$-dimensional delta function 
is $\delta(\x) = \Pi_{i=1}^{d}\delta(x_i)$. Eq.~\eqref{eqn_F_f} converts the Lagrangian 
force density $\F(\s,t)$ into an equivalent Eulerian density $\f(\x,t)$. In the IB literature,
the discretized version of this operation is called \emph{force spreading}.
Using short-hand notation, we denote force spreading operation by $\f = \cS[\X] \, \F$, 
in which $\cS[\X]$ is the \emph{force-spreading operator}.
Eq.~\eqref{eqn_u_interpolation} determines the physical velocity of each Lagrangian 
material point from the Eulerian velocity field, so that the immersed structure moves 
according to the local value of the velocity field $\u(\x,t)$ (Eq.~\eqref{eqn_body_motion}).
This \emph{velocity interpolation} operation is expressed as $\D{\X}{t} = \U = \cJ[\X] \, \u$, 
in which $\cJ[\X]$ is the \emph{velocity-interpolation operator}. It can be shown that
if $\cS$ and $\cJ$ are taken to be adjoint operators, i.e. $\cS = \cJ^{*}$, then 
Lagrangian-Eulerian coupling conserves energy~\cite{Peskin02}.

\subsection{Discrete equations of motion}

We employ a staggered grid discretization for the momentum and continuity equations
(see Fig.~\ref{fig_CV_discrete}). More specifically, Eulerian velocity and force 
variables are defined at face centers while the pressure variable is defined at cell centers. 
Second-order finite difference stencils are used to spatially discretize the Eulerian equations  
on locally refined grids~\cite{Bhalla13,Griffith07}.  The spatial discretization of various operators
are denoted with $h$ subscripts. To discretize equations in time, 
we take $\dt$ as the time step size, and $n$ as the time step number. 
We use the direct forcing method of Bhalla et al.~\cite{Bhalla13}  
for moderate to high Reynolds number cases and the fully constrained IB method
of Kallemov et al.~\cite{Kallemov16} for Stokes flow cases. The two methods differ in how 
the Lagrangian force density or Lagrange multipliers are computed. The time 
integrators are also different for the two methods. Here we briefly describe the 
discretized equations for both methods. We refer 
readers to~\cite{Bhalla13,Kallemov16,Usabiaga17} for more details. 
 
For the direct forcing method, the time stepping scheme reads as~\cite{Bhalla13}
\begin{align}
	\rho \left(\frac{\tilde{\u}^{n+ 1} - \u^n}{\dt} + [\u \cdot \grad_h \u]^{(n+\half)}\right) &= -\grad_h p^{n+\half} + \frac{\mu}{2} \V{\grad}^2_h\left( \tilde{\u}^{n+ 1} + \u^{n}\right), \label{eqn_lm_momentum} \\
	\grad_h \cdot \tilde{\u}^{n+1} &=  \V 0, \label{eqn_lm_continuity} \\
	\F^{n+\half} &= \frac{\rho}{\Delta t} \left(\U_b^{n+1} - \cJ_h[\X^{n+\half}] \tilde{\u}^{n+1} \right) \label{eqn_lm_correction}, \\
	 \rho \left(\frac{\u^{n+ 1} - \tilde{\u}^{n+1}}{\dt} \right) &= \cS_h[\X^{n+\half}] \F^{n+\half}.
\end{align}
Succinctly, we first solve for a velocity field $\tilde{\u}^{n+1}$ and a pressure field $p^{n+\half}$ as a coupled 
system by solving Eqs.~\eqref{eqn_lm_momentum} and  \eqref{eqn_lm_continuity} 
simultaneously. The velocity $\tilde{\u}^{n+1}$, which is correct in the fluid region but not 
in the structure region $\Vbt$, is then corrected by estimating the Lagrange multiplier 
$\F^{n+\half}$ via Eq.~\eqref{eqn_lm_correction}. Here $\U_b^{n+1}$ is the desired 
rigid body velocity of the Lagrangian nodes, and $\X^{n+\half}$ is the midstep 
estimate of Lagrangian node position. Finally, the Lagrange multiplier $\F^{n+\half}$ is 
spread on the background grid to correct the momentum in the structure region to 
$\u^{n+ 1}$. We use Adams-Bashforth to approximate the midstep value of  nonlinear 
convection term $\u \cdot \grad \u$ via
\begin{equation}
	[\u \cdot \grad_h \u]^{(n+\half)} = \frac{3}{2} \u^{n} \cdot \grad_h \u^{n} - \frac{1}{2} \u^{n-1} \cdot \grad_h \u^{n-1}.
\end{equation}

For the fully constrained method, we \emph{simultaneously} solve for the updated 
Eulerian velocity $\u^{n+1}$ and pressure $p^{n+1}$ at time $t^{n+1}$ along 
with the Lagrange multiplier $\F^{n+1}$. The time stepping scheme reads as~\cite{Kallemov16,Usabiaga17}
\begin{align}
	-\grad_h p^{n+1} + \mu \V{\grad}^2_h \u^{n+1} + \cS_h[\X^n]  \F^{n+1} & = \V 0, \label{eqn_stokes_momentum} \\
	\grad_h \cdot \u^{n+1} &=  \V 0, \label{eqn_stokes_continuity} \\
	\cJ_h[\X^n] \u^{n+1} &= \U_b^{n+1}. \label{eqn_stokes_constraint} 
\end{align}
The coupled system of Eqs.~\eqref{eqn_stokes_momentum}-\eqref{eqn_stokes_constraint} 
is solved simultaneously using a preconditioned FGMRES~\cite{Saad93} solver. For 
both methods, we update the Lagrangian node positions $\X^{n+1}$ using rigid body translation 
and rotation. We use Peskin's 4-point regularized delta functions for the $\cS$ and
$\cJ$ operators in all our numerical experiments, unless stated otherwise.  

Our finite Reynolds number fluid solver has support for adaptive mesh refinement, and some cases presented in 
Sec.~\ref{sec_results} make use of multiple grid levels (also known as a \emph{grid hierarchy}).
A grid with $\ell$ refinement levels with grid spacing $\dx_0$, $\dy_0$, and $\dz_0$ on the coarsest grid level
has minimum grid spacing $\dx_\textrm{min} = \dx_0/\nref^{\ell-1}$, 
$\dy_\textrm{min} = \dy_0/\nref^{\ell-1}$, and $\dz_\textrm{min} = \dz_0/\nref^{\ell-1}$ on the finest grid level.
Here, $\nref \in \mathbb{N}$ is the refinement ratio.
In the present work, the refinement ratio in taken to be the same in each direction, although
this is not a limitation of the numerical method. The immersed structure is always placed
on the finest grid level. \REVIEW{For all of the cases considered in this work a constant time step size 
$\dt = \min(\dt^{\ell})$ is chosen, in which the time step size $\dt^{\ell}$ on grid level $\ell$ satisfies 
the convective CFL condition 
$\dt^{\ell} \le C \min \left(\frac{\dx}{\|u_x\|_{\infty}}, \frac{\dy}{\|u_y\|_{\infty}}, \frac{\dz}{\|u_z\|_{\infty}}\right)^{\ell}$. 
In this work, the convective CFL number is set to $C = 0.3$ unless otherwise stated.}


\section{Hydrodynamic force and torque}
\subsection{Moving control volume method}

\subsubsection{Hydrodynamic force}
Letting $\T = \mu\left(\grad \u + \grad \u^T\right)$ denote the
viscous stress tensor, the net hydrodynamic force is defined to be
the force of the \emph{fluid on the body}:

\begin{equation}
\cF(t) = -\oint_{\Sbt} \ndot \left[-p \I + \T\right] \dS,  \label{eq_force}
\end{equation}
in which the integral is taken over the surface of the body $\Sbt = 
\partial \Vbt$, and $\n$ is the unit outward normal to the surface. In practice, 
evaluating Eq.~\eqref{eq_force} is inconvenient in numerical experiments 
because it is often difficult to obtain accurate surface velocity gradients and 
pressure values. Moreover, evaluating Eq.~\eqref{eq_force} also requires computational 
geometry to obtain surface normals and area. Instead we use a control volume 
approach to compute $\cF(t)$ which avoids these requirements. 

\begin{figure}[H]
  \centering
    \includegraphics[scale = 0.3]{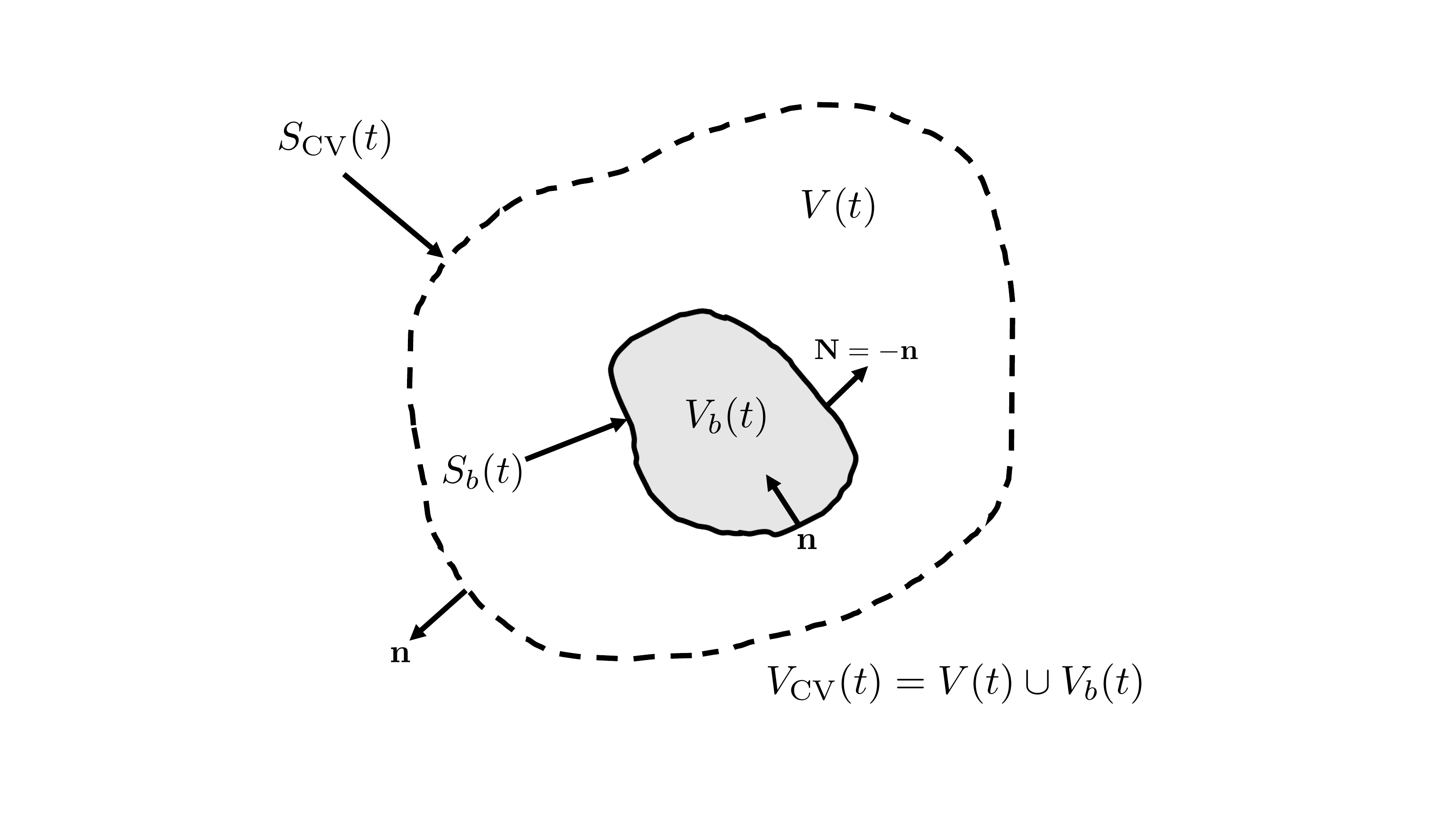}
  \caption{Sketch of the immersed structure (solid line)
  surrounded by an arbitrary control volume (dashed line)}.
  \label{fig_CV_continuous}
\end{figure}

Let an arbitrary (possibly time dependent) domain $\Vcvt$ completely surround
$\Vbt$, i.e. $\Vbt \subset \Vcvt$, as shown in Fig.~\ref{fig_CV_continuous}. By 
considering the change in momentum within the control volume $\Vcvt$ and the 
net momentum flux at its surface $\Scvt = \partial \Vcvt$, a general expression 
for the hydrodynamic force on the body can be obtained

\begin{equation}
\cF(t) = -\d{}{t} \int_{\Vt} \rho \u \dV + \oint_{\Scvt} \ndot \left[-p \I - (\u-\uS)\rho \u + \T \right] \dS  
- \oint_{\Sbt} \ndot (\u-\uS)\rho \u \dS, \label{eq_NocaForce}        
\end{equation}
in which $\Vt = \Vcvt \setminus \Vbt$ is the volume outside the immersed 
body but inside the CV, and $\uS$ is the velocity of the surface over which the (surface) integral 
is evaluated. In general $\uS \neq \u$, and $\uS$ can be arbitrarily chosen so 
that the moving CV always encloses the immersed body. The unit normal $\n$ in 
Eq.~\eqref{eq_NocaForce} points \emph{outward} on $\Scvt$, and \emph{into}
$\Sbt$. The above force equation was first derived by Noca~\cite{Noca97} and 
has been used in various experimental~\cite{Noca99, Unal97, vanOudheusden07, Jardin09} 
and numerical~\cite{Bergmann11, Bergmann16, Shen09, Sallstrom14} studies.
For Cartesian grid based methods, the CV can be chosen as a 
simple rectangular domain, for which the unit normals on $\Scvt$ are aligned with the
Cartesian axes. Finally, the integral over $\Sbt$ vanishes for many applications 
where no-slip ($\uS = \u$) boundary condition can be chosen for $\uS$.
Henceforth, we will analyze cases for no-slip boundary conditions.

When Eq.~\eqref{eq_NocaForce} is discretized, the first term requires
a discrete approximation of the integral at two separate time instances
or at two different locations of the moving control volume. Bergmann et al.
~\cite{Bergmann11} observed spurious 
force oscillations as a result of this time derivative term. We show that by 
manipulating the first term using the Reynolds transport theorem, 
an expression for hydrodynamic force can be obtained that does not require
contributions from control volumes at two different spatial locations. The 
modified equation reads as

\begin{equation}
\label{eq_OurForce}
\cF(t) =  -\int_{\Vcvt} \rho \D{\u}{t} \dV + \d{}{t} \int_{\Vbt} \rho \u \dV + \oint_{\Scvt} \ndot \left[-p \I - \u\rho \u + \T \right] \dS.
\end{equation}

A detailed derivation of Eq.~\eqref{eq_OurForce} is provided in~\ref{app_force}.
Note that although Eqs.~\eqref{eq_NocaForce} and \eqref{eq_OurForce} 
are equivalent formulas to obtain the hydrodynamic force on an immersed body,
their physical interpretations are different. In Eq.~\eqref{eq_NocaForce}, the control 
volume is \emph{moving} with some prescribed velocity, usually chosen to follow 
the structure to ensure it is contained within the CV at all times. Hence, 
one must keep track of the velocity $\uS$ of the control surface $\Scvt$. 
Moreover since the time derivative appears outside of the integral over 
$V(t)$, it requires a discrete approximation of momentum on two 
time-lagged control volumes. On the other hand in Eq.~\eqref{eq_OurForce}, 
all integrals over the control volume are evaluated at a single time instance 
and therefore a discrete evaluation on two separate CVs is never 
needed. Discretely, the control volume is \emph{placed} at a new location 
at time step (${n+1}$) and no information from its previous location at time step
$n$ is used. Hence, $\uS$ never appears in the calculation. This has the numerical
benefit of suppressing force oscillations, which will be shown in 
Sec.~\ref{sec_results}.

\subsubsection{Hydrodynamic torque}

The net hydrodynamic torque on an immersed body is defined to be the net moment of 
hydrodynamic force exerted by \emph{fluid on the body} about a given reference point:

\begin{equation}
\cM(t) =  -\oint_{\Sbt} \rcross \left(\ndot \left[-p \I + \T\right]\right) \dS, \label{eq_torque}
\end{equation}
in which $\r = \x - \x_0$. The torque is computed with respect to
some reference point $\x_0$, which can be fixed at a location or
move with time (e.g. the center of mass of a swimmer). Following Noca's 
derivation for the force expression Eq.~\eqref{eq_NocaForce}, 
one can measure the change in angular momentum within a moving control 
volume to obtain an expression for torque which reads as

\begin{equation}
\label{eq_NocaTorque}
\cM(t) =  -\d{}{t} \int_{\Vt} (\rcross \rho \u) \dV
	+ \oint_{\Scvt} [\rcross (- p \; \n +  \ndot \T) - \ndot (\u-\uS) (\rcross \rho \u)  ] \dS 
	- \oint_{\Sbt} \ndot (\u - \uS)(\rcross \rho \u) \dS. 
\end{equation}

\noindent Note that the torque expression in Eq.~\eqref{eq_NocaTorque} is slightly 
different from Eq.~15b given in Bergmann and  Iolla~\cite{Bergmann11}. 
Once again by applying the Reynolds transport theorem to the first term on 
left-hand side of Eq.~\eqref{eq_NocaTorque}, we obtain a torque 
expression involving a control volume contribution at a single spatial location
and without any $\uS$ terms

\begin{equation}
\label{eq_OurTorque}
\cM(t) = 
-\int_{\Vcvt} \rho \rcross \D{\u}{t}\dV
+ \d{}{t} \int_{\Vbt} \rho (\rcross \u) \dV 
+ \oint_{\Scvt}  [\rcross (- p \; \n + \ndot \T) - (\ndot \u)\rho (\rcross \u) \;]  \dS.
\end{equation}
For the derivation of torque expressions~\eqref{eq_NocaTorque} and~\eqref{eq_OurTorque},
see~\ref{app_torque}.


\subsubsection{Numerical integration}
\begin{figure}[H]
  \centering
    \includegraphics[scale = 0.35]{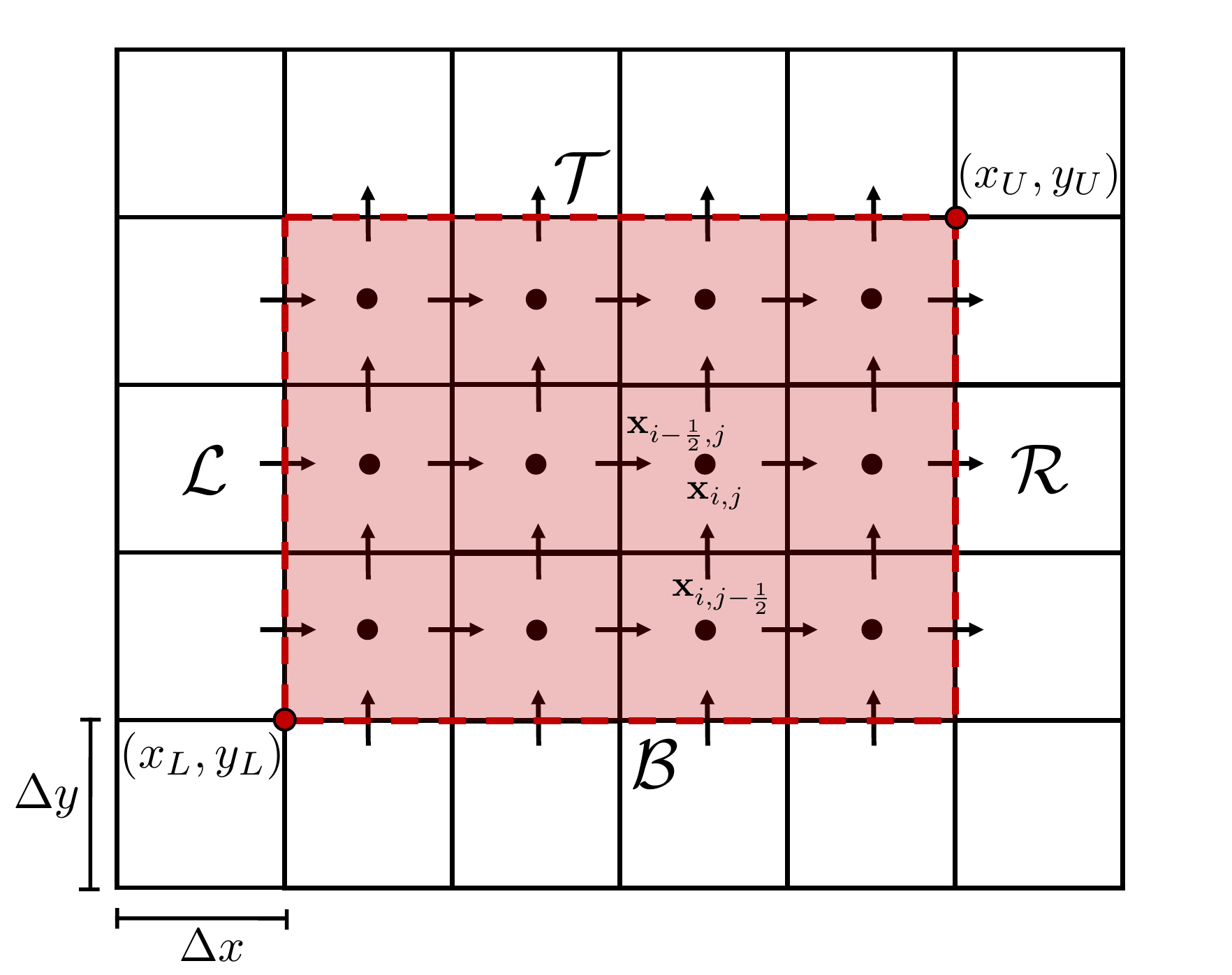}
  \caption{A staggered grid spatial discretization.
  The $x$-velocity component $u$ is solved for at locations depicted with horizontal arrows.
  The $y$-velocity component $v$ is solved for at locations depicted with vertical arrows.
  The pressure $p$ is solved for at locations depicted with a solid black dots.
  The corners $(\xL,\yL)$ and $(\xU,\yU)$ define the control volume $\Vcvt$, which is shaded in red.
  The dashed line represents $\Scvt = \partial \Vcvt$ = $\cL \cup \cT \cup \cR \cup \cB$.}
  \label{fig_CV_discrete}
\end{figure}

We use Riemann summation to evaluate the various integrals in 
Eqs.~\eqref{eq_OurForce} and~\eqref{eq_OurTorque} over a moving rectangular control volume. 
Fig.~\ref{fig_CV_discrete} shows the rectangular control volume marked by its 
lower and upper coordinates $(\xL,\yL)$ and $(\xU,\yU)$, respectively. 

The arbitrary surface velocity $\uS$ is chosen such that the moving CV is forced to
align with Cartesian grid faces.
This greatly simplifies the evaluation of various terms inside the force and torque integrals. 
The linear and angular momentum integrals over $\Vbt$ are 
evaluated in the Lagrangian frame, whereas the rest are computed in the Eulerian frame. 
The details of these computations are given in~\ref{app_numerical_integration}. 

\subsection{Lagrange multiplier method}
The hydrodynamic force and torque on an immersed body can also be computed
in an extrinsic manner from force/torque balance laws. Specifically, for an immersed
body occupying volume $\Vbt$, the force and torque balance laws read as
\begin{align}
 \d{}{t} \int_{\Vbt} \rho \; \U \dV &= \cF + \int_{\Vbt} \F \dV, \label{eq_extrinsic_force} \\
  \d{}{t} \int_{\Vbt} \rho \; \Rcross \U \dV &= \cM + \int_{\Vbt} \Rcross \F \dV, \label{eq_extrinsic_torque}
 \end{align}
in which $\R = \X - \X_0$ is the radius vector from reference point $\X_0$ to Lagrangian node
position $\X$, and $\F$ is the Lagrange multiplier imposing rigidity constraint as defined in 
Eq.~\eqref{eqn_F_f}. Note that for the IB method, the net hydrodynamic force $\cF$ and torque 
$\cM$ on the body can be readily evaluated in the Lagrangian frame as a part of the solution 
process without computing any extra terms. We have used this approach in a previous 
work~\cite{Bhalla13}. We will compare results obtained from both moving CV and LM method 
in Sec.~\ref{sec_results}.

\subsection{Equivalence of the two methods}

Although Eqs.~\eqref{eq_OurForce} and \eqref{eq_extrinsic_force} look different, 
we now show that the use of (regularized) delta functions in the IB method 
make them equivalent expressions. To prove this, consider a single body in the domain 
$\Omega$, occupying a region of space $\Vbt \subset \Omega$.  Let $\F(\s,t)$
denote the Lagrange multiplier field defined on Lagrangian nodes, and $\f(\x,t)$
denote its Eulerian counterpart. Here, $\s \in U$ and $ U \subset \mathbb{R}^d$.
With continuous and Peskin's discrete delta functions, the following 
identity holds 

\begin{align}
\int_{\Omega} \f(\x,t)  \, \Dx &= \int_{\Omega} \left[\int_{U \subset \Omega} \F(\s,t) \, \delta(\x - \X(\s,t)) \, \Ds \right]    \, \Dx   \nonumber\\
					  &= \int_{U \subset \Omega} \left[\int_{\Omega} \F(\s,t) \, \delta(\x - \X(\s,t)) \, \Dx \right]    \, \Ds   \nonumber\\
                                            &= \int_{U \subset \Omega}  \F(\s,t) \left[ \int_{\Omega} \delta(\x - \X(\s,t)) \, \Dx \right] \Ds  \nonumber \\
                                            &= \int_{U} \F(\s,t) \, \Ds \label{eq_lag_eul_equivalence}.
\end{align}  
The above expression is the equivalence of Lagrangian and Eulerian force
densities and is a well known result~\cite{Peskin02, Griffith07}.  Starting with the momentum 
Eq.~\eqref{eqn_momentum} in Eulerian form

\begin{equation}
\rho \D{\u}{t} + \rho \div (\u \u)
= \div (- p \I +  \T) + \f, \label{eqn_conservative_momentum_form}
\end{equation}

\noindent and integrating it over $\Vcvt$, we obtain

\begin{equation}
\int_{\Vcvt} \rho \D{\u}{t} \dV
+ \int_{\Vcvt} \rho \div (\u \u) \dV
= \int_{\Vcvt} \div \left[-p\I + \T \right] \dV 
+ \int_{\Vcvt} \f \dV. \nonumber
\end{equation}

\noindent Applying the divergence theorem to terms with $\div$, and using 
Eq.~\eqref{eq_lag_eul_equivalence} we obtain

\begin{equation}
\label{eq_sumLag_CV}
\int_{\Vcvt} \f \dV = \int_{\Vbt} \F \dV =   \int_{\Vcvt} \rho \D{\u}{t} \dV
- \int_{\Scvt} \ndot \left[-p\I - \u \rho \u+ \T \right] \dS. 
\end{equation}

\noindent Using the above expression for $\int_{\Vbt} \F \dV$ with Eq.~\eqref{eq_extrinsic_force}
we get force expression~\eqref{eq_OurForce}. Similarly, one can prove the equivalence of
torque expressions~\eqref{eq_OurTorque} and \eqref{eq_extrinsic_torque} by noting that
$\rcross \div \S = \div (\rcross \S)$ for any symmetric tensor satisfying $\S = \S^T$.

When multiple bodies exist in the domain, however, there is a subtle difference in the LM
and moving CV expressions. The LM expressions in the form of Eqs.~\eqref{eq_extrinsic_force}
and~\eqref{eq_extrinsic_torque} are restricted to individual bodies; therefore in the presence
of multiple bodies, $\cF$ and $\cM$ can be computed separately for each body. 
For the moving CV method, care must be taken to restrict the control volume to 
a particular body, i.e. it should not enclose other bodies in its vicinity. Otherwise,  
$\cF$ and $\cM$ will contain contribution from multiple bodies. We explore this 
subtlety with the moving CV method by taking flow past two stationary cylinders in 
Sec.~\ref{sec_two_cylinders}, \REVIEW{and drafting-kissing-tumbling of two 
sedimenting cylinders in Sec.~\ref{sec_dkt_cylinders}}. For an AMR framework, the boundary of a control 
volume can span multiple refinement levels. When velocity is reconstructed from old hierarchy
to new, there is generally no guarantee that momentum is conserved because of the 
\emph{inter-level} velocity interpolation. This can lead to jumps in $\cF$ and $\cM$ because 
of the time derivative terms. Eqs.~\eqref{eq_extrinsic_force} 
and~\eqref{eq_extrinsic_torque} are restricted to $\Vbt$, which generally 
lies on the finest grid level. Therefore, one can expect force and torque calculations 
to be relatively insensitive to velocity reconstruction operations, which mostly 
requires \emph{intra-level} interpolation. We explore this issue with translating plate 
example in Sec.~\ref{sec_translating_plate}.


\section{Software implementation} \label{sec_software_implementation}
We use the IBAMR library~\cite{IBAMR-web-page} to implement the moving control volume 
method in this work for our numerical tests. IBAMR has built-in support for direct 
forcing and fully constrained IB methods, among other variants of the IB method.
IBAMR relies on SAMRAI \cite{HornungKohn02, samrai-web-page} for Cartesian grid 
management and the AMR framework. Solver support in IBAMR is provided by 
PETSc library~\cite{petsc-efficient, petsc-user-ref, petsc-web-page}.


\section{Results}
\label{sec_results}

\subsection{Flow past cylinder}
In this section we validate our moving control volume method and Lagrange 
multiplier method for computing hydrodynamic forces and torques on immersed
bodies.
\subsubsection{Stationary cylinder}
\label{sec_stationary_cylinder}

We first consider the flow past a stationary circular cylinder. The cylinder 
has diameter $D = 1$ and is placed in a flow with far-field velocity 
$\U_{\infty} = (U_{\infty}, V_{\infty})=(1,0)$. The computational domain 
is a rectangular channel taken to be of size $18 D \times 12 D$, with 
the center of the cylinder placed at $(x,y) = (0,0)$. The domain is discretized 
by a uniform Cartesian mesh of size $900 \times 600$. The aerodynamic 
drag coefficient $C_D = \cF \cdot \e_x/(\rho D \|\U_{\infty}\|^2/2)$ is calculated numerically 
in two different ways: via Eq.~\eqref{eq_OurForce} and via integrating Lagrange 
multipliers enforcing the rigidity constraint on the cylinder (Eq.~\eqref{eq_extrinsic_force}). The control volume
is taken to be $[-D, 1.5 D] \times [-D, D]$ and does not move from its initial location.
Note that in the case where both the body and the control volume are stationary,
Eqs.~\eqref{eq_NocaForce} and \eqref{eq_OurForce} equivalent and give the same numerical 
solution. The density is set to $\rho = 1$ and the Reynolds number of the flow is $\Re = \rho U_{\infty} D / \mu = 550$. 
This problem has been studied numerically by Bergmann and Iollo~\cite{Bergmann11} 
and by Ploumhans and Winckelmans~\cite{Ploumhans2000}. The temporal behavior 
of $C_D$ matches well with the previous studies~\cite{Bergmann11,Ploumhans2000}.

\begin{figure}[H]
  \centering
    \includegraphics[scale = 0.3]{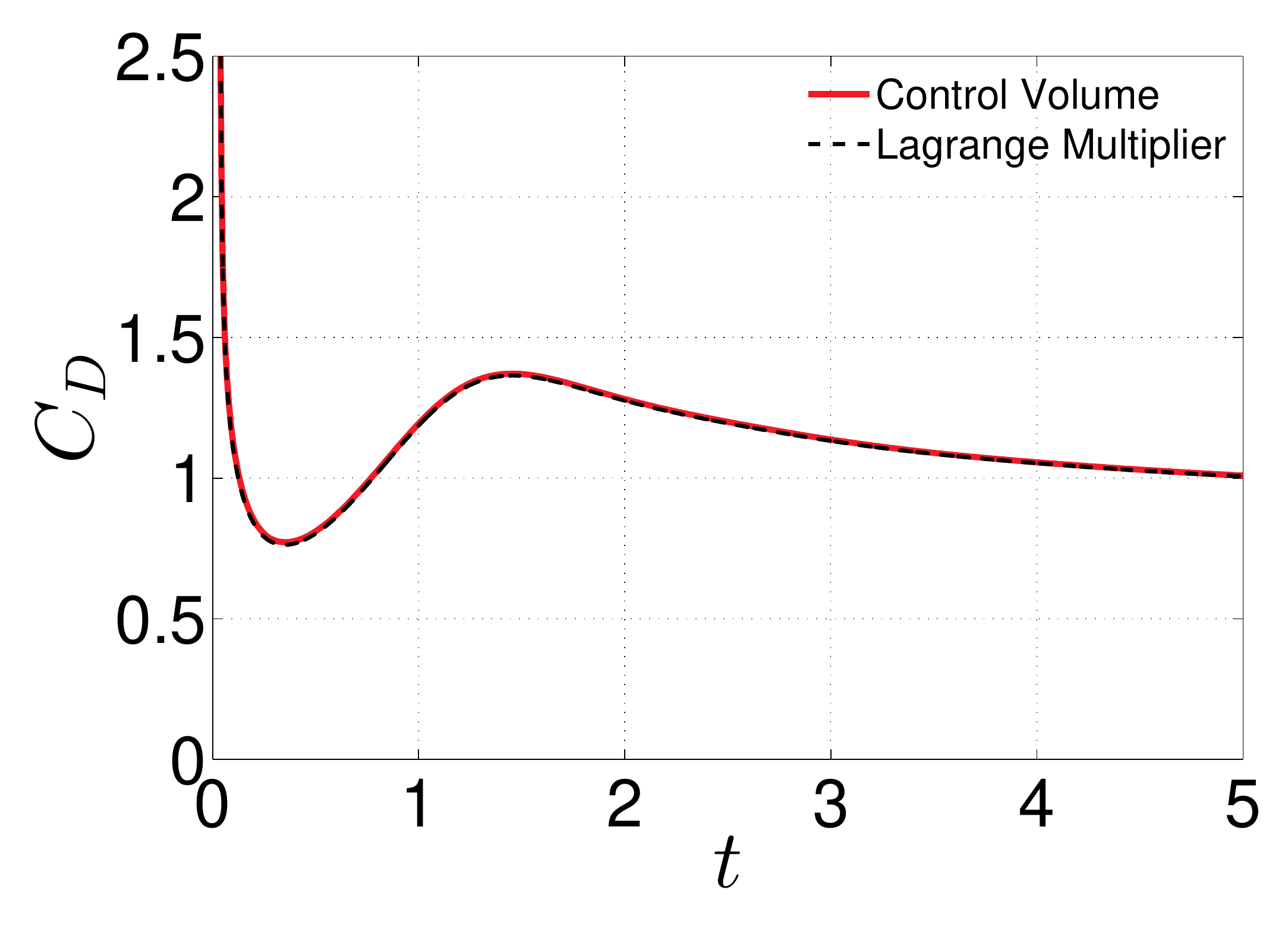}
  \caption{Comparison of the drag coefficient for flow past a
    cylinder at $\Re = 550$ measured in two different ways.
    (---, red): Control volume using Eq.~\eqref{eq_OurForce};
     (\texttt{---}, black): Lagrange multiplier using Eq.~\eqref{eq_extrinsic_force}.}
  \label{fig_stationary_cylinder}
\end{figure}

\subsubsection{Translating cylinder}

Next, we consider the case of a circular cylinder translating with 
prescribed motion. The parameters used in this case are identical 
to those of the flow past a stationary cylinder case, except now 
$\U_{\infty} = (U_{\infty}, V_{\infty})=(0,0)$. The cylinder is dragged 
with speed $\Ub = (U_\text{b}, V_\text{b}) = (-1,0)$. The aerodynamic 
drag coefficient is computed as $C_D = \cF \cdot \e_x/(\rho D \|\Ub\|^2/2)$. 
The control volume is initially set to $[-D, 1.5D] \times [-D, D]$ and 
translates to the left every few time steps to ensure that it always 
contains the cylinder. The density is set to $\rho = 1$ and the 
Reynolds number is $\Re = \rho U_\text{b} D / \mu = 550$.
Periodic boundary conditions are used on all faces of the computational domain.
This problem was also studied numerically by Bergmann and Iollo \cite{Bergmann11}.

\begin{figure}[H]
  \centering
    \includegraphics[scale = 0.3]{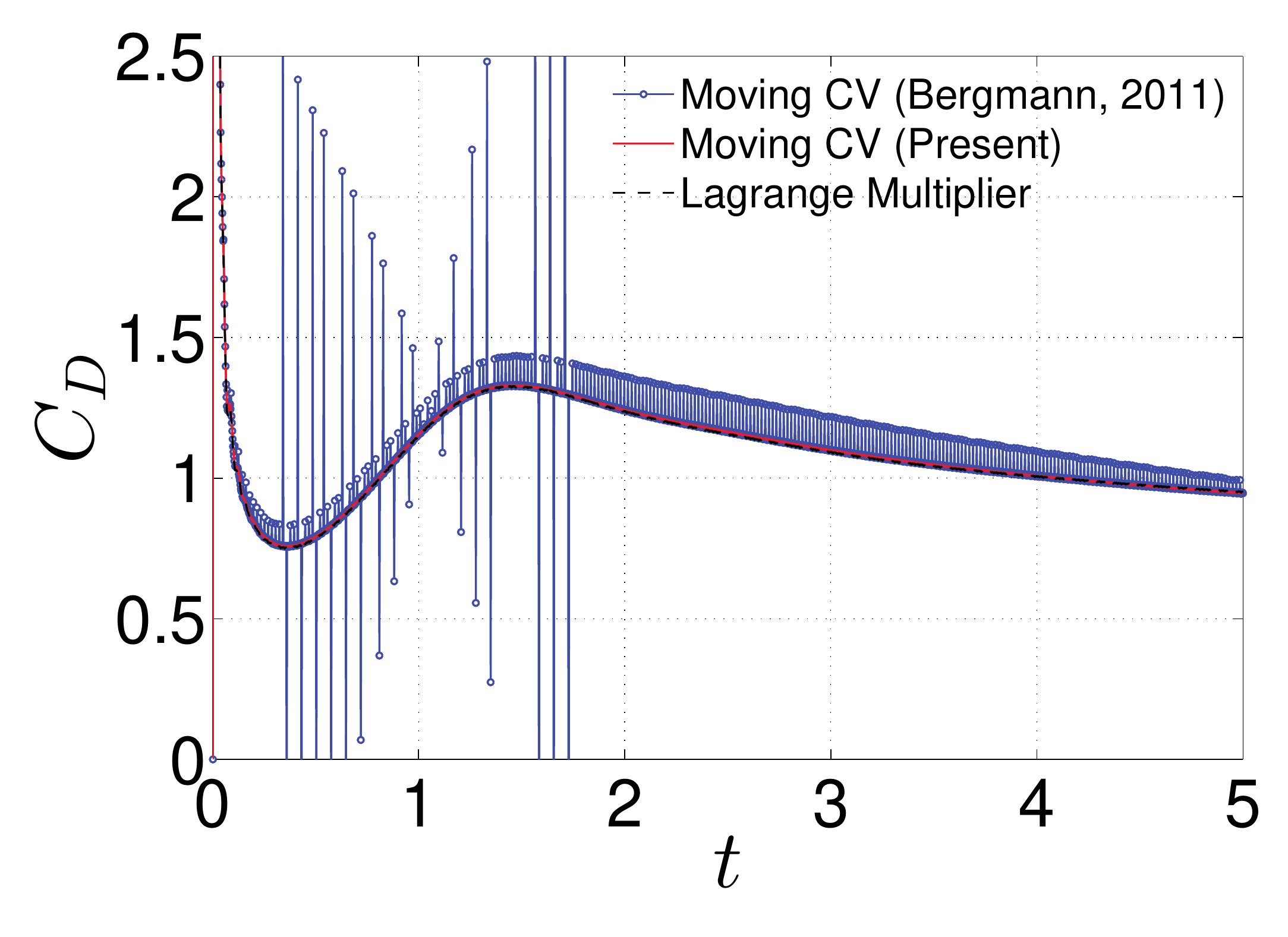}
  \caption{Comparison of the drag coefficient a translating 
    cylinder at $\Re = 550$ measured in three different ways.
    $\bullet$: Control volume using Eq.~\eqref{eq_NocaForce};
    ---: Control volume using Eq.~\eqref{eq_OurForce};
     \texttt{---}: Lagrange multiplier. Note that LM and present moving
     CV curves are on top of each other.}
  \label{fig_translating_cylinder}
\end{figure}

Fig.~\ref{fig_translating_cylinder} shows the time evolution of $C_D$ 
calculated in three different ways, via Eqs.~\eqref{eq_NocaForce} and~\eqref{eq_OurForce}
, and by integrating Lagrange multipliers. The present control volume method
matches well with the Lagrange multiplier approach. However, spurious 
jumps in drag coefficient are seen for the original control volume method
outlined by Noca. These spurious oscillations are also present in the computation
done by Bergmann and Iollo~\cite{Bergmann11}. \REVIEW{We remark that the oscillations 
seen by Bergmann and Iollo~\cite{Bergmann11} are quantitatively different than the 
ones presented here in Fig.~\ref{fig_translating_cylinder} for comparison, although both
use Eq.~\eqref{eq_NocaForce} to compute hydrodynamic force. 
This can be attributed to differences in the numerical method used to impose 
constraint: we use a constraint-based immersed boundary method whereas 
Bergmann and Iollo use a Brinkman penalization method.}

\subsubsection {Two stationary cylinders} \label{sec_two_cylinders}


To study the effect of control volume size in the presence of multiple bodies, 
we consider the case of two stationary circular cylinders each
with $D = 1$ and placed in a flow with far-field velocity 
$\U_{\infty} = (U_{\infty}, V_{\infty})=(1,0)$. The bottom and top cylinders
are centered about position $(x,y) = (0,-2D)$ and $(x,y) = (0,2D)$, respectively within
a rectangular channel of size $18 D \times 12 D$. The domain is discretized 
by a uniform Cartesian mesh of size $900 \times 600$ and the Reynolds number 
is $\Re = \rho U_{\infty} D / \mu = 550$ for each cylinder. The density of the fluid
is set to $\rho = 1$.

Four different (but symmetric) control volume configurations are considered: 
\begin{enumerate}
	\item Two disjoint CVs located at $[-D,D]\times[-3D,-D]$ and $[-D,D]\times[D,3D]$.
	\item Two CVs located at  $[-D,D]\times[-3D,D]$ and $[-D,D]\times[-D,3D]$ that 
		slightly overlap, but do not intersect the other cylinder.
	\item Two CVs located at $[-D,D]\times[-3D,2D]$ and $[-D,D]\times[-2D,3D]$,
		    where each CV holds one full cylinder and half of the second cylinder.
         \item Two CVs located at $[-D,D]\times[-3D,2.7D]$ and $[-D,D]\times[-2.7D,3D]$,
         	 where each CV contains both cylinders. 
\end{enumerate}

\begin{figure}[H]
  \centering
  \subfigure[No overlapping]{
    \includegraphics[scale = 0.2]{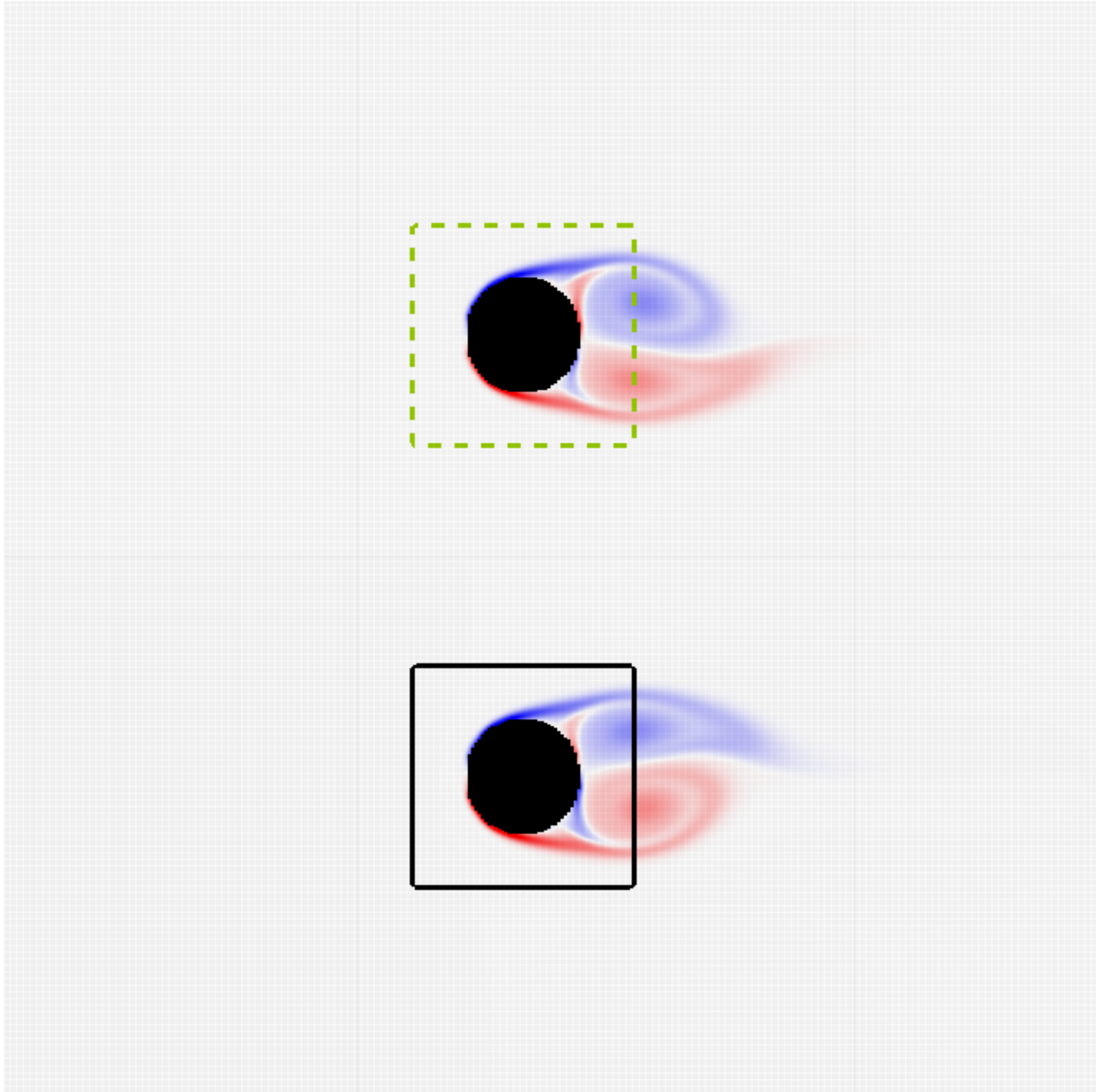}
    \label{fig_no_overlap_viz}
  }
   \subfigure[Slightly overlapping]{
    \includegraphics[scale = 0.2]{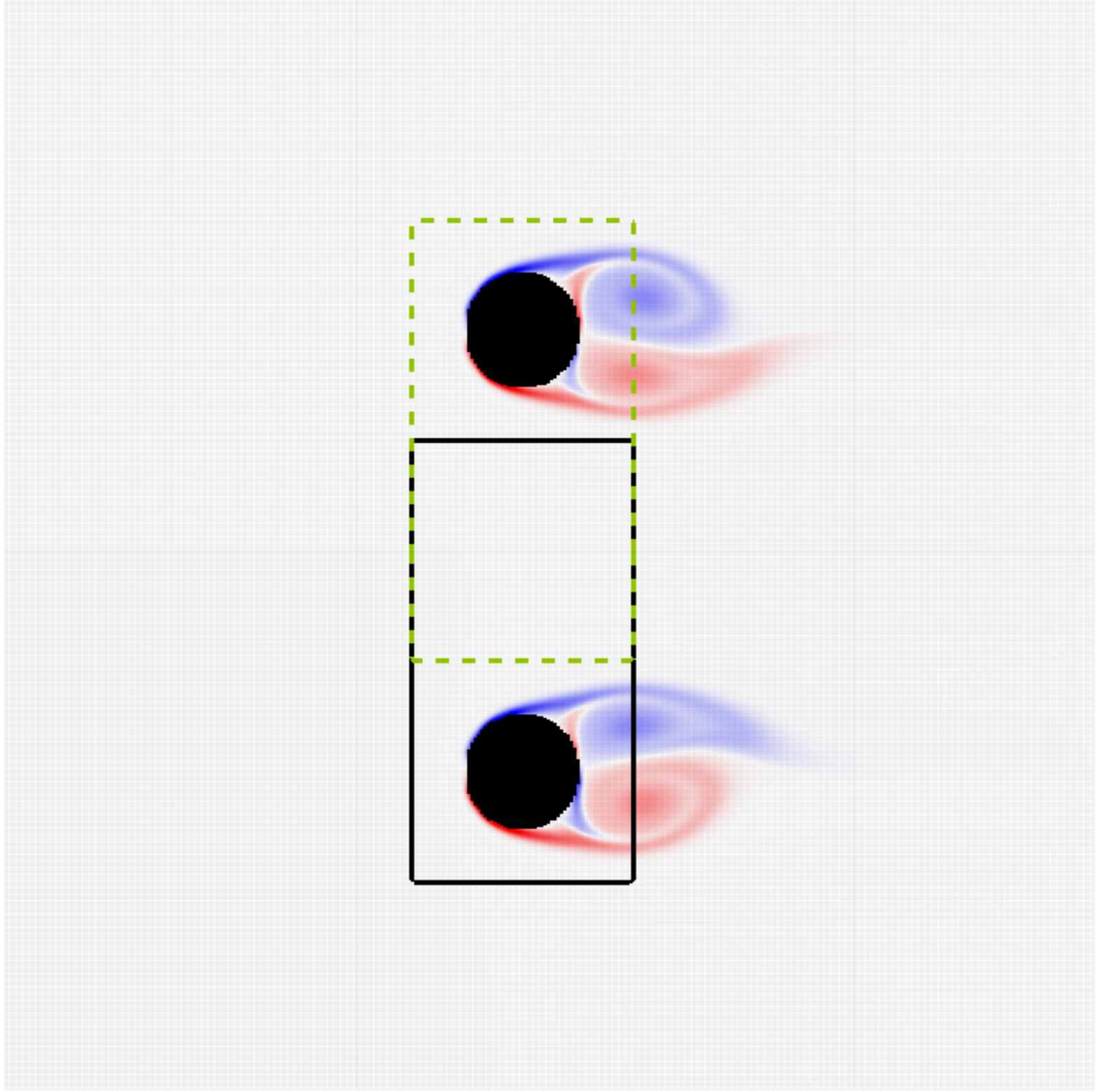}
    \label{fig_small_overlap_viz}
  }
    \subfigure[CVs contain $1.5$ cylinders]{
    \includegraphics[scale = 0.2]{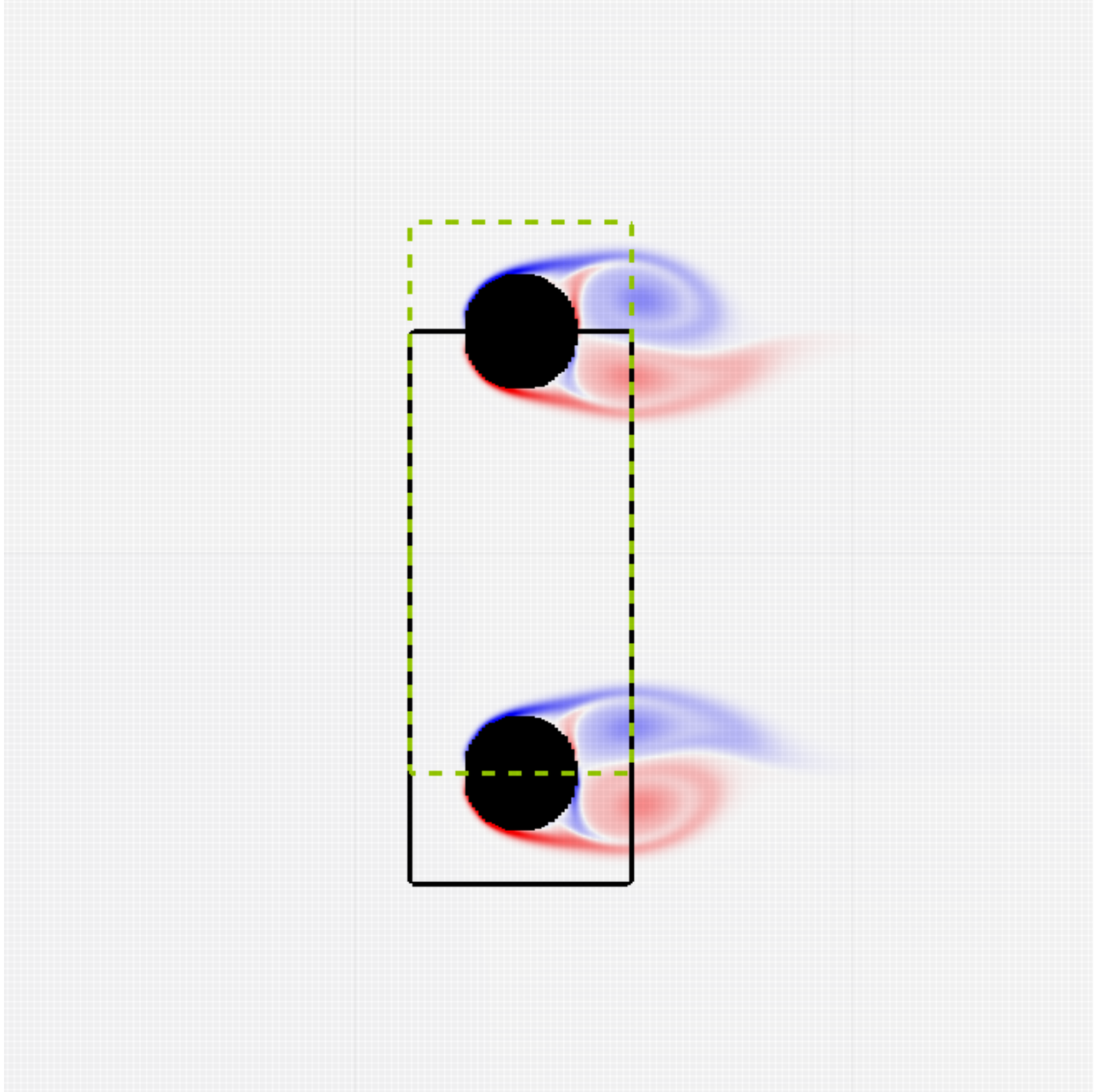}
    \label{fig_med_overlap_viz}
  }
   \subfigure[CVs contain both cylinders]{
    \includegraphics[scale = 0.2]{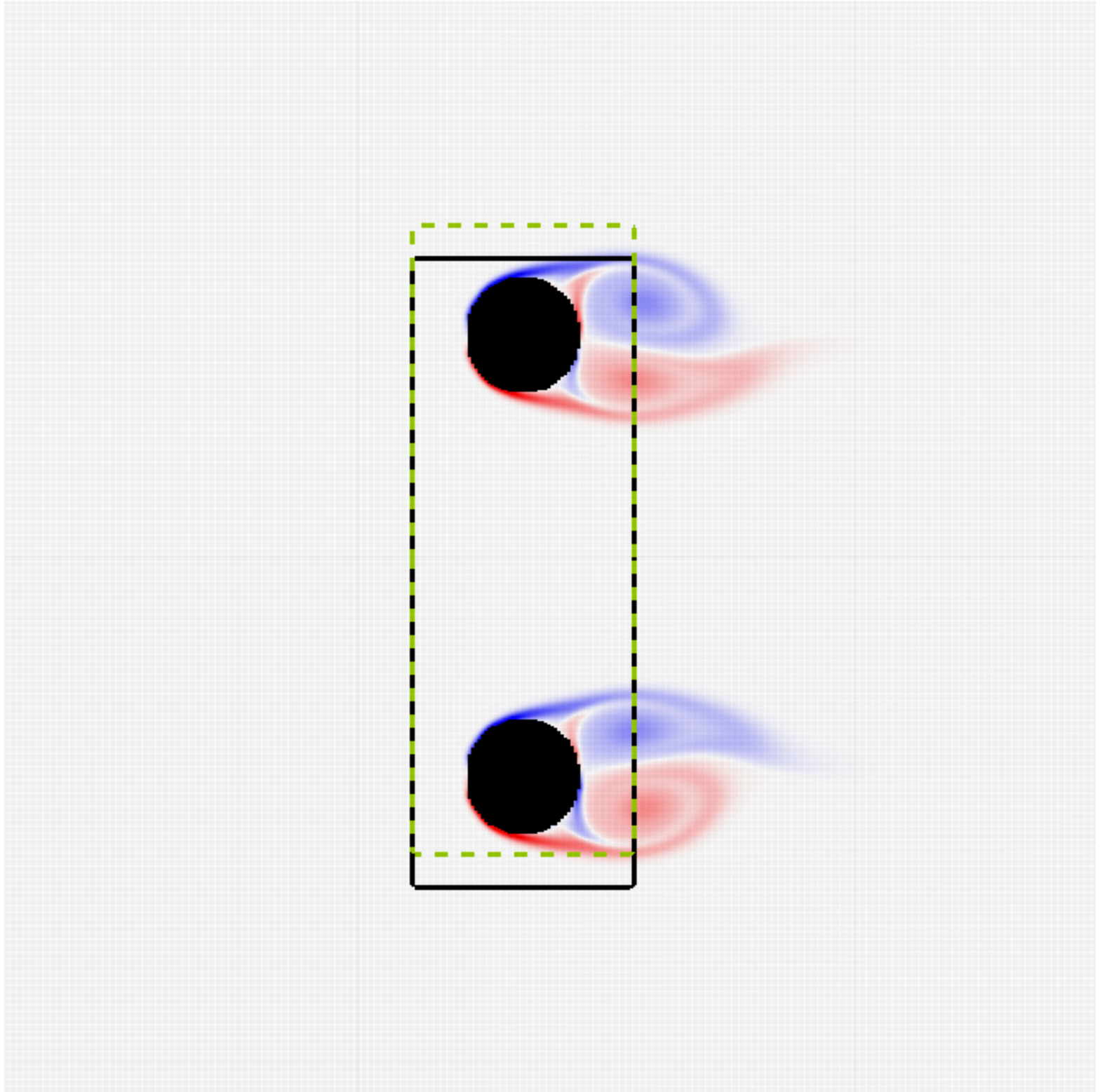}
    \label{fig_big_overlap_viz}
  }
  \caption{Vorticity generated by a circular cylinder at $t=5$ for $\Re = 550$: 
  \subref{fig_no_overlap_viz} two disjoint control volumes;
  \subref{fig_small_overlap_viz} two slightly overlapping control volumes;
  \subref{fig_med_overlap_viz} each control volume contains one and a half cylinders;
  \subref{fig_big_overlap_viz} each control volume contains both cylinders.
   All figures are plotted for vorticity between $-20$ and $20$.
   }
  \label{fig_cv_configurations}
\end{figure}

\begin{figure}[H]
  \centering
    \subfigure[No overlapping]{
    \includegraphics[scale = 0.25]{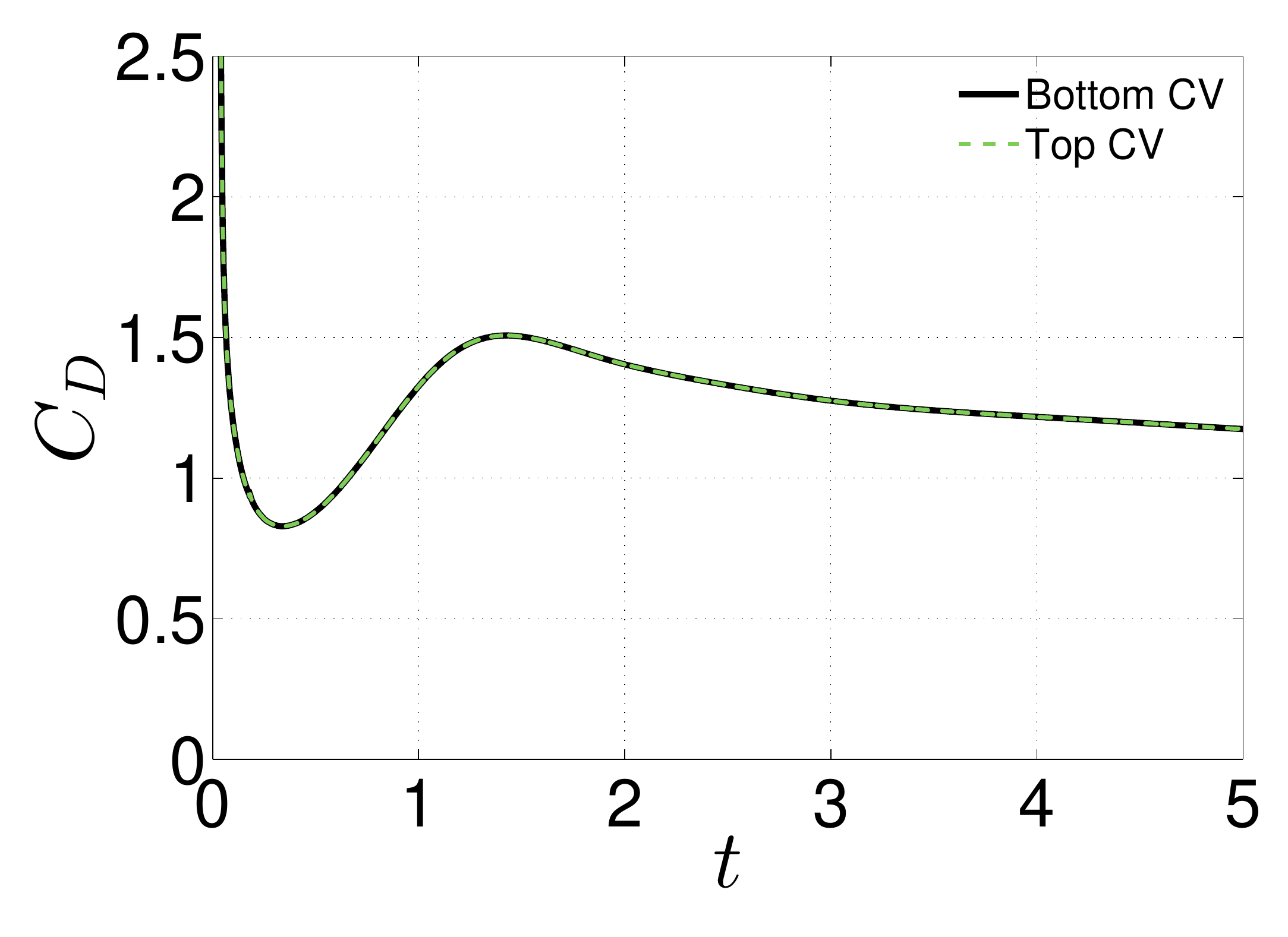}
    \label{fig_no_overlap_CD}
  }
  \subfigure[Slightly overlapping]{
    \includegraphics[scale = 0.25]{Figures/no_overlap_cylinder.pdf}
    \label{fig_small_overlap_CD}
    }
  \subfigure[CVs contain $1.5$ cylinders]{
    \includegraphics[scale = 0.25]{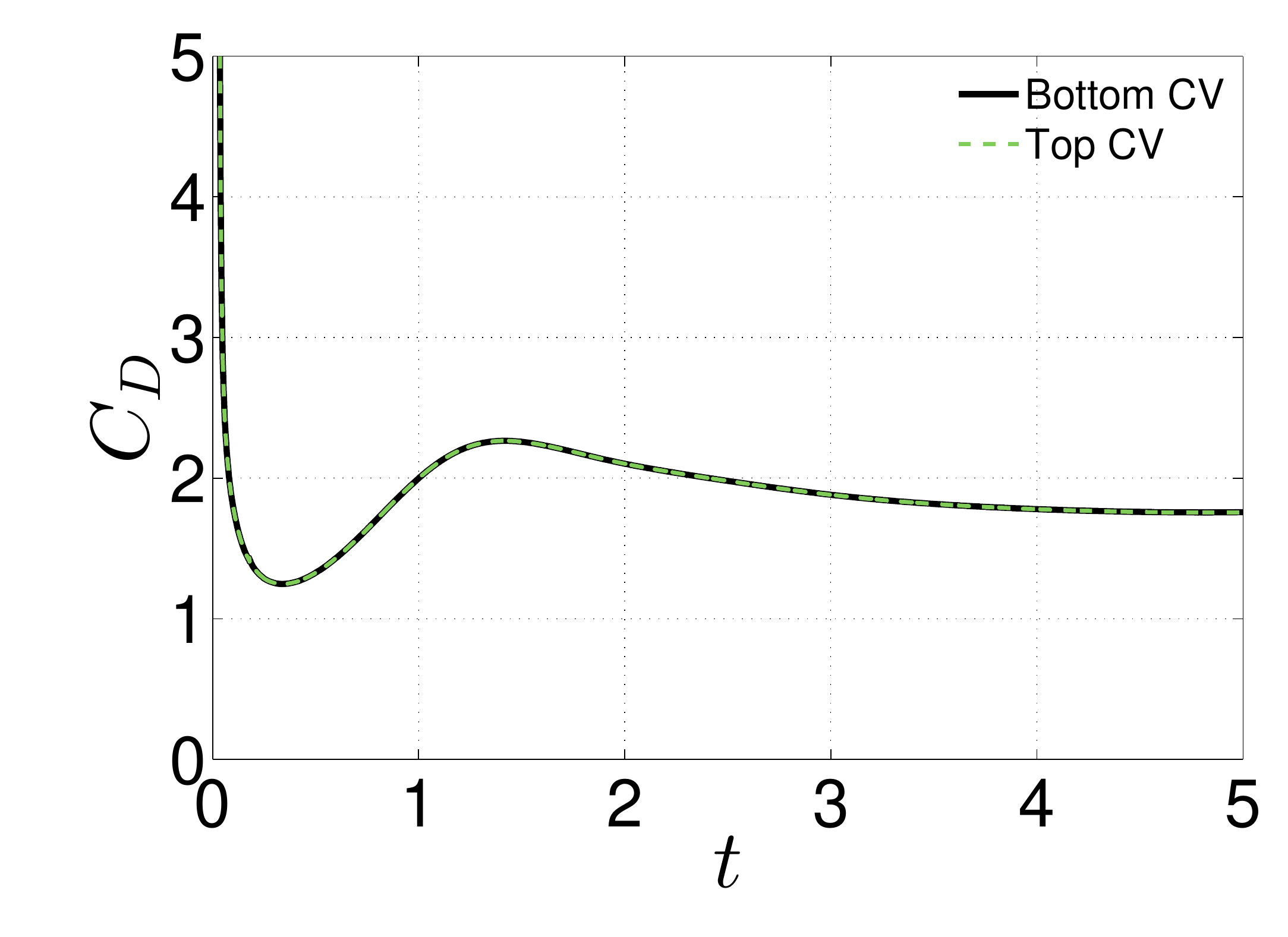}
    \label{fig_med_overlap_CD}
  }
  \subfigure[CVs contain both cylinders]{
    \includegraphics[scale = 0.25]{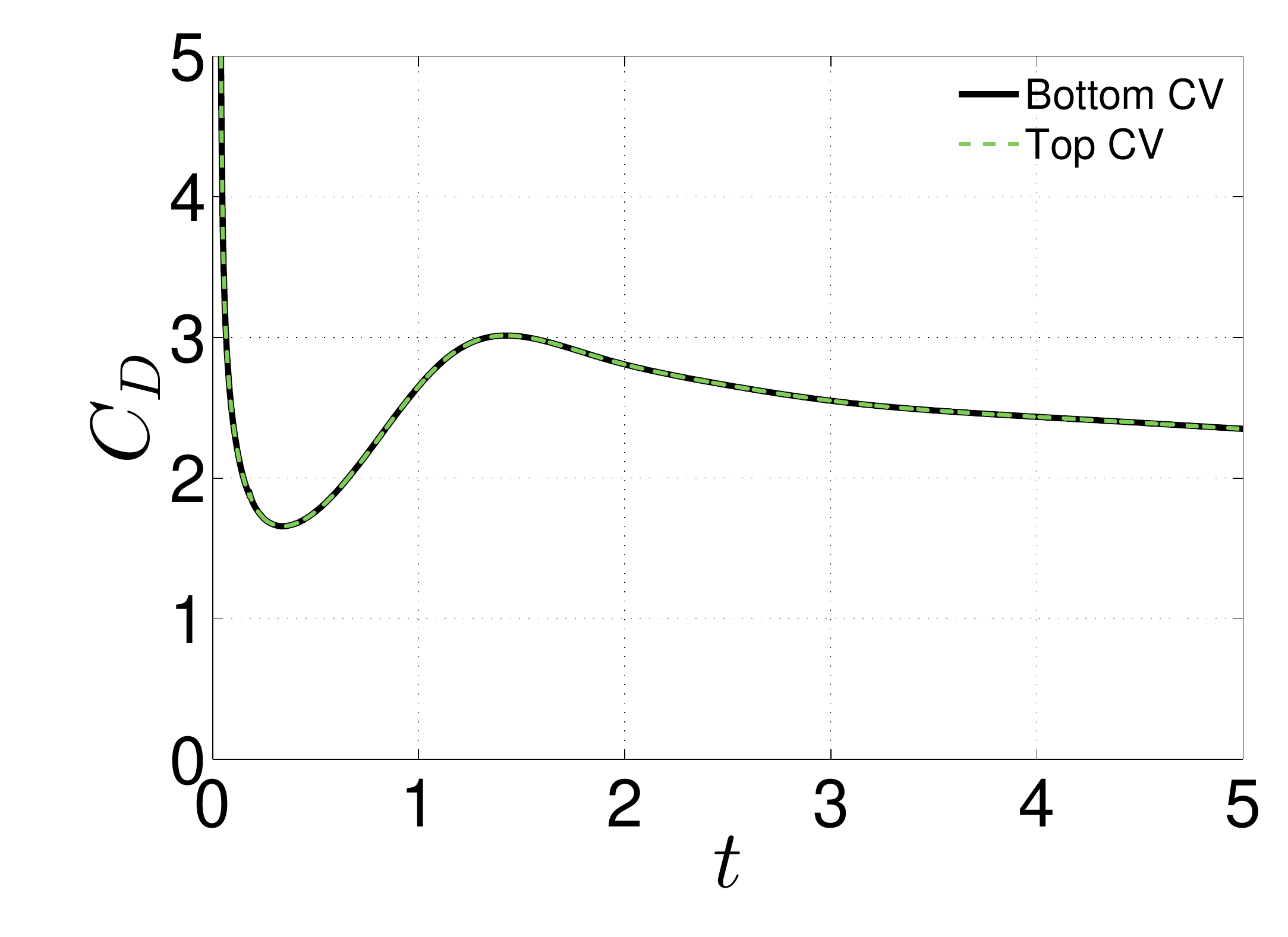}
    \label{fig_big_overlap_CD}
  }
  \caption{Temporal evolution of drag coefficient for $\Re = 550$:
  \subref{fig_no_overlap_CD} two disjoint control volumes; 
   \subref{fig_small_overlap_CD} two slightly overlapping control volumes; 
   \subref{fig_med_overlap_CD} one and a half cylinders;
   \subref{fig_big_overlap_CD} two full cylinders. 
   ---: Bottom CV;
   \texttt{---} (green): Top CV.}
  \label{fig_CD_cylinder_cv_configurations}
\end{figure}

Fig.~\ref{fig_cv_configurations} shows flow visualizations for the four different 
CV configurations. Fig.~\ref{fig_CD_cylinder_cv_configurations} shows the 
drag coefficient over time for each of the four configurations.
In the case where the CVs do not overlap (Figs.~\ref{fig_no_overlap_viz} and
~\ref{fig_no_overlap_CD}), or when they overlap but do not enclose multiple 
bodies either partially or fully (Figs.~\ref{fig_small_overlap_viz} and ~\ref{fig_small_overlap_CD})
the drag coefficients calculated for the top and bottom cylinders are close to the drag 
coefficient calculated for a single cylinder considered in Section~\ref{sec_stationary_cylinder}.

The effect of resizing the control volumes to partially or fully contain other objects
is seen in the hydrodynamic drag force measurement.  
In the case where each CV contains one and a half cylinders 
(Figs.~\ref{fig_med_overlap_viz} and \ref{fig_med_overlap_CD}), 
the drag coefficient deviates significantly from the drag coefficient calculated for 
a single cylinder. Rather, the computed force is the drag on a combined full and 
a half cylinder contained within the CV, which is approximately 
$1.5$ times the drag on a single cylinder. In the case
where each CV contains both cylinders (Figs.~\ref{fig_big_overlap_viz} and 
\ref{fig_big_overlap_CD}), the measured $C_D$ in each control
volume is approximately twice the $C_D$ measured on a single cylinder.
Therefore, in presence of multiple bodies in the domain, care must be 
taken to restrict the CV to an individual body.

\subsection{Oscillating cylinder}
In this section we consider cases from previous studies, some of which have reported 
spurious force oscillations in hydrodynamic drag and lift forces with IB methods. 
We do not observe such spurious force oscillations using LM and CV methods within
an immersed boundary framework.  

\subsubsection{In-line oscillation}
We consider an in-line oscillation of a circular cylinder in a quiescent
flow as done in D\"{u}tsch et al.~\cite{Dutsch98} and Lee et al.~\cite{Lee11}. 
The cylinder has a diameter $D = 1$ and is placed in a domain of size
$[-16D,16D] \times [-8D, 8D]$ with zero velocity prescribed on all boundaries
of the domain. The initial center of mass of the cylinder is placed at $ (0,0)$ and its 
velocity is set to $\Ub = \left(-U_0 \cos(2 \pi f t), 0\right)$, in which
$f$ is the frequency of oscillation. The Reynolds number of the flow is $\Re = \rho U_0 D/\mu = 100$,
and the Keulegan-Carpenter number is $\text{KC} = U_0/(fD) = 5$.
The time period of oscillation of the cylinder is given by $T = 1/f$.
The density of the fluid is set to be $\rho = 1$.
These parameters are chosen to match those reported in~\cite{Lee11, Dutsch98}.

Three levels of mesh refinement
are used, with $\nref = 4$ between each level. The cylinder is embedded in the finest mesh level
at all time instances. At the coarsest level, three different mesh sizes are used:
$50 \times 25$, $100 \times 50$, and $200 \times 100$, which corresponds to finest grid
spacings of $\dx_\textrm{min} = \dy_\textrm{min} = 0.04D$,
$\dx_\textrm{min} = \dy_\textrm{min} = 0.02D$, and $\dx_\textrm{min} = \dy_\textrm{min} = 0.01D$,
respectively. 
The computational time step size is chosen to be $\dt = 0.005D/U_0$, matching that of Lee at al. ~\cite{Lee11}.
\REVIEW{This time step size satisfies the convective CFL condition with $C = 0.7$, which is found to be stable for
all the mesh sizes considered here.}
A stationary control volume is placed at $[-4D, 4D] \times [-2D, 2D]$ in order to contain the entire cylinder at all time instances. The drag coefficient is computed as $C_D = \cF \cdot \e_x/(0.5 \rho U_0^2D)$.

\begin{figure}[H]
  \centering
  \subfigure[Control Volume]{
    \includegraphics[scale = 0.3]{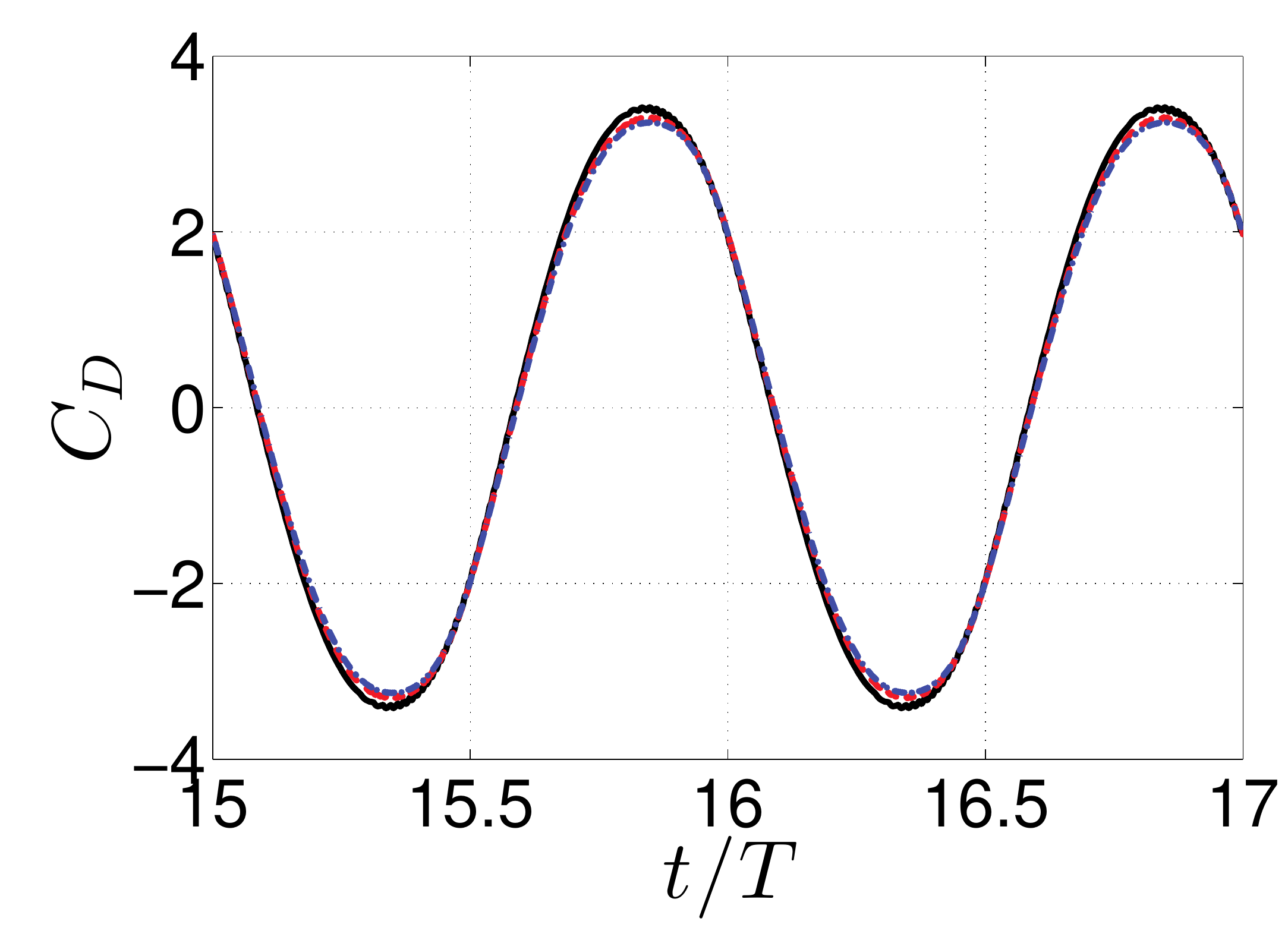}
    \label{fig_oscillating_cylinder_inline_CV}
  }
   \subfigure[Lagrange Multiplier]{
    \includegraphics[scale = 0.3]{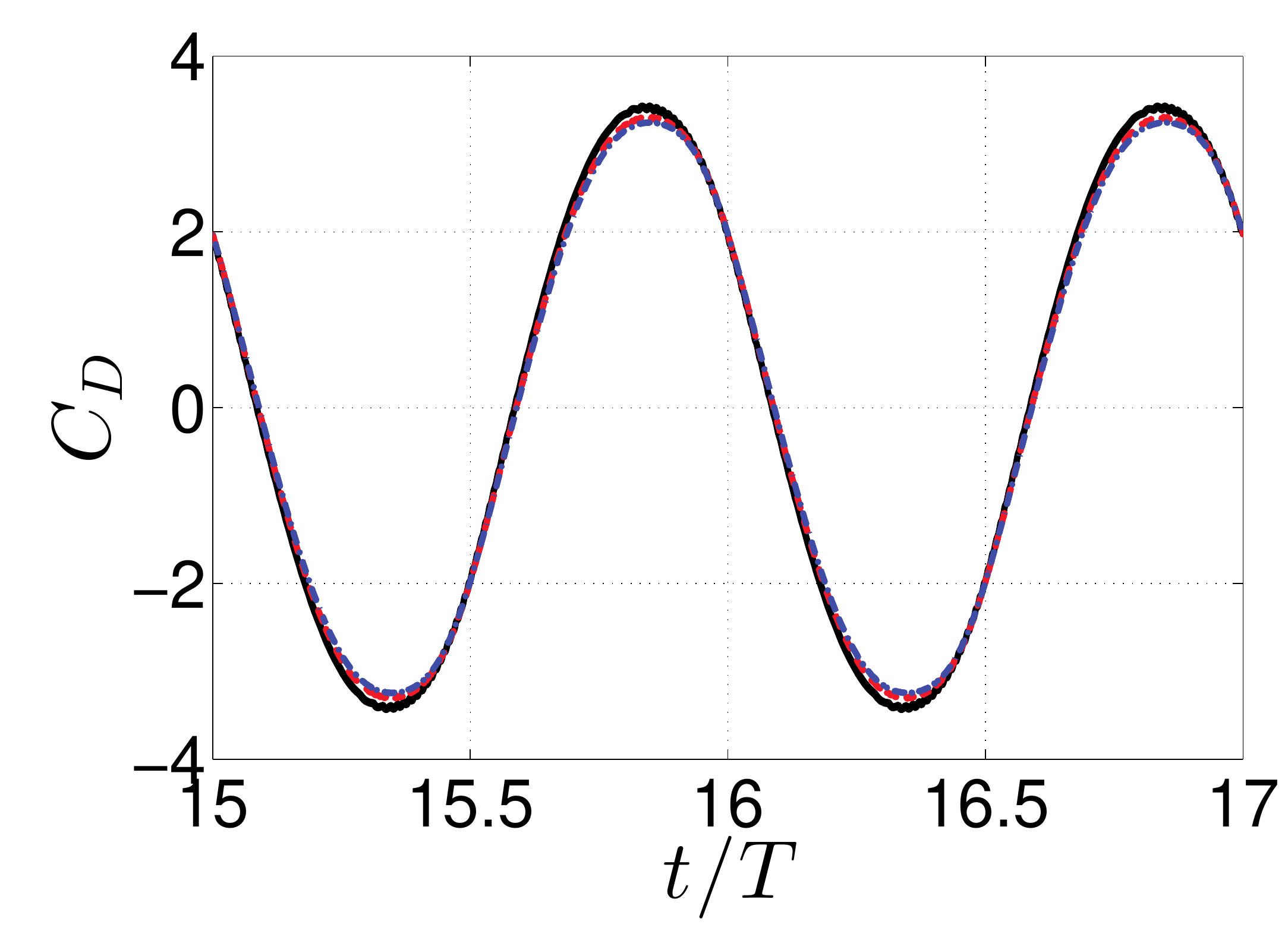}
    \label{fig_oscillating_cylinder_inline_Lag}
  }
  \caption{
   Temporal evolution of the drag coefficient of an oscillating cylinder at $\Re = 100$ measured by the 
   \subref{fig_oscillating_cylinder_inline_CV} present control volume approach and
    \subref{fig_oscillating_cylinder_inline_Lag} Lagrange multiplier approach.
    Here, (---, black): $\dx_\textrm{min} = \dy_\textrm{min} = 0.04D$;
    (\texttt{---}, red): $\dx_\textrm{min} = \dy_\textrm{min} = 0.02D$;
    (\texttt{-}$\cdot$\texttt{-}, blue): $\dx_\textrm{min} = \dy_\textrm{min} = 0.01D$.}
  \label{fig_oscillating_cylinder_inline}
\end{figure}

Fig.~\ref{fig_oscillating_cylinder_inline} shows the time evolution of drag coefficient for the
oscillating cylinder. Both CV and LM approaches yield identical results. Even our coarse resolution
results are in good agreement with the high resolution results of Lee et al.~\cite{Lee11}. Moreover, Lee et al. 
conducted this case within an immersed boundary framework and observed large spurious force 
oscillations (see Fig.~14 in \cite{Lee11}) at coarse grid resolutions. They compute drag on the immersed body 
by evaluating pressure and velocity gradients within the body region $\Vbt$. Only at fine grid resolutions,
where the spatial gradients are more accurate, were they able to suppress the spurious force oscillations. 
Since we do not require spatial gradients of velocity and pressure within the body region $\Vbt$ with our approach, we 
do not observe such spurious oscillations.

\subsubsection{Cross-flow oscillation}
Here we consider a cylinder oscillating in the transverse direction with an imposed
cross-flow. The cylinder has a diameter $D = 1$ and its initial center of mass
is placed at $(x,y) = (0,0.2D)$. The cylinder is placed in a domain of size $[-5D,27D]\times[-8D,8D]$
and oscillates in the transverse direction with velocity $\Ub = \left(0, -V_0 \sin(2 \pi f_e t) \right)$,
where $f_e$ is the frequency of oscillation. An axial free-stream velocity $U_{\infty}$ is set at
the inlet, top, and bottom faces of the computational domain. The transverse traction components
are set to zero on the top and bottom boundaries, and the axial and transverse tractions are set to zero
at the outflow boundary.

\begin{figure}[H]
  \centering
    \includegraphics[scale = 0.3]{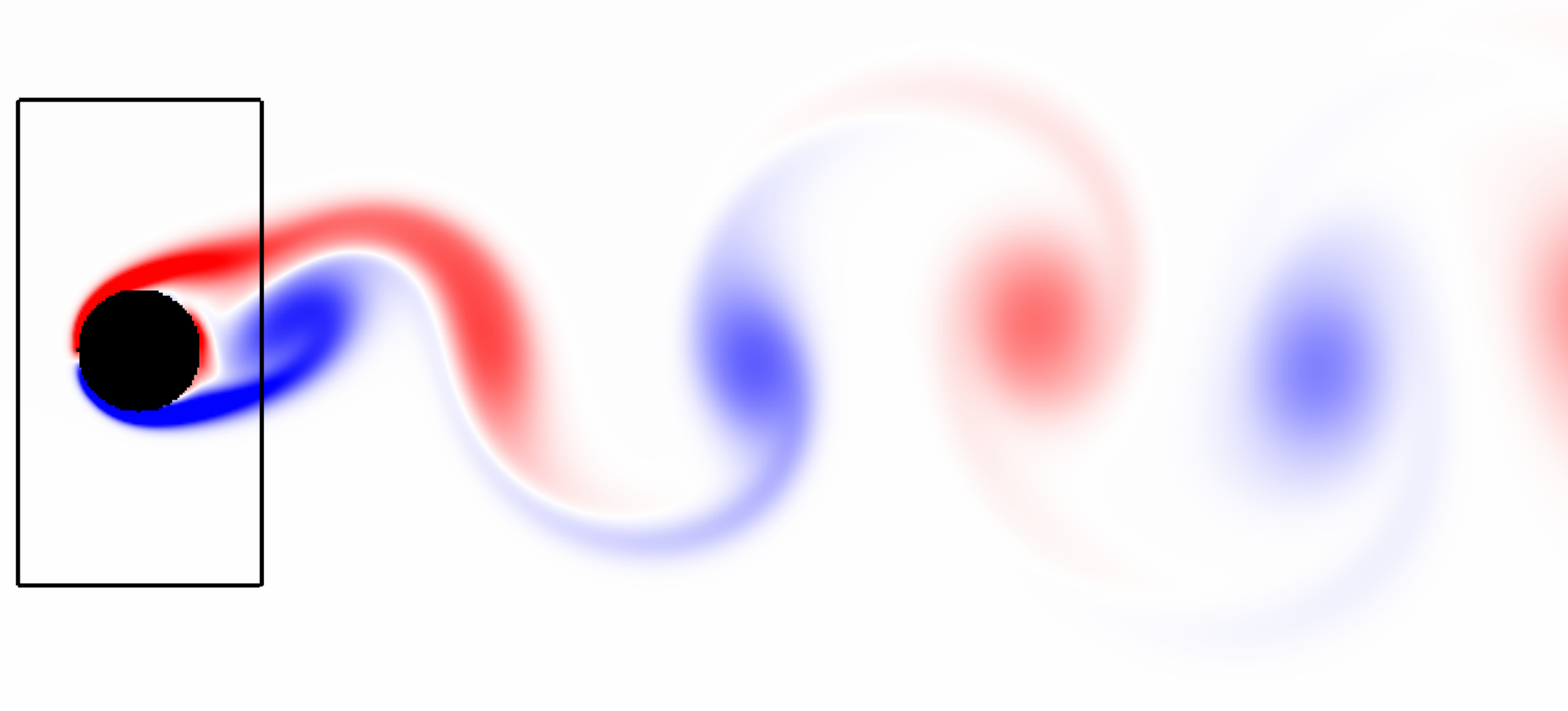}
  \caption{ Vorticity generated by an oscillating cylinder in a cross-flow at $t U_{\infty}/D = 100$ for $\Re = 185$.
  The plotted vorticity is between $-5$ and $5$.}
   \label{fig_oscillating_cylinder_crossflow_viz}
  \end{figure}

The density of the fluid is set to $\rho = 1$ and the Reynolds number based on 
the free-stream velocity is $\Re = \rho U_{\infty} D/\mu = 185$.
Letting $f_0 = 0.193 U_{\infty}/D$ be the natural shedding frequency for a stationary cylinder,
we set $f_e/f_0 = 1$. The maximum oscillation velocity of the cylinder is taken to be
$0.159 V_0/(f_e D) = 0.2$. This case is considered in Lee et al.~\cite{Lee11} and 
Guilmineau and Queutey~\cite{Guilmineau02}. A stationary CV is placed at $[-D,D]\times[-2D,2D]$, 
which contains the cylinder at all time instances (see Fig.~\ref{fig_oscillating_cylinder_crossflow_viz}).

\begin{figure}[H]
  \subfigure[Control Volume]{
    \includegraphics[scale = 0.3]{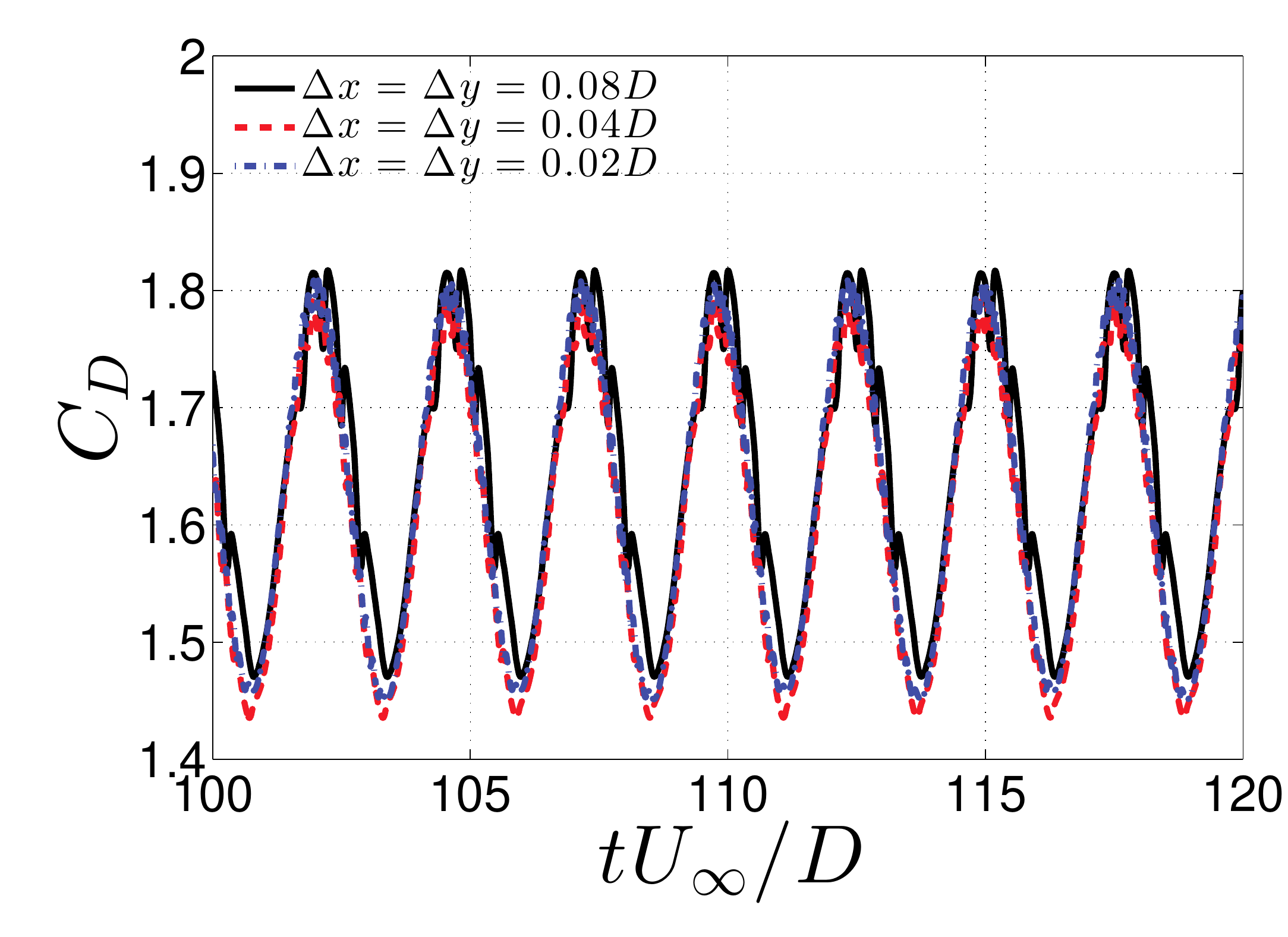}
    \label{fig_oscillating_cylinder_crossflow_CV}
  }
   \subfigure[Lagrange Multiplier]{
    \includegraphics[scale = 0.3]{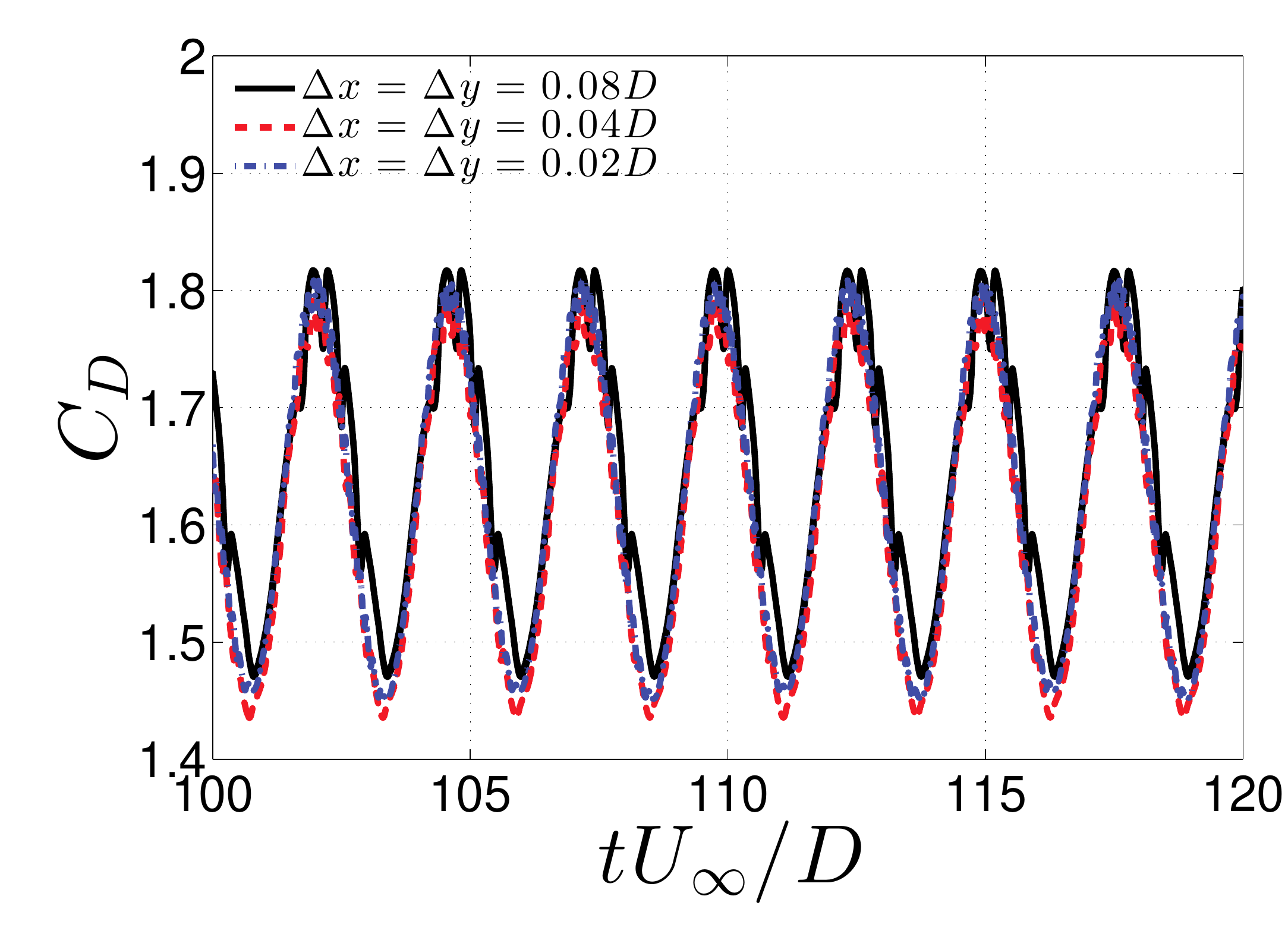}
    \label{fig_oscillating_cylinder_crossflow_Lag}
  }
  \caption{
   Temporal evolution of drag coefficient measured by the 
   \subref{fig_oscillating_cylinder_crossflow_CV} present control volume approach and
    \subref{fig_oscillating_cylinder_crossflow_Lag} Lagrange multiplier approach.
    Here, (---, black): $\dx = \dy = 0.08D$;
    (\texttt{---}, red): $\dx = \dy = 0.04D$;
    (\texttt{-}$\cdot$\texttt{-}, blue): $\dx = \dy = 0.02D$.}
  \label{fig_oscillating_cylinder_crossflow}
\end{figure}

The domain is discretized with three different uniform meshes of sizes $400 \times 192$, $800
\times 384$, and $1600 \times 768$, which corresponds to grid spacings of $\dx = \dy = 0.08D$,
$\dx = \dy = 0.04D$, and $\dx = \dy = 0.02D$, respectively. A constant time step size of $\dt  = 0.005 
D/U_{\infty}$ is used for all the computations, matching that of Lee at al. ~\cite{Lee11}. 
\REVIEW{The time step size satisfies the convective CFL condition with $C = 0.7$, which is found to be stable for
all the mesh sizes considered here.}
Fig.~\ref{fig_oscillating_cylinder_crossflow} shows
the time evolution of drag coefficient $C_D = \cF \cdot \e_x/(0.5 \rho U_{\infty}^2D)$. Again, both
the CV and LM approach produce identical results  and do not incur spurious force oscillations 
even at coarse resolutions as observed by Lee et al. (see Fig.~16 in \cite{Lee11}).

\subsubsection{Rotational oscillation}
As a last example of this section, we consider the case of a cylinder undergoing a rotational 
oscillation about its center of mass in a quiescent flow. The diameter of the cylinder is taken to 
be $D = 1$ and is placed in a domain of size $[-20D, 20D] \times [-20D, 20D]$, with zero velocity 
prescribed on all computational boundaries. The initial center of mass of the cylinder is placed 
at $(x,y) = (0,0)$ and it rotates about its center with a velocity $\omega_b = A_m \sin(2 \pi f t)$, 
in which $f$ is the frequency of oscillation and $T = 1/f$ is the time period of oscillation.
The density of the fluid is set to be $\rho = 1$.
The Reynolds number of the flow is $\Re = \rho U_m D/\mu = 300$, in which $U_m = A_m D/2$. The cylinder
rotates with frequency $f = 0.1$ and has maximum angular velocity $A_m = 10 f D$. These parameters are
chosen to match Borazjani et al.~\cite{Borazjani13}.

\begin{figure}[H]
  \centering
    \includegraphics[scale = 0.32]{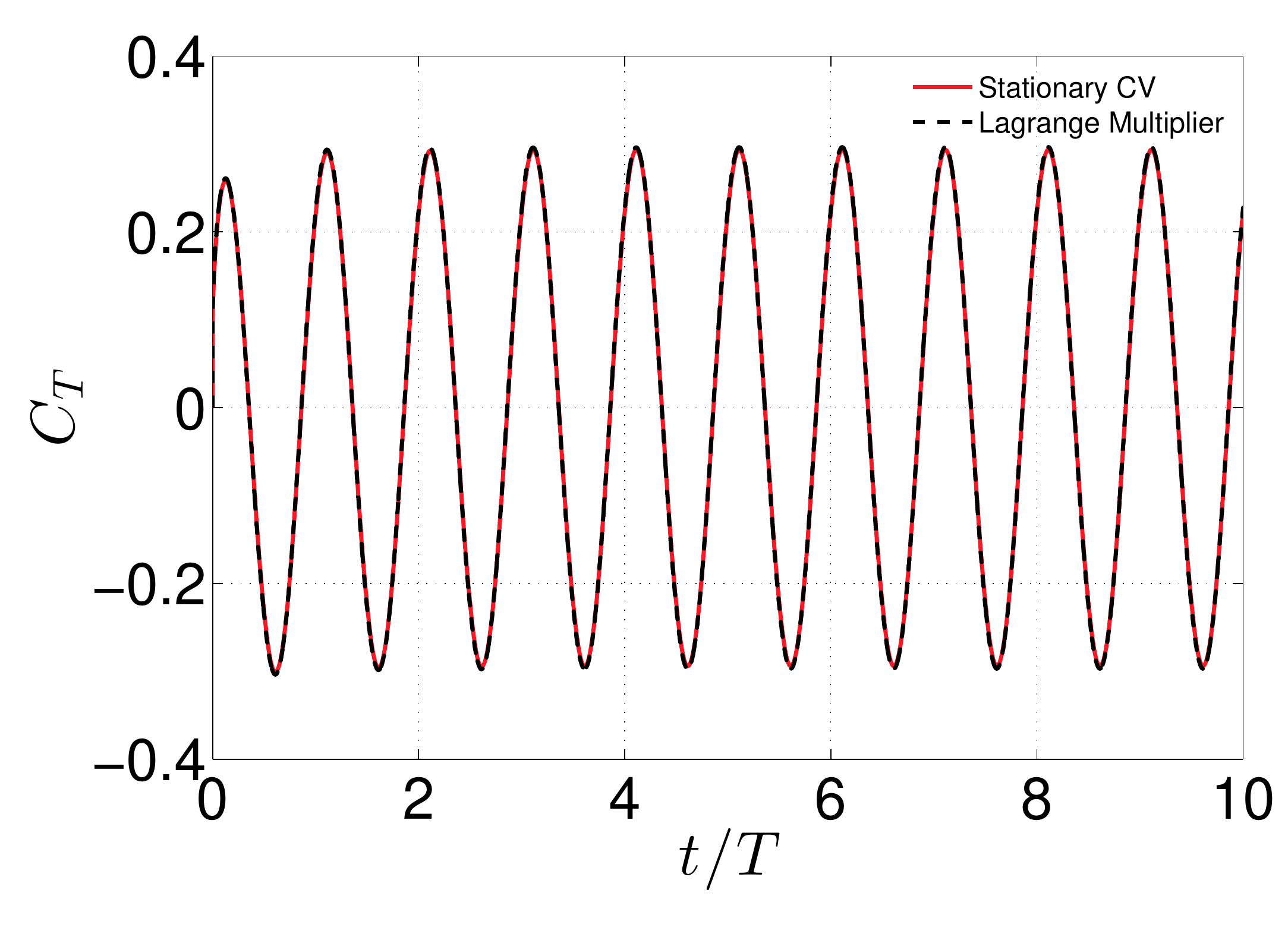}
  \caption{Comparison of the torque coefficient for a cylinder undergoing 
  rotational oscillation at $\Re = 300$ measured in two different ways.
   Here, (---, red): Control volume using Eq.~\eqref{eq_OurTorque};
     (\texttt{---}, black): Lagrange multiplier using Eq.~\eqref{eq_extrinsic_torque}.}
  \label{fig_rotating_cylinder}
\end{figure}

The time step size is chosen to be $\Delta t = 1\times 10^{-4} T$ and the grid is discretized by a two level
mesh, which consists of a coarse mesh of size $512 \times 512$, and an embedded fine mesh with 
refinement ratio $\nref = 4$. The minimum grid spacing at the finest level is 
$\dx_\textrm{min} = \dy_\textrm{min} = 0.0195D$. The structure remains on the finest grid level for at all time instances. 
A stationary control volume is placed at $[-1.01562D,1.01562D]^2$ and the torque coefficient is computed
as $C_T = \cM \cdot \e_z/(0.5 \rho U_m^2 D)$. Fig.~\ref{fig_rotating_cylinder} shows the temporal evolution
of $C_T$, which is in excellent agreement with sharp-interface CURVIB method results of Borazjani et 
al.~\cite{Borazjani13}.

\subsection{Moving plate}
In this section we consider the effect of regridding on hydrodynamic force calculations 
using the control volume approach. In the context of an immersed body AMR framework, 
regridding occurs when the body has moved some distance, or when new flow 
features of interest have appeared in the computational domain that require additional
mesh refinement to resolve them adequately.  We consider a moving plate example to 
understand the jumps in the drag coefficient due to regridding, and to provide 
some strategies to mitigate them.
 
\subsubsection{Moving plate on a uniform mesh}
\label{sec_translating_plate_uniform_mesh}
First we consider a single level domain case in two spatial dimensions. A finite plate of 
height $b = 1$ is dragged perpendicular to itself with constant velocity 
$\Ub = (U_\text{b}, V_\text{b}) = (-1,0)$ in an infinite fluid at rest. The plate is modeled 
as a thin line of points, separated by grid cell size distance in the transverse direction.
The physical domain is a periodic box of dimension $32b \times 22b$ and the domain
is discretized by a uniform Cartesian grid of size $1024 \times 1024$. The initial 
location of the control volume is $[-2b,2b] \times [-b,b]$, and it moves with an arbitrary 
speed to enclose the plate at all time instances. The density of the fluid is set to be $\rho  = 1$.
The Reynolds number of the flow 
is $\Re = \rho U_\text{b} b/\mu = 20$. The drag coefficient is calculated as 
$C_D = \cF \cdot \e_x/(\rho U_\text{b}^2 b/2)$. An asymptotic
solution $C_D \approx 2.09$ was derived by Dennis et al.~\cite{Dennis93}, and this
problem was also studied numerically by Bhalla et al.~\cite{Bhalla13}.

\begin{figure}[H]
  \centering
    \includegraphics[scale = 0.3]{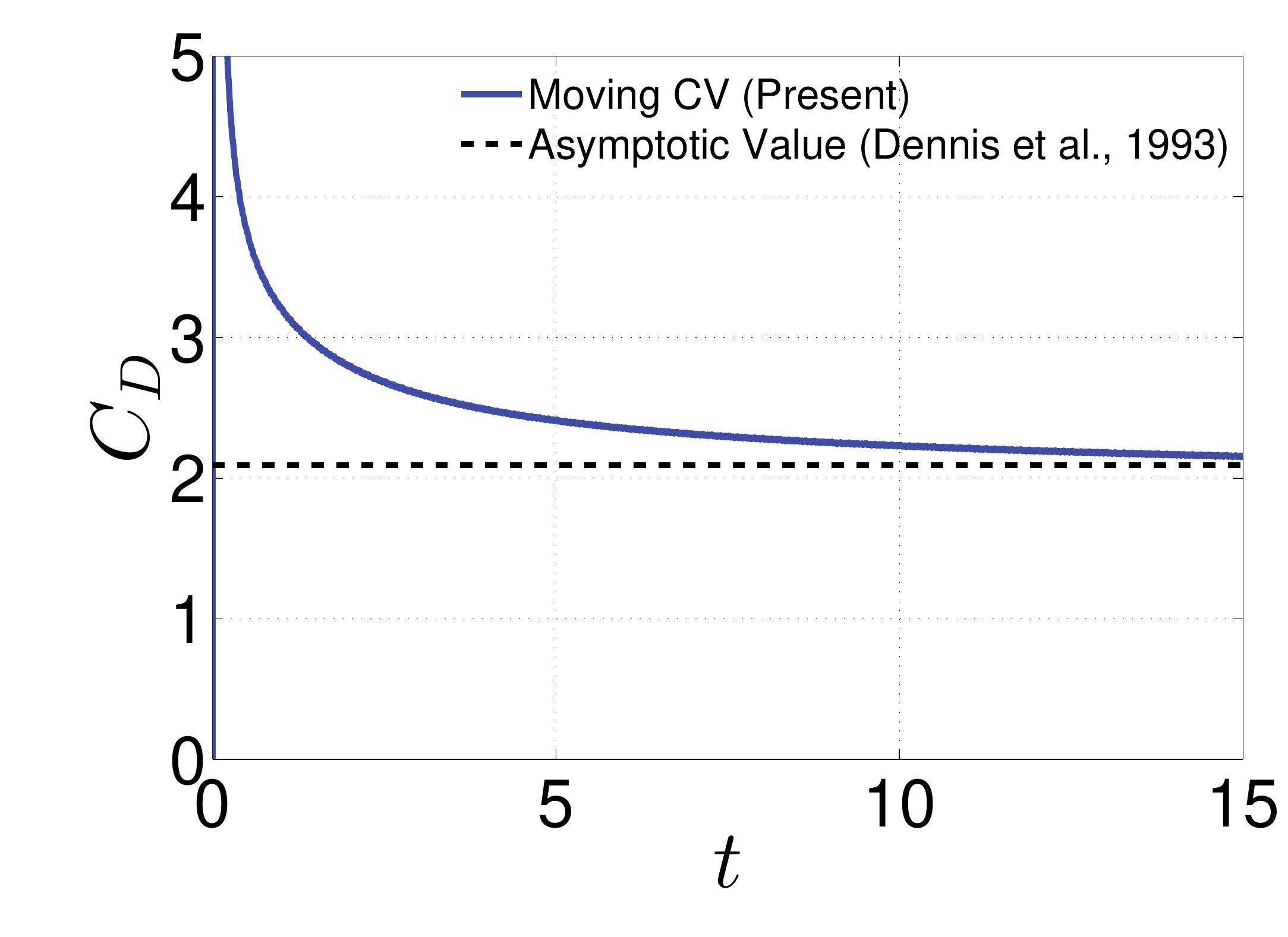}
  \caption{Comparison of the drag coefficient for a translating 
    plate at $\Re = 20$ measured by the present moving control
    volume approach (---, blue); 
     \texttt{---}: asymptotic value from Dennis et al.~\cite{Dennis93}.}
  \label{fig_moving_plate_uniform}
\end{figure}

Fig.~\ref{fig_moving_plate_uniform} shows the time evolution of drag coefficient
for the moving plate. We see that the numerical solution obtained by the moving CV
computation matches well with the asymptotic value derived in~\cite{Dennis93} and does 
not contains any spurious force oscillations or jumps.

\subsubsection{Moving plate on an adaptive mesh} \label{sec_translating_plate}
Next, we consider the same moving plate example but with locally refined grids. The domain 
is discretized by a coarse grid of size $256 \times 256$. A more refined mesh immediately 
surrounds the plate with refinement ratio $\nref = 4$, giving the finest level an 
equivalent grid size of $1024 \times 1024$. The mesh is adaptive in the sense that it 
selectively refines in areas with large velocity gradients and where the immersed body is
located in the domain. Apart from the locally refined grids, we use the same parameters  
of Sec.~\ref{sec_translating_plate_uniform_mesh} for this case.

\begin{figure}[H]
  \centering
  \subfigure[]{
    \includegraphics[scale = 0.27]{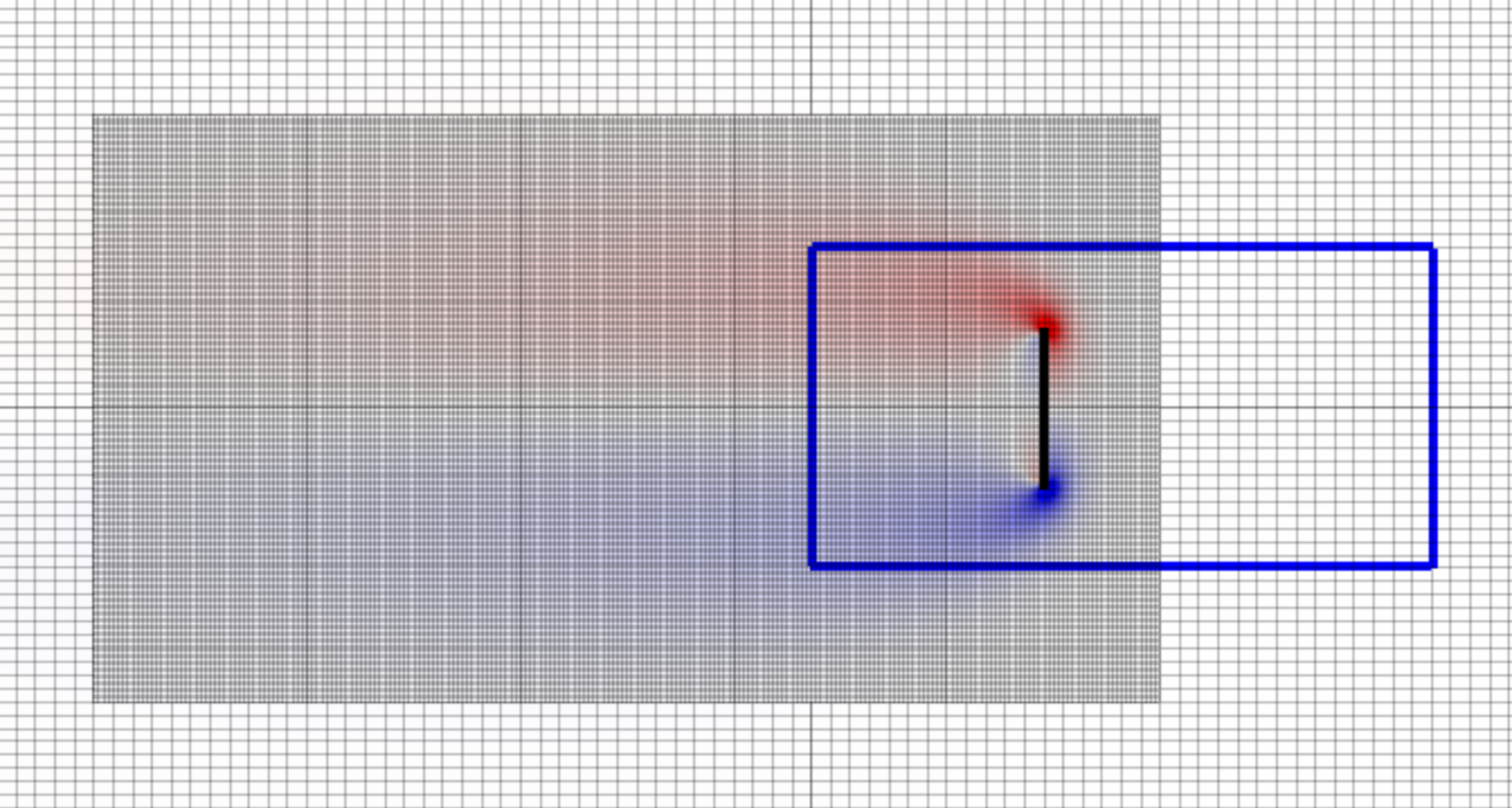}
    \label{fig_amr_plate_multiple_levels_moving_viz}
  }
   \subfigure[]{
    \includegraphics[scale = 0.28]{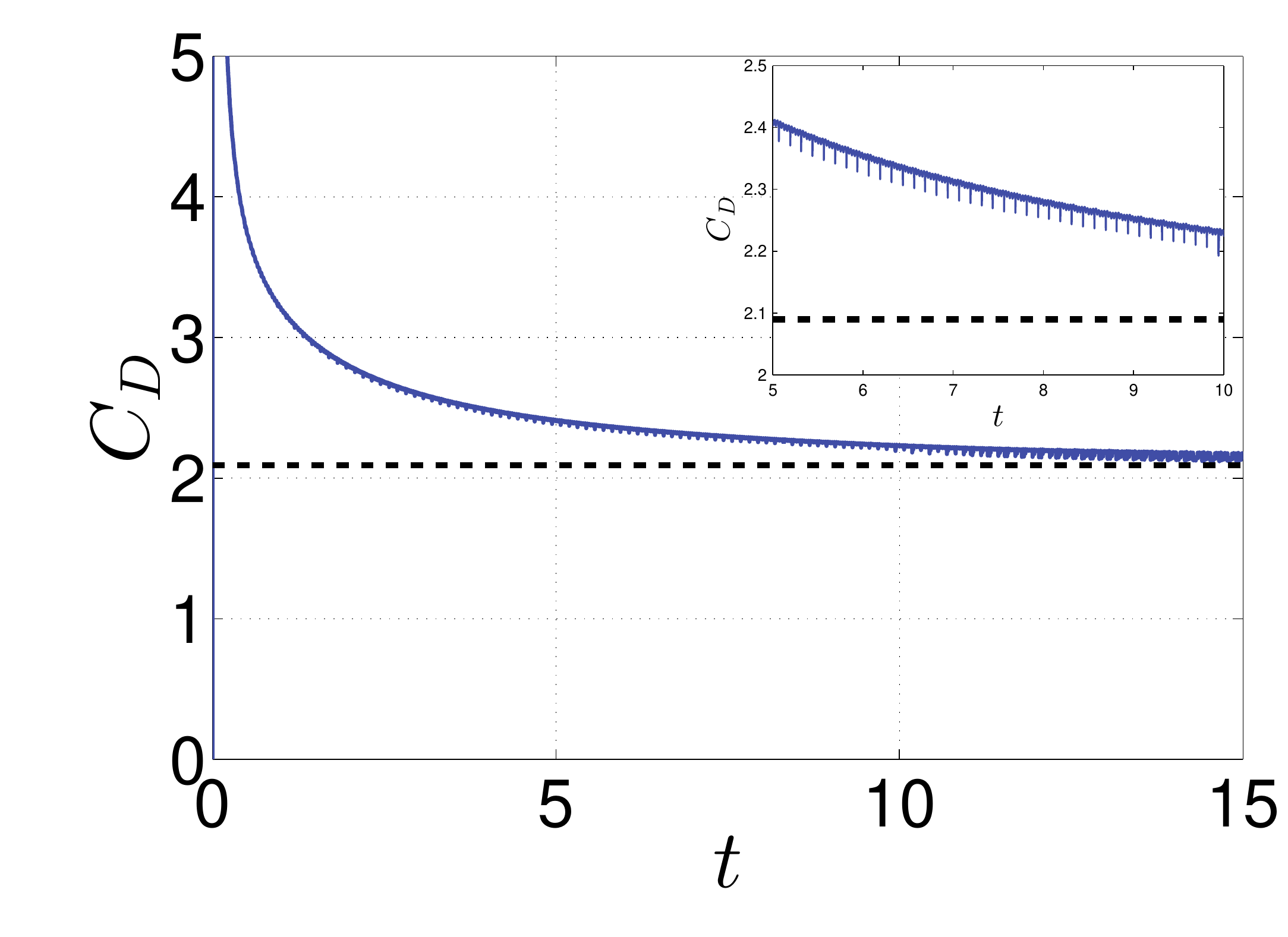}
    \label{fig_amr_plate_multiple_levels_moving_CD}
  }
  \caption{\subref{fig_amr_plate_multiple_levels_moving_viz} 
  Vorticity generated by a moving flat plate at $t=7.5$ for $\Re = 20$. The \emph{moving} control
  volume spanning both the coarsest and finest grid level is shown in blue.
  The plotted vorticity is between $-7.5$ and $7.5$.
  \subref{fig_amr_plate_multiple_levels_moving_CD}
   Temporal evolution of drag coefficient measured by the present \emph{moving} control
    volume approach (---, blue); 
     \texttt{---}: asymptotic value from Dennis et al.~\cite{Dennis93}.}
  \label{fig_amr_plate_moving_multiple_levels}
\end{figure}

\begin{figure}[H]
  \centering
  \subfigure[]{
    \includegraphics[scale = 0.4]{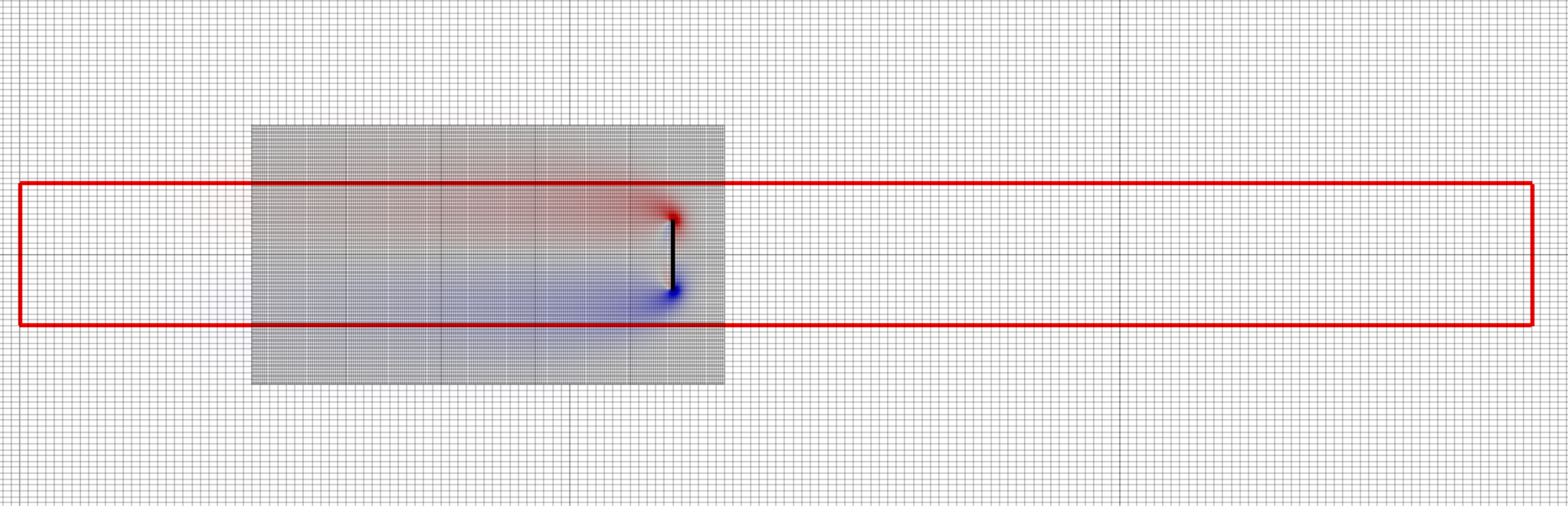}
    \label{fig_amr_plate_multiple_levels_stationary_viz}
  }
   \subfigure[]{
    \includegraphics[scale = 0.28]{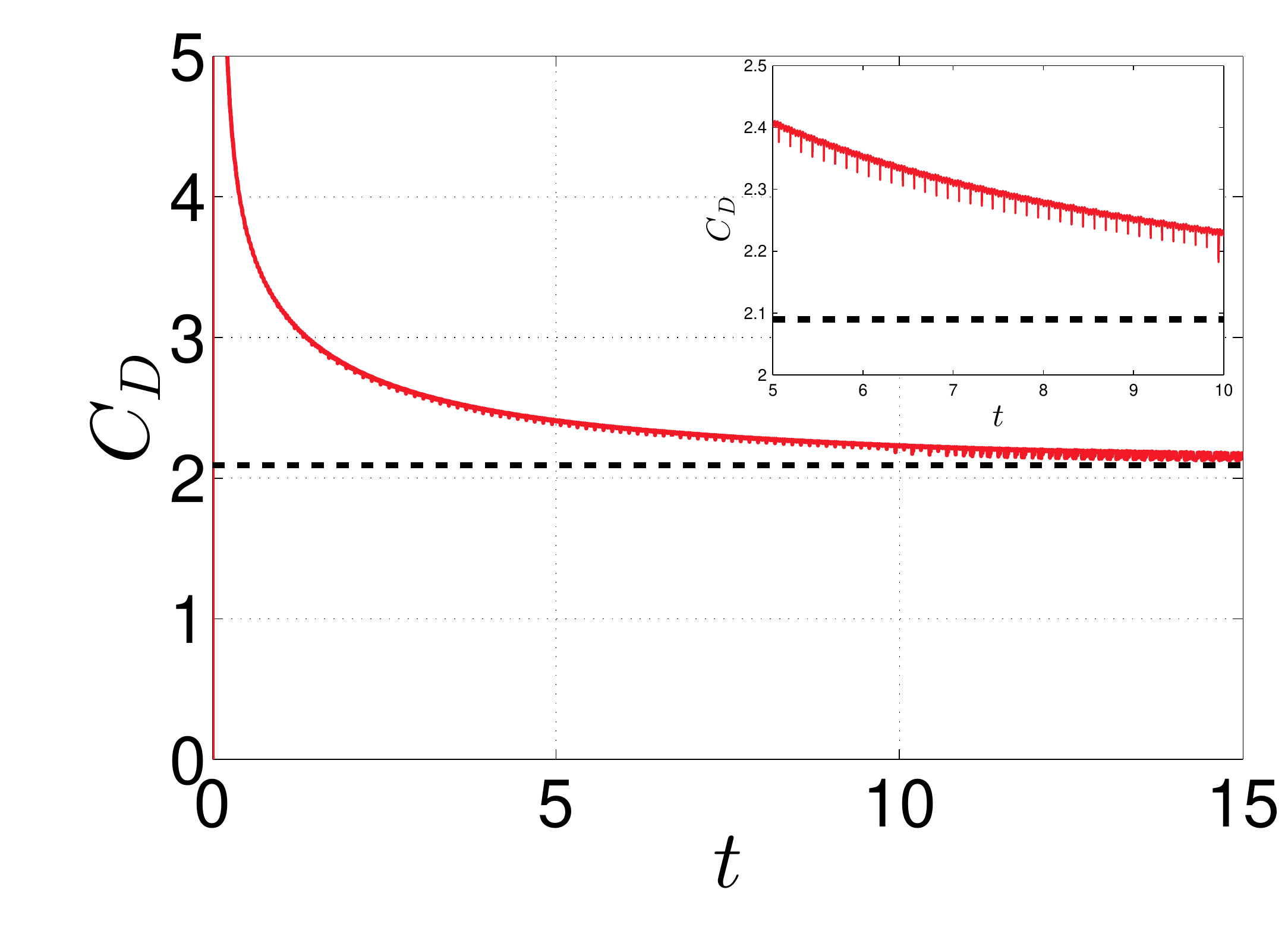}
    \label{fig_amr_plate_multiple_levels_stationary_CD}
  }
  \caption{\subref{fig_amr_plate_multiple_levels_moving_viz} 
  Vorticity generated by a moving flat plate at $t=7.5$ for $\Re = 20$. The \emph{stationary} control
  volume spanning both the coarsest and finest grid level is shown in red.
  The plotted vorticity is between $-7.5$ and $7.5$.
  \subref{fig_amr_plate_multiple_levels_stationary_CD}
   Temporal evolution of drag coefficient measured by the present control
    volume approach (---, red); 
     \texttt{---}: asymptotic value from Dennis et al.~\cite{Dennis93}}
  \label{fig_amr_plate_stationary_multiple_levels}
\end{figure}

\begin{figure}[H]
  \centering
  \subfigure[]{
    \includegraphics[scale = 0.3]{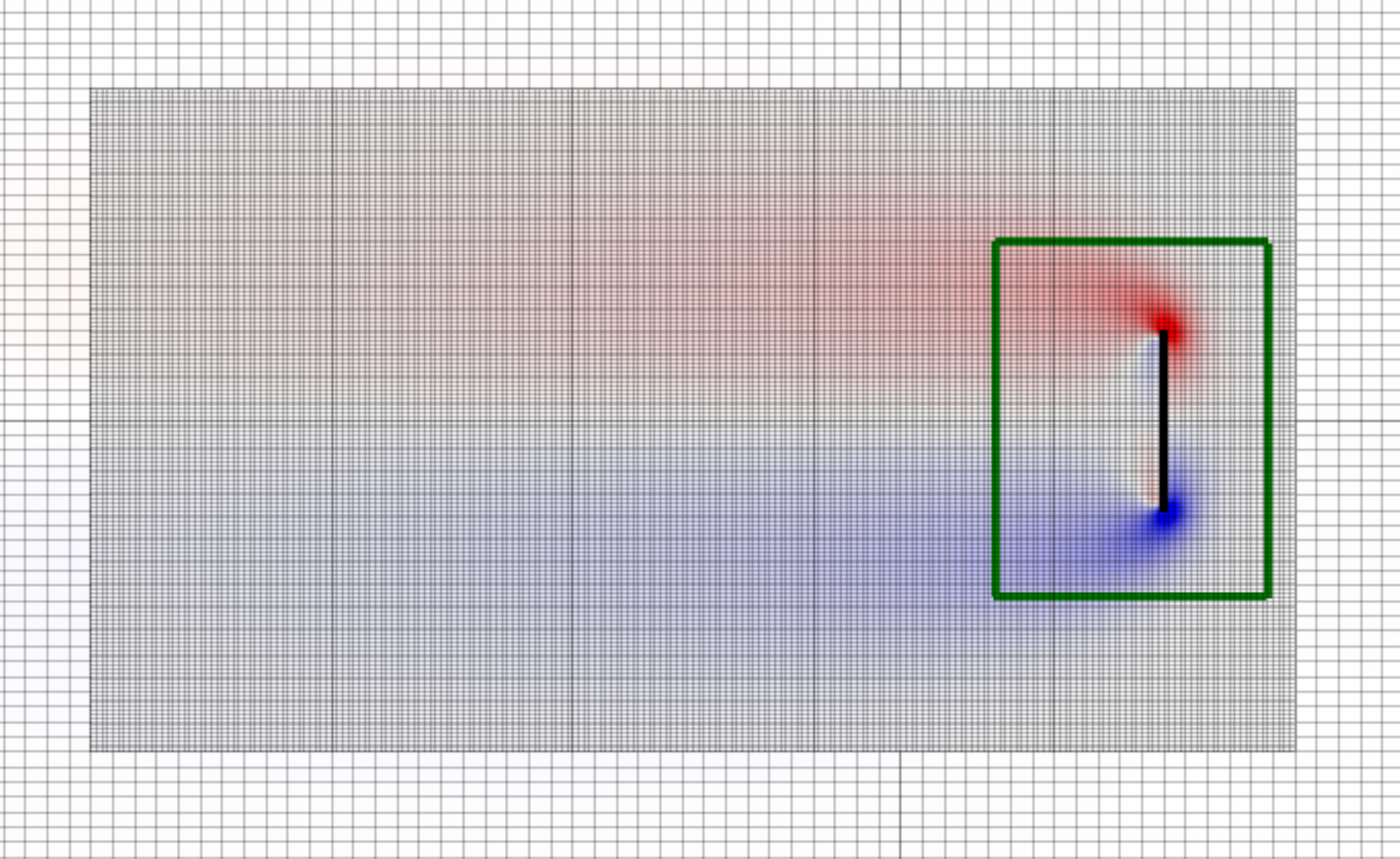}
    \label{fig_amr_plate_single_level_moving_viz}
  }
   \subfigure[]{
    \includegraphics[scale = 0.3]{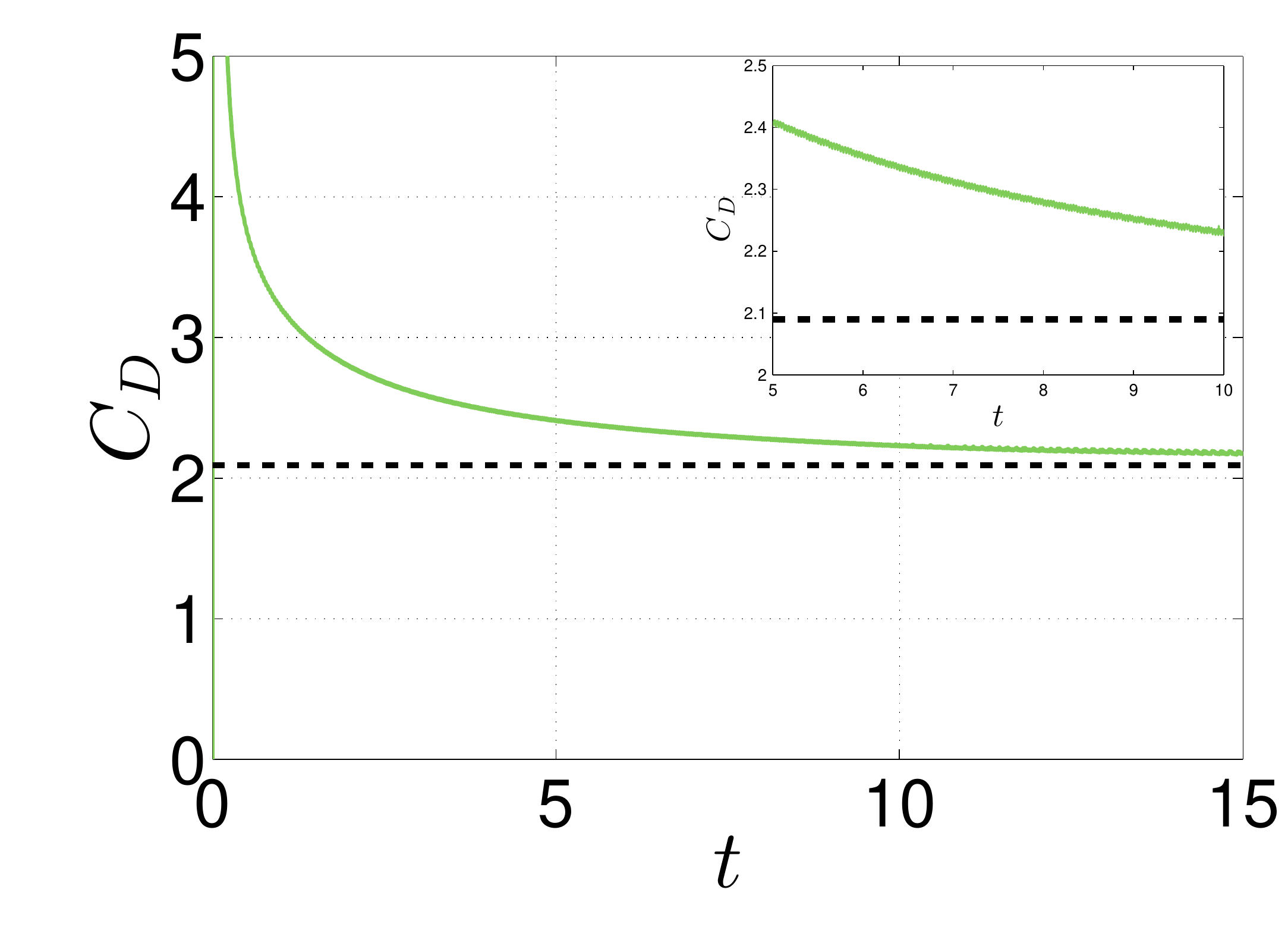}
    \label{fig_amr_plate_single_level_moving_CD}
  }
  \caption{\subref{fig_amr_plate_single_level_moving_viz} 
  Vorticity generated by a translating flat plate at $t=7.5$ for $\Re = 20$. The 
  \emph{moving} control volume spanning both the coarsest and finest grid level 
  is shown in green. The plotted vorticity is between $-7.5$ and $7.5$.
  \subref{fig_amr_plate_single_level_moving_CD}
   Temporal evolution of drag coefficient measured by the present \emph{moving} control
    volume approach (---, green); 
     \texttt{---}: asymptotic value from Dennis et al.~\cite{Dennis93}.}
  \label{fig_amr_plate_moving_single_level}
\end{figure}

Three different control volumes are used. First, the CV is set initially to $[-2b,2b]\times
[-1.03125b,1.03125b]$ and it translates to the right along with the plate. In this
configuration, the CV spans multiple levels of the grid hierarchy.
Fig.~\ref{fig_amr_plate_multiple_levels_moving_CD} shows the measured drag coefficient
over time. Again, $C_D$ evolves towards the asymptotic value derived in~\cite{Dennis93}.
However, there are small jumps in $C_D$ over time, which are absent from
the uniform mesh case of Sec.~\ref{sec_translating_plate_uniform_mesh}.
These jumps can be attributed to the mesh hierarchy regridding to follow the moving
plate or due to the moving CV itself. To rule out the possibility of jumps due to the motion 
of control volume, we consider a second CV configuration in which the CV is held stationary all 
times. The stationary CV again spans multiple grid levels and is big enough to contain the moving
plate at all time instances. Fig.~\ref{fig_amr_plate_stationary_multiple_levels} shows 
the CV configuration and measured drag coefficient over time for this CV 
$ = [-2b, -1.03125b] \times [20b, 1.03125b]$.  The jumps in $C_D$ are again observed.
In our third configuration, we limit the moving CV on the finest grid level. The CV configuration
and $C_D$ temporal profile is shown in Fig.~\ref{fig_amr_plate_moving_single_level}.
The initial location of this control volume is $[-b,0.625b] \times [-1.03125b,1.03125b]$.

\noindent From Fig.~\ref{fig_amr_plate_single_level_moving_CD}, we note 
that the jumps due to regridding can be substantially mitigated if the moving 
CV is restricted to the finest mesh level. The CV translates to the right along 
with the plate, but remains on the finest mesh level throughout the simulation.
There are no longer jumps in the computed drag values for $t \le 10$.

The jumps in the hydrodynamic forces due to regridding occur due to 
velocity reconstruction following a regridding operation. A common velocity
reconstruction strategy employed in an AMR framework is to first interpolate velocities 
from the coarser level to the new fine level, and then directly copy velocities from the old fine level 
in the spatial regions where old and new fine levels intersect. We refer readers to
Griffith et al.~\cite{Griffith07} for details.
When we restrict the moving CV to the finest grid level, we reduce errors in 
momentum change due to velocity reconstruction from coarse to fine level interpolation. 
The contribution of $\u^{n}$ in the time derivative term in Eqs.~\eqref{eq_OurForce} 
and~\eqref{eq_OurTorque} can be evaluated either before or after regridding. In our empirical 
tests, we have observed that evaluating the time derivative term after regridding helps in 
mitigating the jumps even further (comparison data not shown). In all our results shown 
above, we evaluate the contribution from $\u^{n}$ after regridding, i.e, using the old velocity 
at the new hierarchy configuration. 

\subsection{Swimming eel}
In this section we demonstrate that the moving control volume approach can be used 
to determine the hydrodynamic forces and torques on a free-swimming body. 
We consider a two-dimensional undulating eel geometry, which is adapted 
from~\cite{Kern06,Bhalla13}. The eel's reference frame is aligned with the $x$-axis 
and in this refrence frame the lateral displacement along $0 \le x \le L$ over its 
projected length $L$ is given by

\begin{equation}
\label{eq_eel_kinematics}
y(x,t) = 0.125 \frac{x+0.03125}{1.03125} \sin \left[2 \pi (x-t/T)\right].
\end{equation}
A backwards-traveling wave of the above form having a time period $T$ 
causes the eel to self-propel. The swimmer is taken to have a projected length 
$L = 1$, and time period $T = 1$. The Reynolds number based on 
$V_\textrm{max} = 0.785 L/T$, the maximum undulation velocity at the tail tip, 
is $\Re = \rho V_\textrm{max} L/\mu = 5609$. The undulations travel in the positive 
$x$-direction, thereby propelling the eel in the negative $x$-direction.
The density of the fluid is taken to be $\rho = 1$.

The eel's total velocity $\Ub = \Ur + \Wr \wedge \R + \Uk$ in the Lagrangian frame 
is given by the following components: its rigid linear $\Ur$ and angular $\Wr$ center of mass velocities
and its deformational velocity $\Uk$, which we assume to have zero net linear and angular momentum. 
The self-propulsion velocities are obtained using conservation of linear and angular momentum in 
the body domain $\Vbt$

\begin{align}
\Mb \Ur^{n+1} & = \int_{\Vbt} \rho \; \left(\cJ_h[\X^{n+\half}] \, \tilde{\u}^{n+1} \right) \dV, 
\label{eq_conserve_linear_mom}\\
\Ib \Wr^{n+1} & = \int_{\Vbt} \rho \; \R^{n+\half} \wedge \left(\cJ_h[\X^{n+\half}] \, \tilde{\u}^{n+1}  \right) \dV, \label{eq_conserve_angular_mom}
\end{align}
in which $\Mb$ and $\Ib$ are the mass and moment of inertia tensor of the body, respectively.
Having obtained these rigid body velocities, the body velocity at time $t^{n+1}$ required in 
Eq.~\eqref{eqn_lm_correction} is obtained as $\Ub^{n+1} = \Ur^{n+1} + \Wr^{n+1} \wedge \R^{n+\half} + \Uk^{n+1}$,
in which $\R^{n+\half} = \X^{n+\half} - \X^{n+\half}_0$ is the radius vector from (an estimated) midstep  
center of mass $\X_0^{n+\half}$ to (an estimated) midstep Lagrangian node position $\X^{n+\half}$.
We refer readers to~\cite{Bhalla13} for more details.

\begin{figure}[H]
  \centering
  \subfigure[$t/T= 1$]{
    \includegraphics[scale = 0.35]{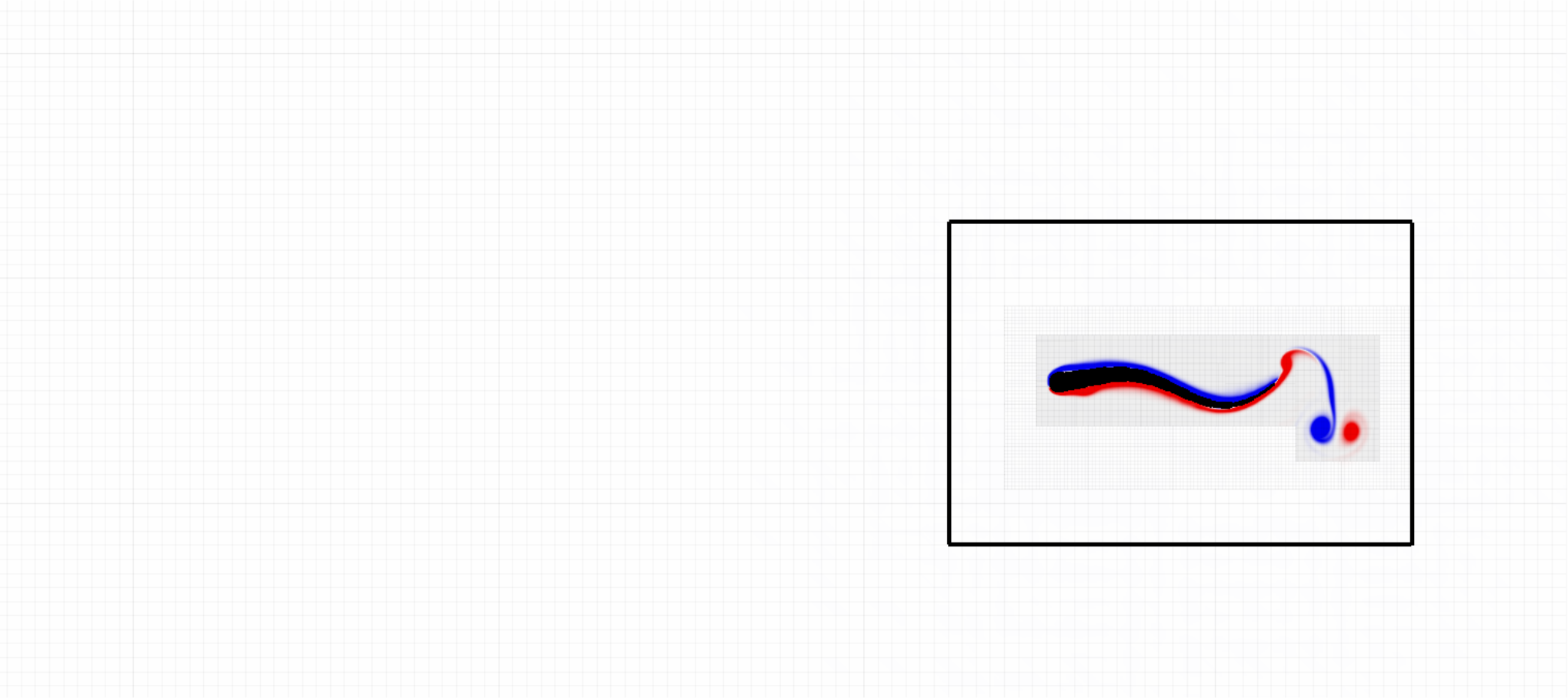}
    \label{fig_undulating_eel_viz1}
  }
   \subfigure[$t/T = 3$]{
    \includegraphics[scale = 0.35]{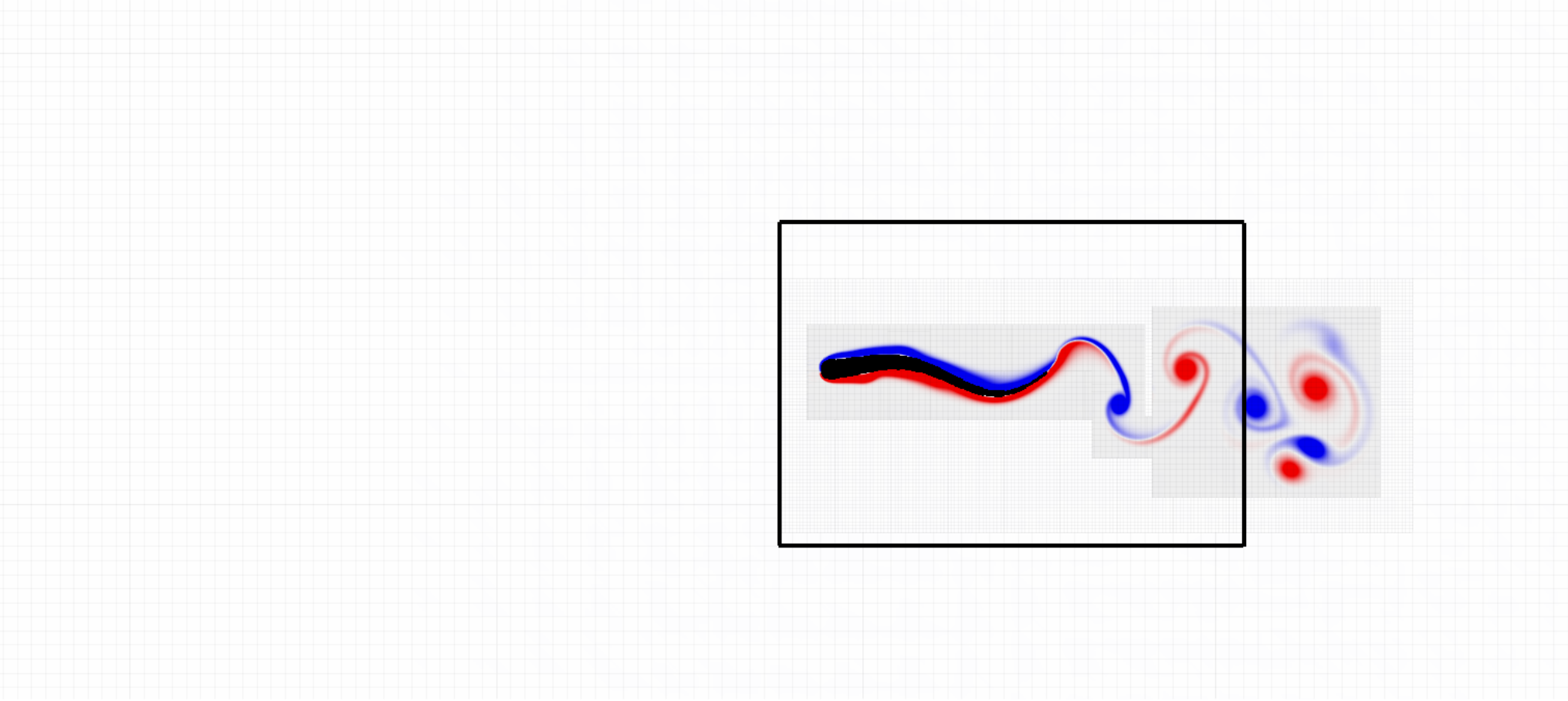}
    \label{fig_undulating_eel_viz3}
  }
    \subfigure[$t/T = 5$]{
    \includegraphics[scale = 0.35]{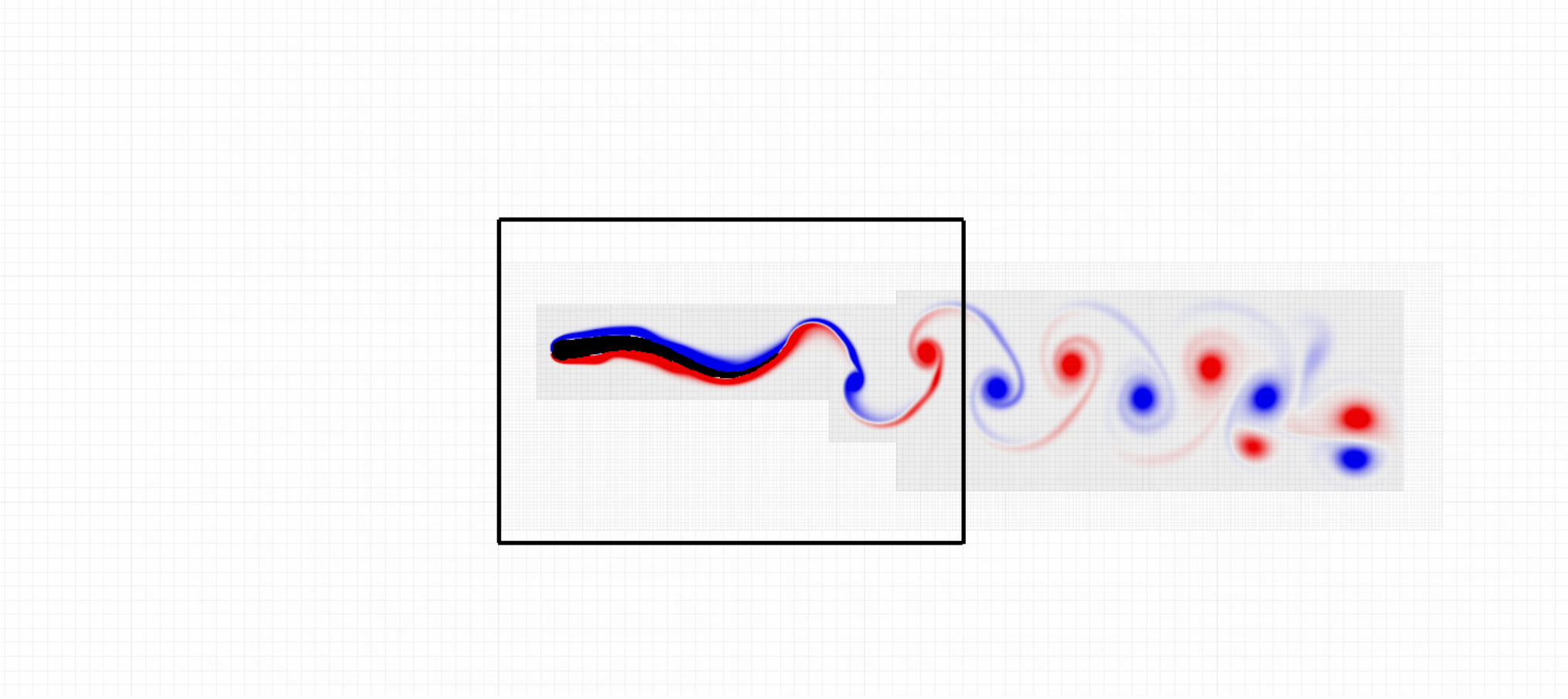}
    \label{fig_undulating_eel_viz5}
  }
   \subfigure[$t/T = 7$]{
    \includegraphics[scale = 0.35]{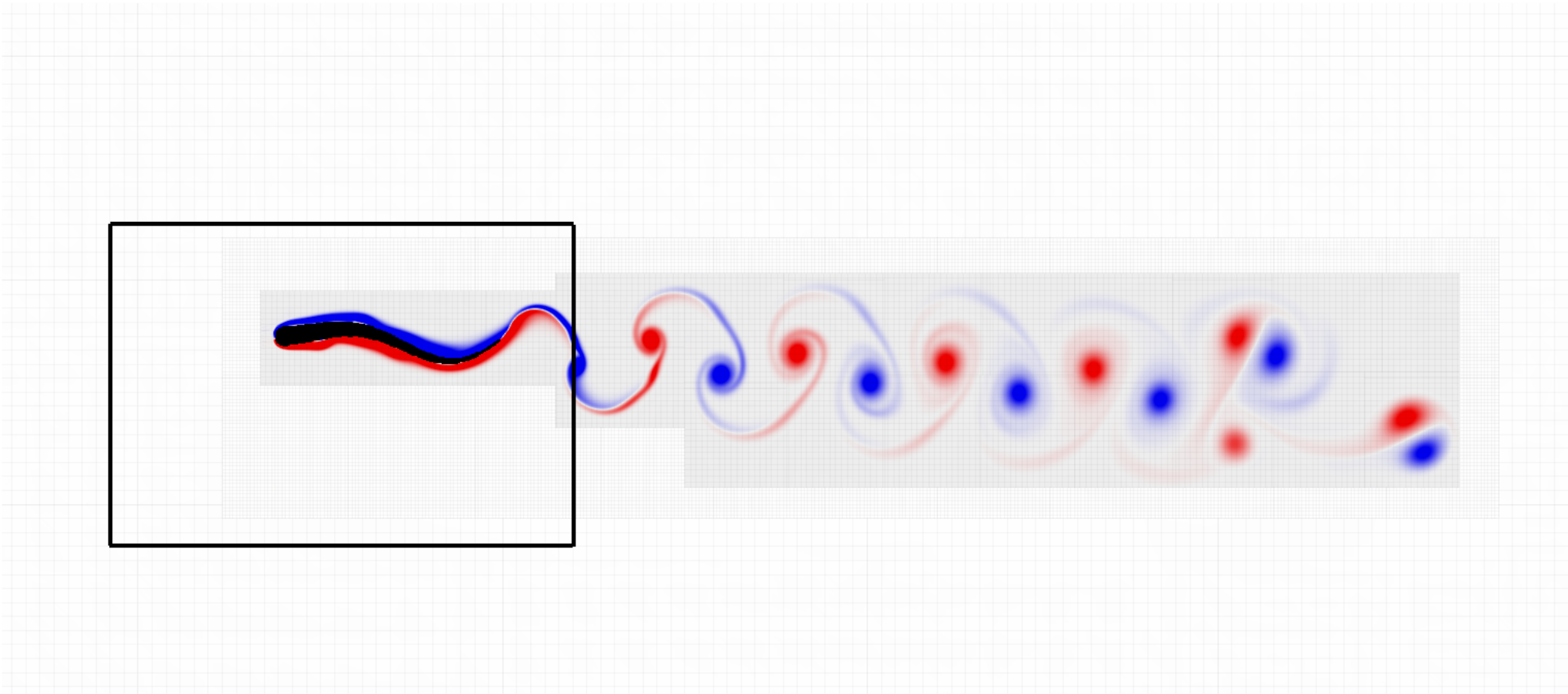}
    \label{fig_undulating_eel_viz7}
  }
  \caption{Vorticity generated by the two-dimensional eel model at $\Re = 5609$,
  along with the moving CV location at four different time instances on an adaptive mesh hierarchy. 
  All figures are plotted for vorticity between $-10$ and $10$.
   }
  \label{fig_eel_viz}
\end{figure}

The fully periodic domain is taken to be of size $8L \times 4L$ and is discretized with a 
three-level hierarchy of Cartesian grids. The size of the coarsest grid is  
$128 \times 64$ grid cells and $\nref = 4$ is taken for subsequent finer grids. Hence, the 
finest grid, with spacing equivalent to that of a uniform mesh of size $2048 \times 1024$, embeds 
the undulatory swimmer at all times. A time step size of $\dt = 1\times 10^{-4}\;T$ is employed. 
The head of the swimmer is initially centered at $(x,y) = (0,0)$ and its body extends in the 
positive $x$-direction. The CV is initially located at $[-1.02L, 1.0425L] \times [-0.7075L, 0.73L]$,
which encompasses the entire swimmer. The CV is allowed to span multiple grid levels.
Whenever the eel's center of mass translates
a distance $\dx$, the CV moves with velocity $\uS =  (-\dx/\dt,0)$ in order to remain aligned with
the grid lines. Fig.~\ref{fig_eel_viz} shows the vortical structures generated by the eel
at four separate time instances, along with the locations of the moving control volume.

\begin{figure}[H]
  \centering
 \subfigure[]{
    \includegraphics[scale = 0.32]{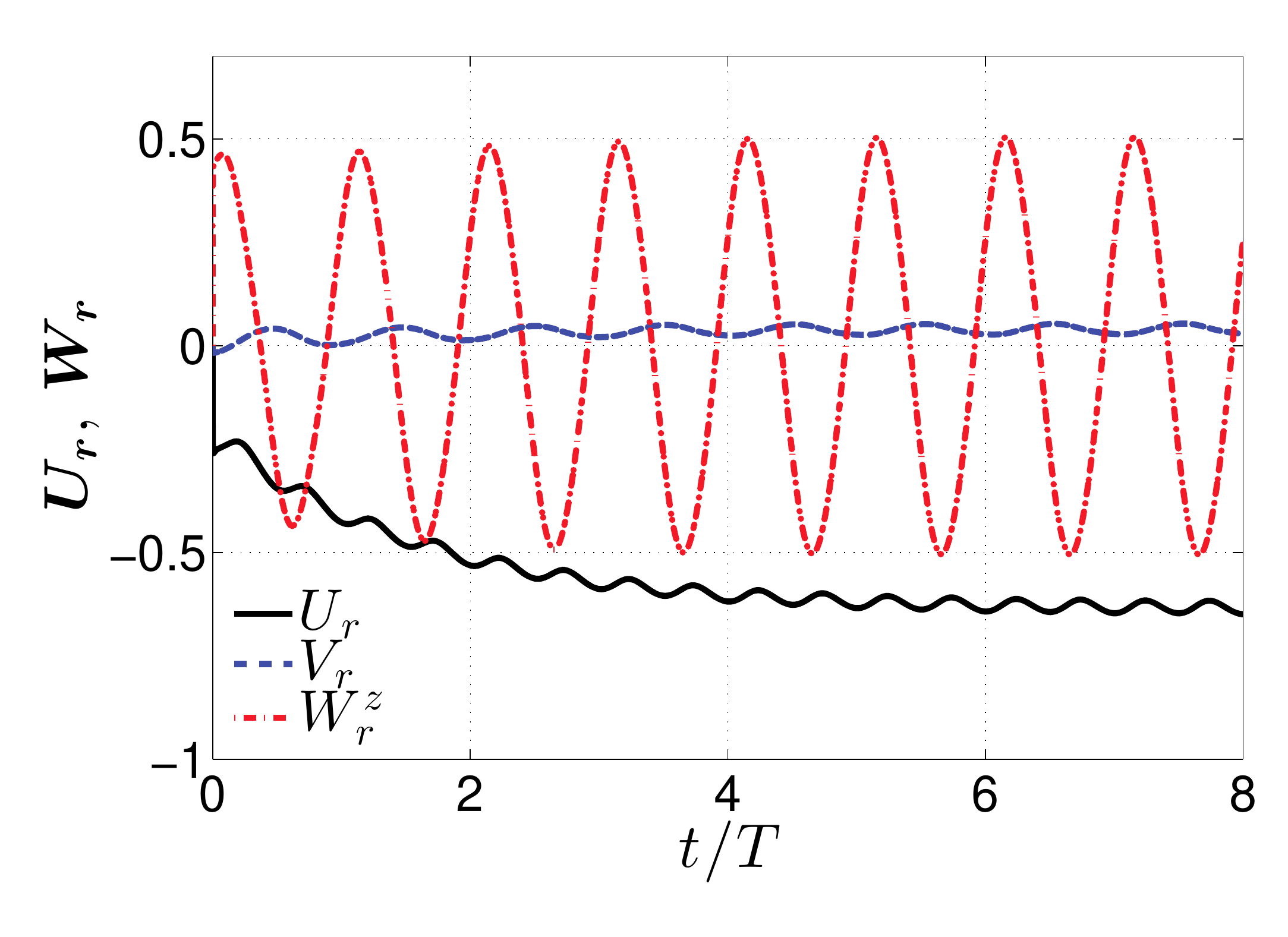}
    \label{fig_undulating_eel_speed}
  }
  \subfigure[]{
    \includegraphics[scale = 0.32]{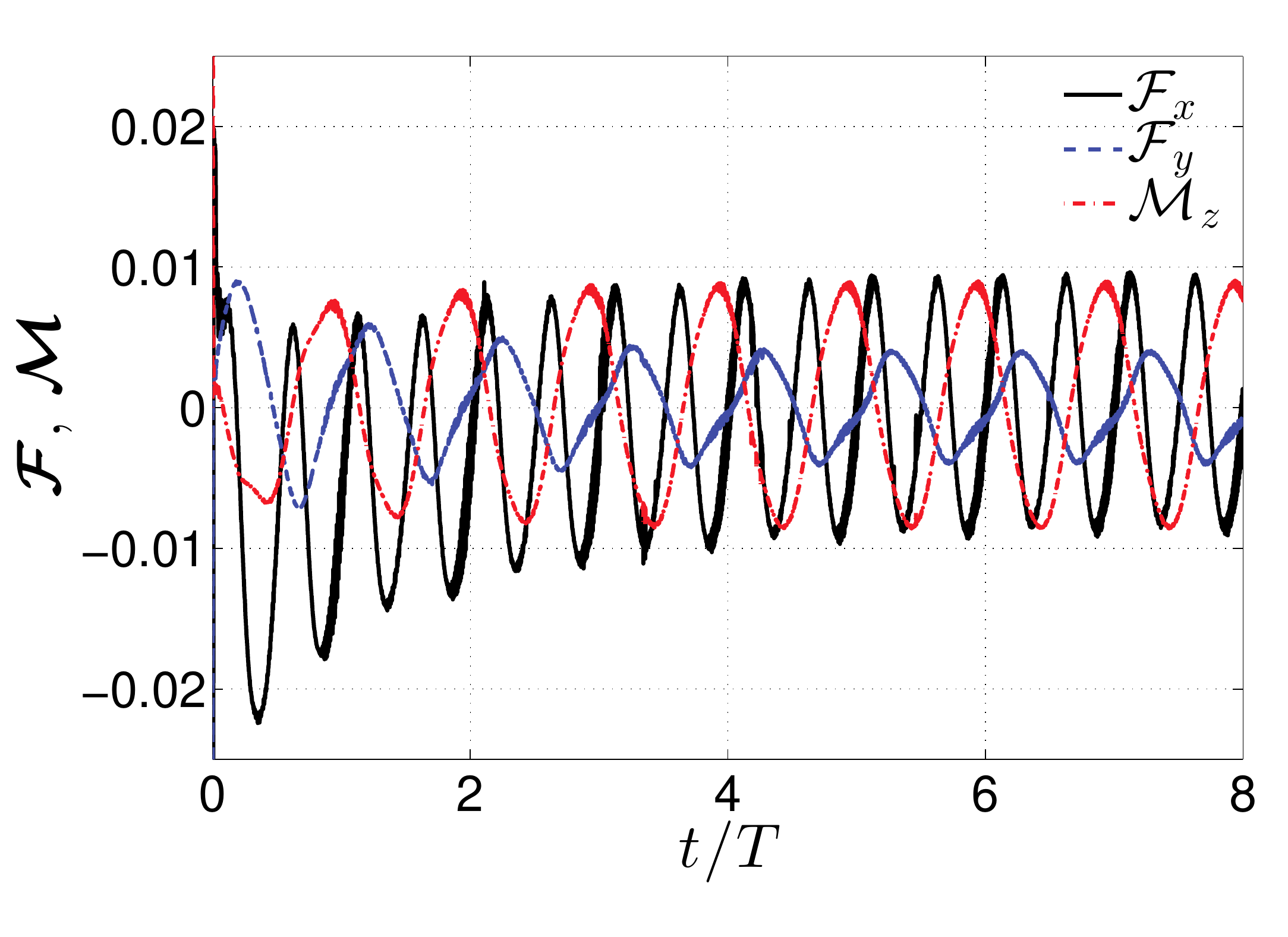}
    \label{fig_undulating_eel_force}
  }
   \subfigure[]{
    \includegraphics[scale = 0.32]{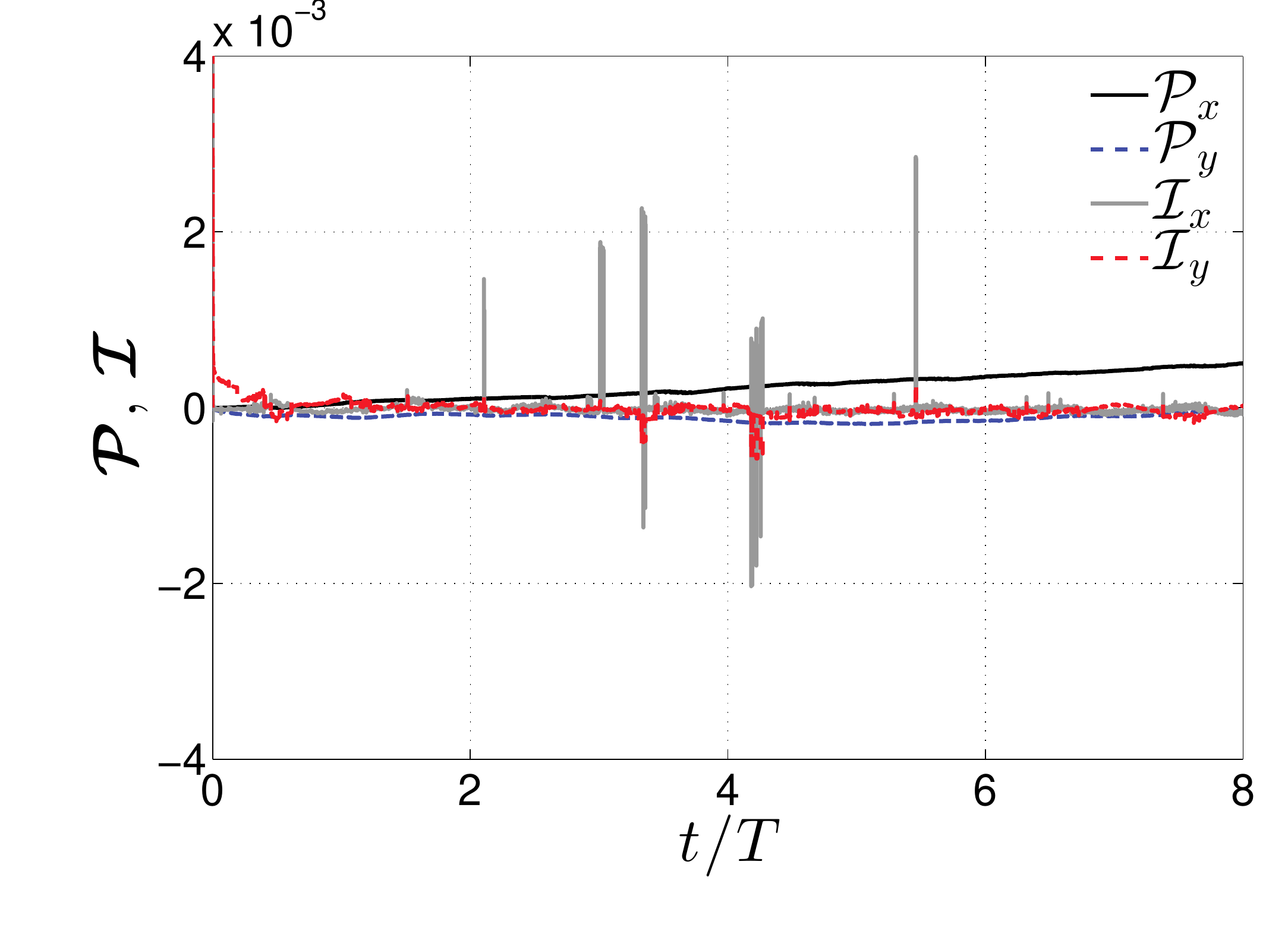}
    \label{fig_undulating_eel_momentum}
  }
  \caption{
  \subref{fig_undulating_eel_speed}
  Temporal evolution of (---, black) axial $U_r$ \& (\texttt{---}, blue) lateral $V_r$ swimming velocity,
  and (\texttt{-}$\cdot$\texttt{-}, red) rotational $W_r^z$ velocity of the eel calculated by Eqs.~\eqref{eq_conserve_linear_mom} 
  and \eqref{eq_conserve_angular_mom}.
  \subref{fig_undulating_eel_force}
   Temporal evolution of net (---, black) axial $\mathcal{F}_x$ \& (\texttt{---}, blue) lateral $\mathcal{F}_y$ forces,
   and (\texttt{-}$\cdot$\texttt{-}, red) net torque $\mathcal{M}_z$ on the body of the eel 
   measured by the present \emph{moving} control volume approach.
  \subref{fig_undulating_eel_momentum}
   Temporal evolution of (---, black) axial $\cP_x$ and (\texttt{---}, blue) lateral $\cP_y$ momentum
   of the entire fluid domain. Temporal evolution of change in linear momentum within
   the control volume: (---, gray) $\mathcal{I}_x$ and (\texttt{---}, red) $\mathcal{I}_y$.}
  \label{fig_undulating_eel}
\end{figure}

Fig.~\ref{fig_undulating_eel_speed} shows the time evolution of axial ($U_r = \Ur \cdot \e_x$) and lateral
($V_r = \Ur \cdot \e_y$) swimming velocities, along with the rotational velocity $W^z_r = \Wr \cdot \e_z$. The eel
is shown to travel in the $-\e_x$ direction, eventually reaching a steady state speed. The angular velocity
oscillates about a zero mean value, while the lateral velocity has small non-zero mean due to initial transients.
Fig.~\ref{fig_undulating_eel_force} shows the time evolution of net axial ($\mathcal{F}_x = \cF \cdot \e_x$) and lateral
($\mathcal{F}_y = \cF \cdot \e_y$) forces acting on the eel's body, along with the net torque $\mathcal{M}_z = \cM 
\cdot \e_z$, which is measured from the eel's center of mass. The forces and torque are computed using  
the moving CV approach, although identical estimates are obtained from the LM approach (data not shown).
Both $\cF$ and $\cM$ oscillate about a mean value of zero, which is expected during free-swimming as there
are no external forces and torques applied on the swimmer.

For an object initially at rest in a periodic and quiescent fluid, the net linear momentum over 
the entire computational domain should remain zero~\cite{Bhalla13}, i.e., $\vcP(t) = \int_{\Omega} \rho \u \dV = \V0$.
Similarly, the net angular momentum of the system should also remain zero at all times, i.e.,
$\vcL(t) = \int_{\Omega} \rho \rcross \u \dV = \V0$. 
In other words, all of the momentum generated due to the eel's vortex shedding should be
redistributed to the eel's translational and rotational motion. Moreover, the change in linear momentum
of the body should be equal to the net force on the body during free-swimming. Hence, the net force
on the body should be given by $\cF(t) = \d{}{t} \int_{\Vb} \rho \u \dV$, implying that

\begin{equation}
\label{eq_CV_other}
\cI(t) = - \int_{\Vcvt} \rho \D{\u}{t} \dV 
+ \oint_{\Scvt} \ndot \left[-p \I - \u\rho \u + \T \right] \dS = \V0.
\end{equation}
The above statement also implies the conservation of linear momentum for the control volume.
Eq.~\eqref{eq_CV_other} also implies that the sum of the Langrange multipliers enforcing the
rigidity constraint within $\Vbt$ is zero during free-swimming (see Eq.~\eqref{eq_sumLag_CV}).
Fig.~\ref{fig_undulating_eel_momentum} shows the temporal evolution of $\cP_x = \vcP \cdot \e_x$, 
$\cP_y = \vcP \cdot \e_y$, $\mathcal{I}_x = \cI \cdot \e_x$, and $\mathcal{I}_y = \cI \cdot \e_y$; it is indeed seen
that these quantities are nearly zero for all time instances.
The slight increase in $\cP_x$ is attributed to spatial and temporal discretization errors,
whereas the jumps in $\mathcal{I}_x$ correspond to time steps at which regridding occurs.
Similar observations are made for angular momentum conservation for the entire system (data not shown here).

\REVIEW{
\subsection{Drafting, kissing, and tumbling} \label{sec_dkt_cylinders}
In this section we simulate the dynamic interactions between two sedimenting cylindrical particles and use 
the moving control volume and Lagrange multiplier approaches to determine the hydrodynamic forces. 
The cylinders are identically shaped with diameter $D = 0.2$ cm and are placed in a domain 
of size $[-5D, 5D] \times$ [0, 40D], with zero velocity prescribed on the left and right boundaries, 
and with axial and transverse tractions set to zero at the top and bottom boundaries. The density
and viscosity of the fluid are set to $\rho = 1.0$ g/cm$^3$ and $\mu = 0.01$ g/(cm $\cdot$ s), 
respectively. Each particle is subject to a gravitational body force 
$\F^{g} = -(\rho_s - \rho) g V_p \e_y$, where $g=980$ cm/s$^2$ is the gravitational
constant, $\rho_s = 1.01 \rho$ is the density of the solid, 
and $V_p = \pi (D/2)^2$ is the volume of each particle. This is realized through an Eulerian
body force $\f^g$ added to the right-hand side the momentum Eq.~\eqref{eqn_momentum}, 
which is nonzero only in the particle domains. Similar to the free-swimming eel case, 
the particles' translational and rotational velocities are obtained
via Eqs.~\eqref{eq_conserve_linear_mom} and~\eqref{eq_conserve_angular_mom}.

The domain is discretized with a two-level hierarchy of Cartesian grids. The size of the coarsest grid is
$64 \times 256$ grid cells and $\nref = 4$ is taken for subsequent finer grids. Hence, the finest grid, with
spacing equivalent to that of a uniform mesh of size $256 \times 1024$, embeds each particle at all times.
The minimum grid spacing on the finest level is $\dx_\textrm{min} = \dy_\textrm{min} = 0.0390625D$.
A time step size of $\dt = 5\times 10^{-4}$ s is used. Particle $1$ is placed with initial center of mass
$\x_1 = (X_p, Y_p) = (-0.005D, 36D)$, while particle $2$ is placed below particle $1$ with initial 
center of mass location at $\x_2 = (X_p, Y_p) = (0, 34D)$. Under these conditions, the two particles 
start to accelerate downwards due to gravity. Particle $1$ travels through a low pressure wake 
created by the leading particle $2$, which causes particle $1$ to fall faster; this stage is called drafting.
Eventually, particle $1$ catches up to and nearly contacts particle $2$, a process termed as kissing in 
literature. This kissing stage is unstable and eventually the particles are left to tumble separately. 
The parameters here are chosen to match with previous numerical studies on drafting, kissing, 
and tumbling done by Feng et al.~\cite{Feng04}, Jafari et al.~\cite{Jafari11}, and Wang et al.~\cite{Wang14}.

Artificial repulsive forces are added to avoid numerical issues due to overlapping particles. The functional form
of this force on particle $i$ due to particle $j$ is given by

\begin{equation}
\label{eq_repulsive_force}
\F_{ij}^P = 
\begin{cases} 
      0,  & \|\x_i - \x_j \| > R_i + R_j + \zeta \textrm{ or } i = j\\
      \frac{c_{ij}}{\epsilon_P} \left(\frac{\|\x_i - \x_j \| - R_i - R_j - \zeta}{\zeta}\right)^2 
      \left(\frac{\x_i - \x_j}{\|\x_i - \x_j \| }\right),  & \|\x_i - \x_j \| \le R_i + R_j + \zeta
\end{cases}
\end{equation}
in which $R_i = R_j = R$ is the radius of both particles, $c_{ij} = \rho \pi R^2 g$ is a force scale parameter,
$\epsilon_P = 2.0$ g cm/s$^2$ is a stiffness parameters for collisions, and $\zeta = \dy_\textrm{min}$ is a
mesh threshold parameter indicating how far away the two particles need to be in order to feel particle-particle
interaction force. This particular repulsive force was used by Feng et al.~\cite{Feng04}. No repulsive force
between the particle and the wall are used for the case considered here. Similar to the gravitational force,
the particle interaction force is realized through an Eulerian body force $\f_P$ added to the right-hand side of
the momentum Eq.~\eqref{eqn_momentum}, which is nonzero only in the particle domains.

The strategy for placing the control volume for this example is different than the previous examples.
Rather than setting the CV in motion with some prescribed velocity to enclose the body, the CVs are 
chosen to surround each cylinder based on its center of mass  location: 
$[X_p - 5 \dx_\textrm{min}, X_p + 5 \dx_\textrm{min}] \times [Y_p - 5 \dy_\textrm{min}, Y_p +5 \dy_\textrm{min}]$. 
This means that when the particles are close to each other during the kissing stage, 
each CV contains the second body partially, leading to inaccurate force measurements.
This is a limitation of using simple rectangular control volumes. Fig.~\ref{fig_dkt_viz} shows 
the vortical structures generated by the two particles at four separate time instances, 
along with the locations of the moving control volumes.

\begin{figure}[H]
  \centering
  \subfigure[$t= 0.5$]{
    \includegraphics[scale = 0.5]{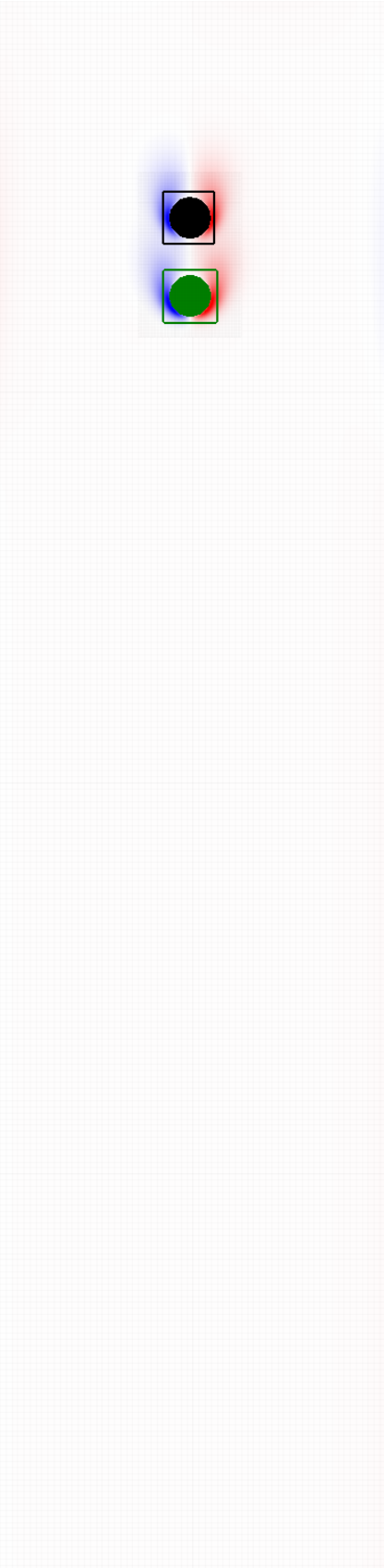}
    \label{fig_dkt_viz05}
  }
   \subfigure[$t = 1.5$]{
    \includegraphics[scale = 0.5]{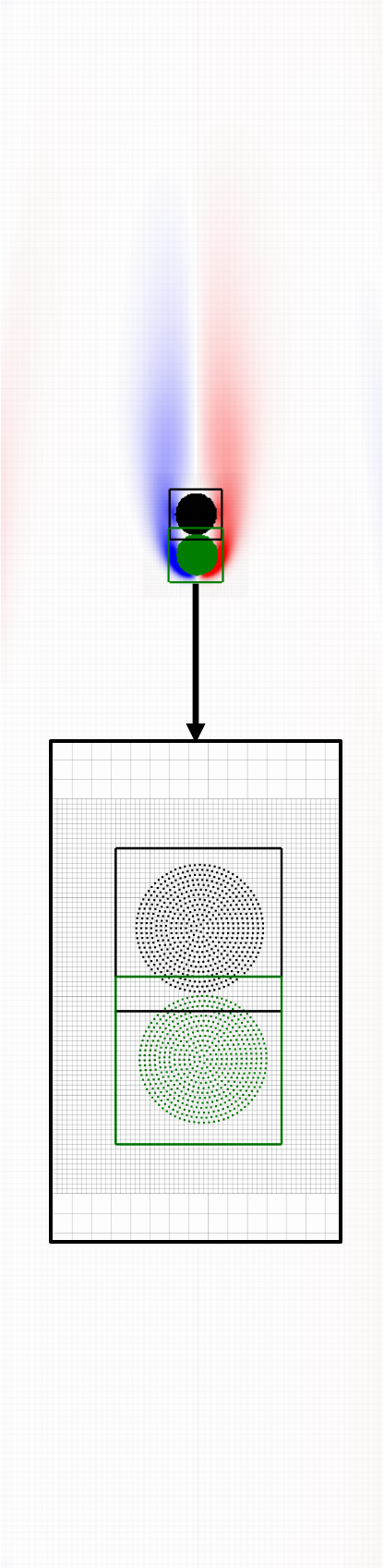}
    \label{fig_dkt_viz15}
  }
    \subfigure[$t = 2.0$]{
    \includegraphics[scale = 0.5]{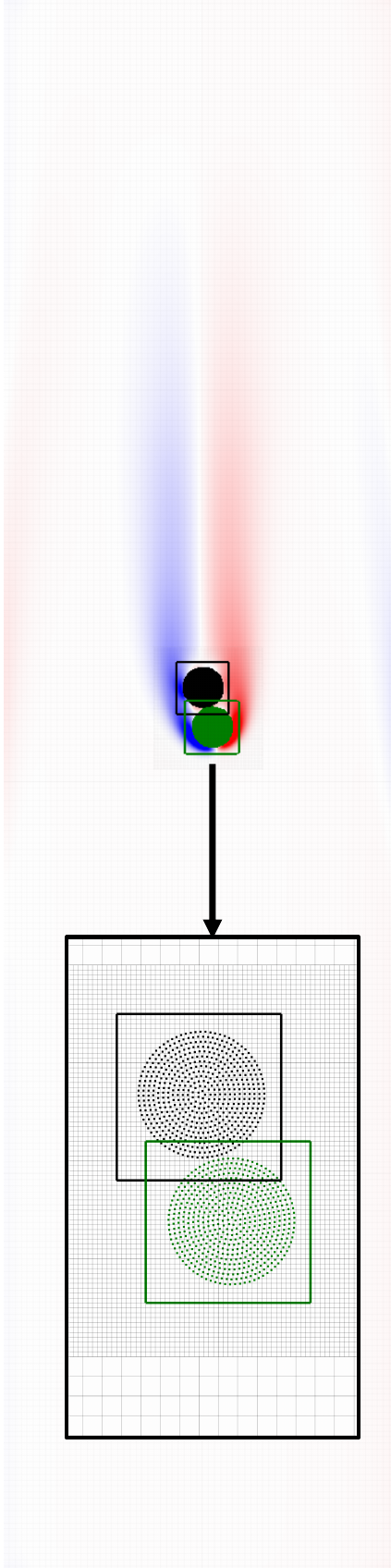}
    \label{fig_dkt_viz2}
  }
   \subfigure[$t = 3.0$]{
    \includegraphics[scale = 0.5]{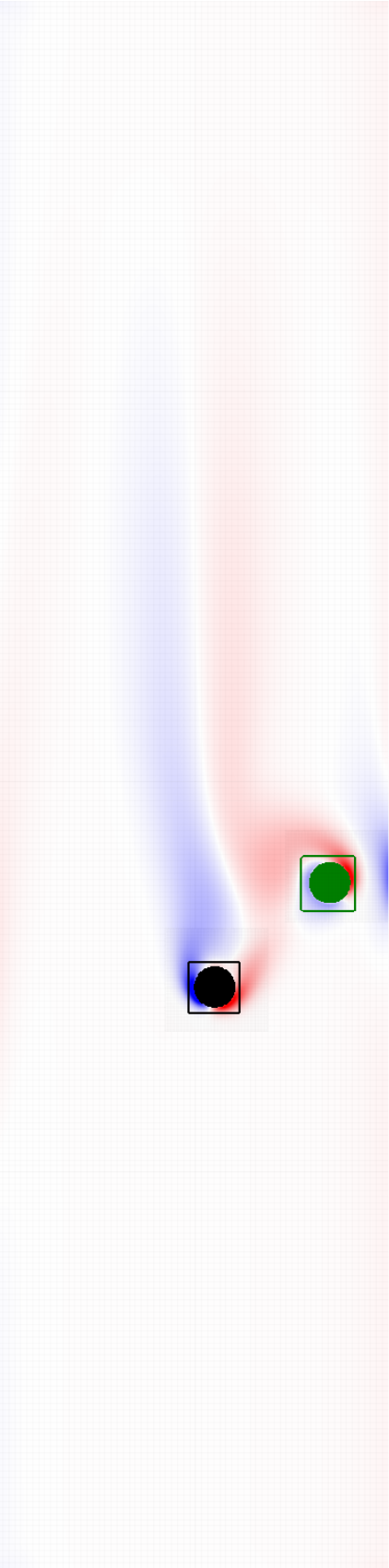}
    \label{fig_dkt_viz3}
  }
  \caption{\REVIEW{
  Vorticity generated by two-dimensional drafting, kissing, and tumbling of two cylindrical particles,
  along with the \emph{moving} CV location at four different time instances on an adaptive mesh hierarchy. 
  All figures are plotted for vorticity between $-30$ and $30$.}
   }
  \label{fig_dkt_viz}
\end{figure}

The expression for hydrodynamic force based on the Lagrange multiplier method needs to be modified to
account for the presence of additional body forces in the solid region (additional to the Lagrange multiplier 
constraint forces). These forces need to be added to the right-hand side of Eq.~\eqref{eq_extrinsic_force}, yielding
\begin{equation}
 \cF_i = \d{}{t} \int_{V_b^i(t)} \rho \; \U_i \dV - \int_{V_b^i(t)} (\F_i + \F_{i}^g + \F_{ij}^P) \dV, \label{eq_extrinsic_force_DKT}
\end{equation}
in which $\cF_i$, $\U_i$, $\F_i$, $\F_{i}^g$and $V_b^i(t)$ are the net hydrodynamic force, center of mass velocity, 
Lagrange multiplier force, gravitational force, and domain of particle $i$, respectively.
The expression for hydrodynamic force based on control volume analysis remains unchanged 
in the presence of the additional body forces in the solid region, and Eq.~\eqref{eq_OurForce} remains 
valid in this scenario\footnote{This is because $\int_{\Vcvt} \f \dV$ that contains contribution of all body forces 
in Lagrangian domain is evaluated via the right-hand side of Eq.~\eqref{eq_sumLag_CV}}.  

Figs.~\ref{fig_dkt_xcom} and~\ref{fig_dkt_ycom} show the time evolution of the centers of mass $(X_p, Y_p)$ of
both particles. Particle $1$ gradually approaches particle $2$ up until near $t = 2.5$ s when the particles kiss and
eventually separate. Upon separation, both particles over time reach a terminal velocity. This temporal behavior matches 
well with the results of Jafari et al.~\cite{Jafari11}. Figs.~\ref{fig_dkt_Fx} and~\ref{fig_dkt_Fy} show the time evolution
of net axial ($\mathcal{F}_x = \cF \cdot \e_x$) and lateral ($\mathcal{F}_y = \cF \cdot \e_y$) forces acting on each particle,
calculated using both the control volume and Lagrange multiplier approaches. Between $t = 1$ s and $t = 2.5$ s, the
particles are close to each other and it is not possible to create rectangular CVs that contain only a single 
particle (see insets in Figs.~\ref{fig_dkt_viz15} and ~\ref{fig_dkt_viz2}). Hence, the CV force 
calculations during this time period are inaccurate. However, the forces calculated by the LM method 
remain accurate at all times. Outside of this time period, both the CV and LM approaches are in 
excellent agreement. Eventually, the net hydrodynamic force on each particle is $-\F^g$, indicating a 
terminal velocity has been achieved.

\begin{figure}[H]
  \centering
 \subfigure[]{
    \includegraphics[scale = 0.32]{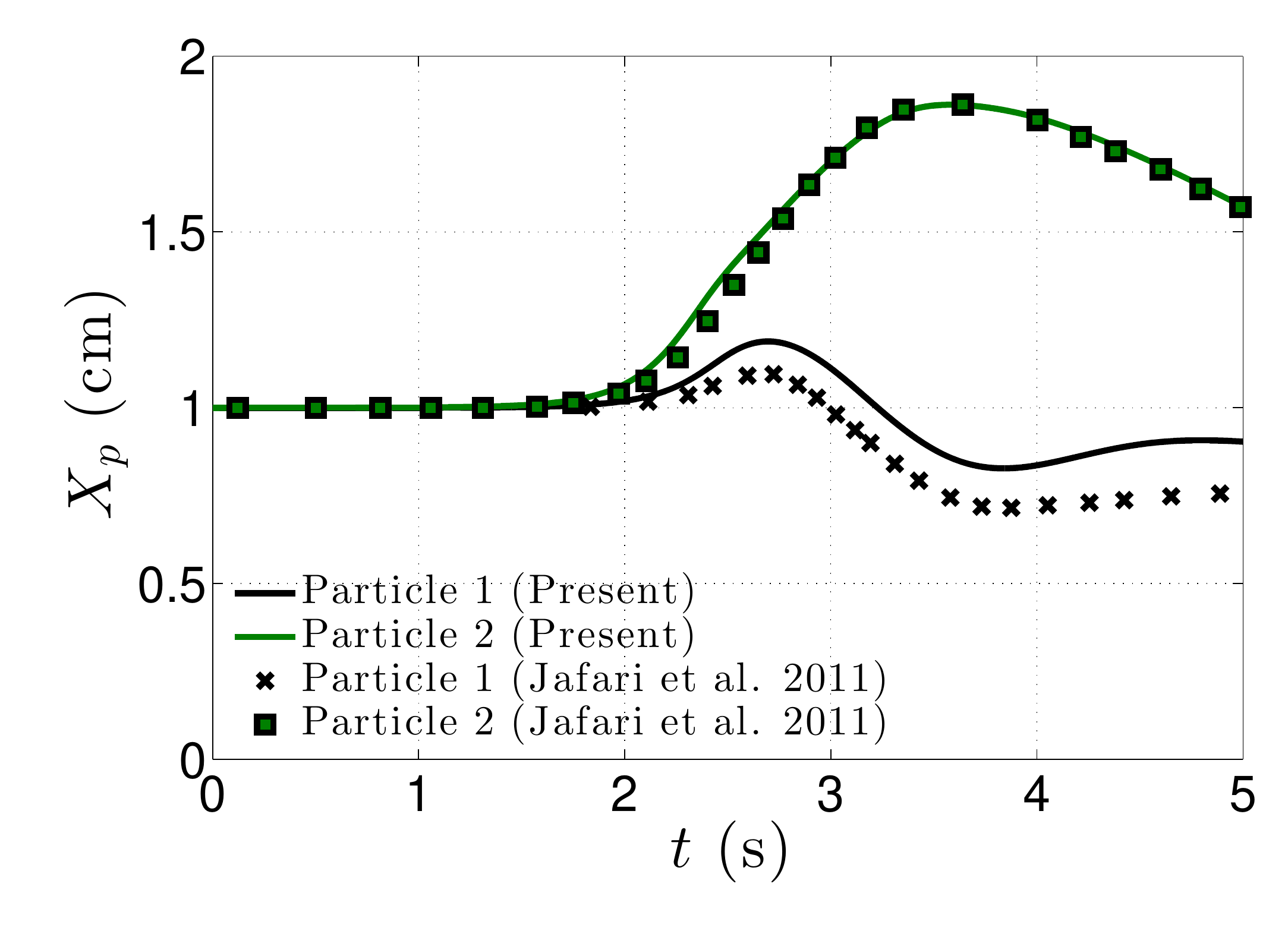}
    \label{fig_dkt_xcom}
  }
  \subfigure[]{
    \includegraphics[scale = 0.32]{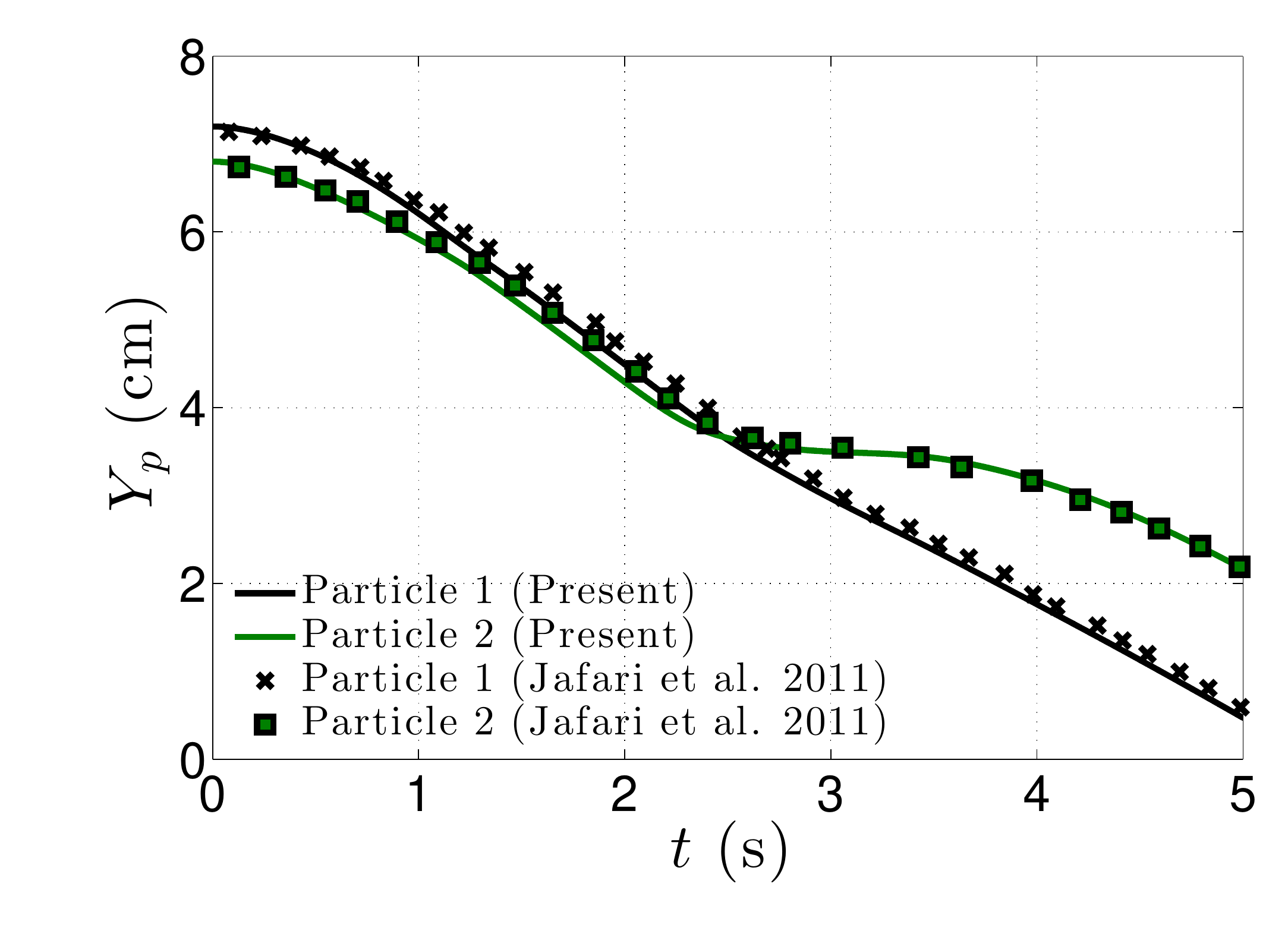}
    \label{fig_dkt_ycom}
  }
   \subfigure[]{
    \includegraphics[scale = 0.32]{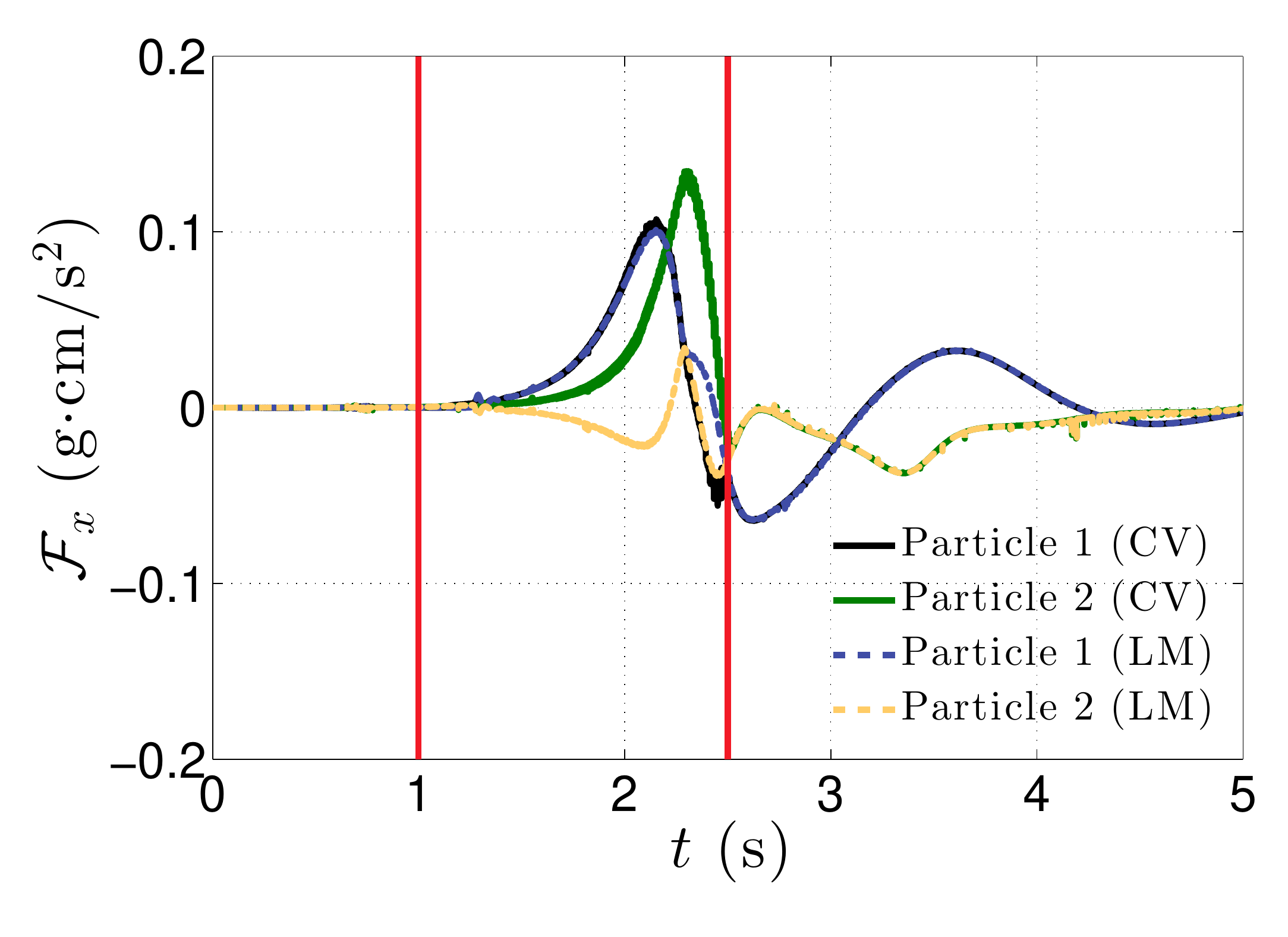}
    \label{fig_dkt_Fx}
  }
   \subfigure[]{
    \includegraphics[scale = 0.32]{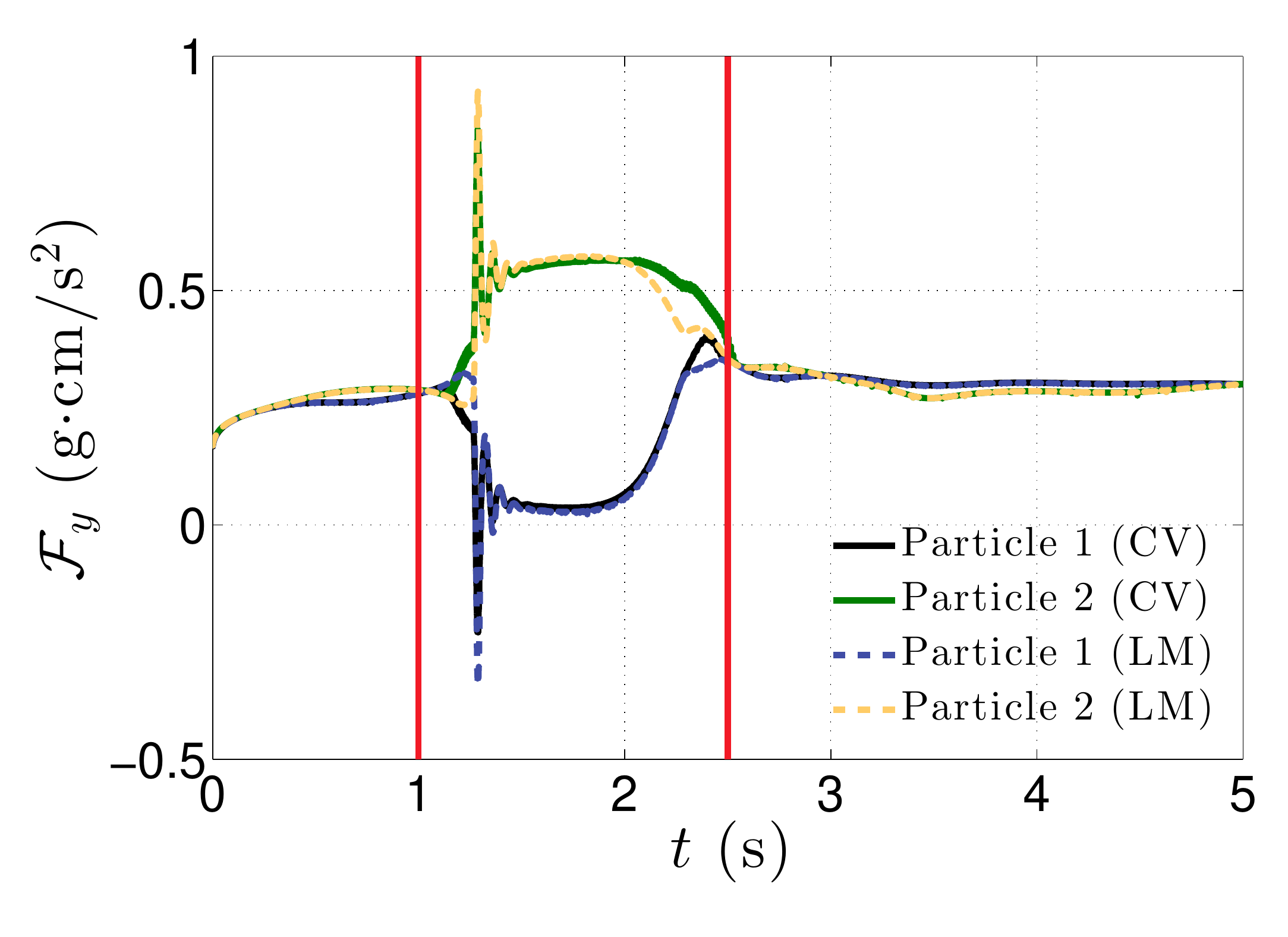}
    \label{fig_dkt_Fy}
  }
  \caption{ \REVIEW{
  \subref{fig_dkt_xcom}
  Temporal evolution of the $x$ coordinate of center of mass $X_p$ for (---, black) particle $1$ \& (---, green) particle $2$.
  \subref{fig_dkt_ycom}
   Temporal evolution of the $y$ coordinate of center of mass $Y_p$ for (---, black) particle $1$ \& (---, green) particle $2$.
   Center of mass data from Jafari et al. \cite{Jafari11} for ($\times$, black) particle $1$ \& ($\blacksquare$, green) particle $2$.
  \subref{fig_dkt_Fx}
   Temporal evolution of net $\mathcal{F}_x$ via the moving control volume approach on (---, black) particle $1$  \& (---, green) 
   particle $2$. Temporal evolution of net $\mathcal{F}_x$ via the Lagrange multiplier approach approach on (\texttt{---}, blue) 
   particle $1$  \& (\texttt{---}, yellow) particle $2$.
   \subref{fig_dkt_Fy}
   Temporal evolution of net $\mathcal{F}_y$ via the moving control volume approach on (---, black) particle $1$  \& (---, green) 
   particle $2$. Temporal evolution of net $\mathcal{F}_y$ via the Lagrange multiplier approach approach on (\texttt{---}, blue) 
   particle $1$  \& (\texttt{---}, yellow) particle $2$. The time period between the two red lines indicates inaccurate 
   CV force measurements since each rectangular CV contains multiple bodies.}
   }
  \label{fig_dkt}
\end{figure}

} 

\subsection {Stokes flow}
In the previous sections, we considered finite Reynolds number cases simulated 
using a direct forcing IB method in which Lagrange multipliers were approximated 
in the body domain. Here we consider a fully constrained IB method in which we 
compute Lagrange multipliers exactly. We consider Stokes flow examples here,
although the fully constrained method also work equally well at finite Reynolds numbers
as shown in~\cite{Kallemov16}.

For steady Stokes flow in the absence of inertia ($\rho = 0$), the momentum equation reads as

\begin{equation}
\label{eq_stokes}
-\grad p(\x) + \mu \lap \u(\x) +\f(\x) =  \V 0.
\end{equation}

\noindent Since Eq.~\eqref{eq_stokes} is a steady state problem, it cannot
be solved numerically with the split fluid-structure solver described in~\cite{Bhalla13}.
Rather, the discretized system of Eqs.~\eqref{eqn_stokes_momentum}-\eqref{eqn_stokes_constraint} 
are solved by the monolithic fluid-structure solver described in~\cite{Kallemov16} to obtain 
a numerical solution to the constrained Stokes system. The rigidity constraint is
enforced on the surface of the body, and the body is discretized
only by surface nodes and not by a volumetric mesh. This is because enforcing the rigidity 
constraint on the surface also imposes rigid body motion inside the body for 
Stokes flow. Setting the inertial terms to zero in Eqs.~$\eqref{eq_OurForce}$ and
$\eqref{eq_OurTorque}$ yield

\begin{align}
\cF &=  \oint_{\Scv} \ndot \left[-p \I + \T \right] \dS, \label{eq_cv_stokes_force} \\ 
\cM  &= \oint_{\Scv}  [\rcross (- p \; \n + \ndot \T) \;]  \dS. \label{eq_cv_stokes_torque}
\end{align}

\noindent Hence, the net hydrodynamic force and torque on an object in Stokes
flow is simply the pressure and viscous fluxes through the control surface. For 
the LM method, the Lagrange multiplier force density $\F(\X)$ is computed (exactly) on  
body surface $\X \in \Sb$, and the net force and torque is given by

\begin{align}
\cF &=  -\oint_{\Sb} \F(\X) \dS, \label{eq_lag_stokes_force}\\ 
\cM &= -\oint_{\Sb} \R \wedge \F(\X) \dS. \label{eq_lag_stokes_torque}
\end{align}

\subsubsection{Flow between two concentric shells}
We first consider the case of two concentric shells, which was studied numerically
by Kallemov et al.~\cite{Kallemov16}. The inner and outer shell have geometric
radii $R_1^g = 1.807885$ and $R_2^g = 4 R_1^g$, respectively. The computational
domain is a cube of size $[0,L]^3 = [0,4.15 R_2^g]^3$, which is discretized by a 
uniform grid of size $60^3$. The center of both shells is placed at $(x,y,z) = (L/2,L/2,L/2)$.
Uniform velocity $\u = (1,0,0)$ is prescribed on each wall of the computational domain, and the inner and
outer shell are set to have rigid--body velocity $(U_1,V_1,W_1) = (0,0,0)$ and 
$(U_2,V_2,W_2) = (1,0,0)$, respectively. The inner shell is discretized with $42$ surface markers
while the outer shell is discretized with $642$ surface markers to ensure that the markers
are about $2$ grid cells apart. The viscosity is set to $\mu = 1$.

Each spherical shell has an effective \emph{hydrodynamic} radius $R^h$ due to the immersed
boundary kernel used to discretize the delta-function~\cite{Vazquez-Quesada14}.
For the 6-point kernel considered by Kallemov et al., it was found that $R_1^h = 1.22 R_1^g$
and $R_2^h = 0.96R_2^g$ for the numerical parameters chosen here \cite{Kallemov16}. It was
also found that as both the Eulerian and Lagrangian meshes are refined, $R^h/R^g \rightarrow
1$. These hydrodynamic radii can be used in the analytical expression for the drag on the 
inner sphere \cite{Happel65}, given by

\begin{equation}
\label{eq_analytical_concentric_sphere}
\cF_\textrm{exact} \cdot \e_x = -6 \pi \mu R_1^h U_2 K,
\end{equation} 

\noindent in which $K = (1-\lambda^5)/\alpha$, $\alpha = 1-9\lambda/4 + 5 \lambda^3/2
-9\lambda^5/4 + \lambda^6$ and $\lambda = R_1^h/R_2^h$.

\begin{figure}[H]
  \centering
    \includegraphics[scale = 0.3]{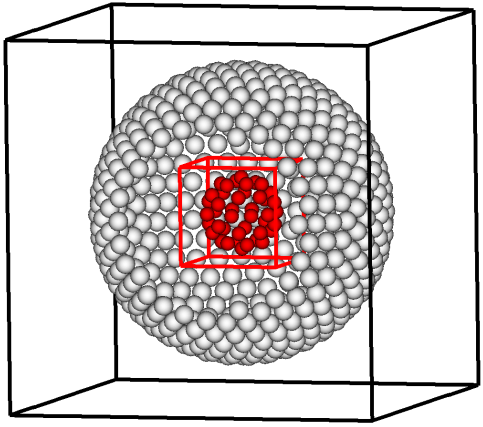}
  \caption{Surface nodes on the inner shell (red) and outer shell (gray)
  for Stokes flow between concentric spheres.
  Location of a control volume containing only the inner shell (red) and
  containing both the inner and outer shell (black).}
  \label{fig_concentric_stokes_CV}
\end{figure}

Depending on the control volume size, the drag on either the inner shell, or on both
shells can be obtained. First, a CV of dimension $[L/2-1.659R_1^g,L/2+1.659R_1^g]^3$
is chosen to surround the inner shell, but to exclude the outer shell. Next, a CV of dimension
$[L/2-1.383R_2^g,L/2+1.383R_2^g]$ is chosen to include both inner and outer shells. Fig.~\ref{fig_concentric_stokes_CV} shows the configuration of the concentric shells/spheres and the
two control volumes.

Table~\ref{table_F_inner_compare} shows the drag measurements from the analytical
expression, from integrating surface Lagrange multipliers, and from the control volume
analysis. The middle column shows that all three methods are in agreement for the
drag on the inner shell. Moreover, the last column shows that the combined drag
on both inner and outer shells is the same when computed from Lagrange multipliers
and control volume analysis.

\begin{table}[H]{}
  \caption{Comparison of force measurements
  for Stokes flow around two concentric spheres.}
  \centering
  \begin{tabular}{lll}
    \hline \hline
    Method &  $|\cF \cdot \e_x|$ (Inner) & $|\cF \cdot \e_x|$ (Inner \& Outer)\\ [0.5ex]
    \hline
    Analytical Eq.~\eqref{eq_analytical_concentric_sphere}  &  115.409 & N/A \\
    Lagrange multiplier  Eq.~\eqref{eq_lag_stokes_force}     &  114.297 & 17.350 \\
    Control volume Eq.~\eqref{eq_cv_stokes_force}              &  114.298 & 17.350  \\ [1ex]
    \hline
  \end{tabular}
  \label{table_F_inner_compare}
\end{table}

\subsubsection{Single rotating shell}
Next, we consider a single shell rotating in a bounded domain with no exterior flow.
The shell is taken to be the same as the inner sphere in the previous example, with
geometric radius $R_1^g = 1.807885$ in a computational domain of size $[0,L]^3 = 
[0, 16.6R_1^g]^3$. There is no outer shell in this example. The center of the shell is 
placed at the centroid of the cube $(x,y,z) = (L/2,L/2,L/2)$ and $\u = (0,0,0)$ is set at 
all computational boundaries. Periodic boundary conditions were also used and yielded nearly identical
torque measurements (data not shown). The viscosity is set to $\mu = 1$. The shell rotates about a 
diameter with angular velocity $\vomega = (0,1,0)$. In an \emph{unbounded} flow at rest, 
Fax\'{e}n's law states that the torque on the sphere by the fluid is

\begin{equation}
\cM = -8 \pi \mu (R_1^h)^3 \vomega, \label{eq_FaxenTorque}
\end{equation}
in which the hydrodynamic radius of the sphere $R_1^h$ is used~\cite{Happel65}. Although the domain
in the numerical method is bounded, we still get decent agreement between our numerical results
and Eq.~\eqref{eq_FaxenTorque}.

A CV of dimension $[L/2-1.659R_1^g, L/2+1.659R_1^g]^3$ is chosen to surround the shell.
Three different grid sizes are used to discretize the domain: $60^3$, $120^3$, and $240^3$, which
corresponds to $42$, $162$, and $642$ surface markers on the shell respectively (to ensure that 
the markers are approximately 2 grid cells apart).

Table~\ref{table_T_compare} shows the torque measurements $\cM \cdot \e_y$ from 
the analytical expression, from integrating the moments of surface Lagrange multipliers, and from 
the control volume analysis, for the three different grid resolutions. As expected, the LM and 
CV measured torque do not match exactly with the analytical expression (presumably because 
of finite domain effects). However, the LM and CV torque values are in excellent agreement
with each other. 

\begin{table}[H]{}
  \caption{Comparison of torque measurements
  $\cM \cdot \e_y$ for Stokes flow around single rotating sphere.}
  \centering
  \begin{tabular}{llllll}
    \hline \hline
    Grid size &  Number of markers & $R_1^h/R_1^g$~\cite{Kallemov16} & 
    Analytical Eq.~\eqref{eq_FaxenTorque} & LM Eq.~\eqref{eq_lag_stokes_torque} &
    CV Eq.~\eqref{eq_cv_stokes_torque} \\ [0.5ex]
    \hline
    $60^3$   &  42   & 1.22 & -269.669 & -236.474 & -236.474 \\
    $120^3$ &  162 & 1.09 & -192.323 & -183.662 & -183.662 \\
    $240^3$ &  642 & 1.04 & -167.052 & -164.446 & -164.446 \\[1ex]
    \hline
  \end{tabular}
  \label{table_T_compare}
\end{table}


\section*{Conclusions}
In the present study we presented a moving control volume (CV) approach to compute 
the net hydrodynamic forces and torques on a moving body immersed in a fluid.
This approach does not require evaluation of (possibly) discontinuous spatial velocity
or pressure gradients within or on the surface of the immersed body.
The analytical expressions for forces and torques were modified from those initially
presented in \cite{Noca97}, and this modification has been shown to eliminate spurious
jumps in drag \cite{Bergmann11}. Our implementation treats the control volume
as a rectangular box whose boundary is forced to remain on grid lines, which greatly
simplifies the evaluation of surface integrals.

The approach is shown to accurately compute the forces and torques on a wide array
of fluid-structure interaction problems, including flow past stationary and moving objects, 
Stokes flow, and high Reynolds number free-swimming. Spurious momentum gain or loss
due to adaptive mesh refinement can produce jumps in the computed forces in the CV
approach, although forcing the CV to remain on the finest grid level can ameliorate this issue. 

We also show the equivalence between the Lagrange multiplier (LM) approach and the CV approach.
The main advantage of the CV approach over the LM approach is that it is applicable
to situations where explicit Lagrange multipliers are not available, for example, in the embedded
boundary/cut-cell approach to FSI. 

The control volume approach implemented here assumes a no-slip boundary condition on the fluid-structure
interface. However, a generalization for transpiration boundary conditions can be derived as well (see
Eq.~\eqref{eq_our_general_force}). Use of such a boundary condition would required richer geometric 
information and data structures to evaluate surface quantities of the immersed body. Finally, our approach 
can be easily extended to cases where additional body forces are present in the momentum equation.

\section*{Acknowledgements}
A.P.S.B and N.N acknowledge helpful discussions related to software design in IBAMR 
with Boyce E. Griffith (UNC-Chapel Hill) over the course of this work.  
N.N, N.A.P, and A.P.S.B~acknowledge computational resources
provided by Northwestern University's Quest high performance computing
service.  N.N~acknowledges research support from the National Science Foundation 
Graduate Research Fellowship Program (NSF award DGE-1324585).
N.A.P~acknowledges support from the National Science Foundation (NSF award SI2-SSI-1450374).
A.P.S.B and H.J acknowledge support from the U.S. Department of Energy, Office of Science, 
ASCR (award number DE-AC02-05CH11231). A part of this work was carried at UNC-Chapel Hill for 
which A.P.S.B gratefully acknowledges support from awards NIH HL117163 and NSF ACI 1450327 
(awarded to Boyce E. Griffith). 

\appendix 

\section{Derivation of the new hydrodynamic force expression} \label{app_force}

Here we present a detailed derivation of Eq.~\eqref{eq_OurForce} from  
Eq.~\eqref{eq_NocaForce}. The ultimate goal is to obtain an expression
for hydrodynamic force which involves integral contributions from a 
single CV rather than two time-lagged CVs. The Reynolds transport 
theorem (RTT)~\cite{Kundu14,Pozrikidis11} gives an expression for the time 
derivative of an arbitrary quantity $\q(\x,t)$ on a time dependent region $\Omegat$

\begin{equation}
\label{eq_RTT}
\d{}{t} \int_{\Omegat} \q \dV =
\int_{\Omegat} \D{\q}{t} \dV
+ \oint_{\partial \Omegat} (\ndot \uS) \q \dS,
\end{equation}

\noindent in which $\uS$ is the velocity and $\n$ is the outward pointing
unit normal vector of the boundary $\partial \Omegat$.
Applying the RTT to Eq.~\eqref{eq_NocaForce} yields the expression

\begin{align}
\cF(t) &= 
-\int_{\Vt} \rho \D{\u}{t} \dV 
- \oint_{\partial \Vt} (\ndot \uS) \rho \u \dS  \nonumber \\
& + \oint_{\Scvt} \ndot \left[-p \I - (\u-\uS)\rho \u + \T \right] \dS 
- \oint_{\Sbt} \ndot (\u-\uS)\rho \u \dS.
\end{align}

\noindent Recall that $\Vt = \Vcvt \setminus \Vbt$, where $\Vcvt$ is the entire
control volume and contains the body domain $\Vbt$; $\Scvt = \partial \Vcvt$ is the
boundary of the CV, and $\Sbt = \partial \Vbt$. Hence, $\partial \Vt = \Scvt \cup \Sbt$.
and the second integral can be split into two boundary integrals

\begin{align}
\cF(t) = \nonumber
& -\int_{\Vt} \rho \D{\u}{t} \dV \\
&- \oint_{\Scvt} (\ndot \uS) \rho \u \dS 
+ \oint_{\Scvt} \ndot \left[-p \I - (\u-\uS)\rho \u + \T \right] \dS  \\
& - \oint_{\Sbt} (\ndot \uS) \rho \u \dS \nonumber
- \oint_{\Sbt} \ndot (\u-\uS)\rho \u \dS,
\end{align}

\noindent which can be simplified to obtain

\begin{equation}
\label{eq_simplified_drag}
\cF(t) = 
 -\int_{\Vt} \rho \D{\u}{t} \dV
 + \oint_{\Scvt} \ndot \left[-p \I - \u \rho \u + \T \right] \dS 
 - \oint_{\Sbt} \ndot (\u \rho \u) \dS.
\end{equation}

\noindent Recall that in Eq.~\eqref{eq_simplified_drag}, the unit normal vector
points outward on $\Scvt$ and inward on $\Sbt$ since the integral is considered with
respect to the boundary of $\Vt$. Next, notice that by definition
$\Vcvt = \Vt \cup \Vbt$ is the union of disjoint regions. Hence, the integral over
$\Vt$ can be split into $\int_{\Vt} \q \dV = \left(\int_{\Vcvt} - \int_{\Vbt}\right) \q \dV$.
Applying the split to Eq.~\eqref{eq_simplified_drag} yields

\begin{equation}
\cF(t) =  
 -\int_{\Vcvt} \rho \D{\u}{t} \dV 
 + \int_{\Vbt} \rho \D{\u}{t} \dV
 + \oint_{\Scvt} \ndot \left[-p \I - \u \rho \u + \T \right] \dS 
 - \oint_{\Sbt} \ndot (\u \rho \u) \dS.
\end{equation}

\noindent Finally, we can apply the Reynolds transport theorem (Eq.
~\eqref{eq_RTT}) to the integral over $\Vbt$ above. Letting $\N$ be the
\emph{outward} pointing unit normal vector to $\Sbt$, we obtain

\begin{align}
\cF(t) = 
 &-\int_{\Vcvt} \rho \D{\u}{t} \dV  \nonumber
 + \d{}{t} \int_{\Vbt} \rho \u \dV \\
 &+ \oint_{\Scvt} \ndot \left[-p \I - \u \rho \u + \T \right] \dS \\
 &  - \oint_{\Sbt} (\Ndot \uS) \rho \u \dS 
 - \oint_{\Sbt} \ndot (\u \rho \u) \dS. \nonumber
\end{align}

\noindent Substituting the fact that $\N = -\n$ yields a general expression
for the hydrodynamic force on an immersed body

\begin{equation}
\label{eq_our_general_force}
\cF(t) = 
 -\int_{\Vcvt} \rho \D{\u}{t} \dV  
 + \d{}{t} \int_{\Vbt} \rho \u \dV 
 + \oint_{\Scvt} \ndot \left[-p \I - \u \rho \u + \T \right] \dS 
 - \oint_{\Sbt} \ndot (\u \rho \u - \uS \rho \u) \dS.
\end{equation}

\noindent Eq.~\eqref{eq_our_general_force} can be further simplified if we assume
no-slip boundary conditions at the fluid-structure interface by setting $\uS = \u$, which gives
 
\begin{equation}
\cF(t) = 
-\int_{\Vcvt} \rho \D{\u}{t} \dV
+ \d{}{t} \int_{\Vbt} \rho \u \dV 
+ \oint_{\Scvt} \ndot \left[-p \I - \u\rho \u + \T \right] \dS.
\end{equation}

\noindent This is the net hydrodynamic force expression as written in Eq.~\eqref{eq_OurForce}.


\section{Derivation of the new hydrodynamic torque expression} \label{app_torque}
Conservation of angular momentum for a material volume $\Vmt$ (a volume 
that moves with the local fluid velocity) that shares its boundary with 
an arbitrary moving control volume $\Vt$ at time $t$ can be written as~\cite{Kundu14,Pozrikidis11}

\begin{align}
	\frac{d}{dt} \int_{\Vmt} (\rcross \rho \u) \dV &= \oint_{\partial \Vmt = \partial \Vt}   \rcross [ \ndot \Sigma] \dS \nonumber \\ 
	 &=  \oint_{\Scvt} \rcross [ \ndot \Sigma ] \dS + \oint_{\Sbt} \rcross [ \ndot \Sigma ] \dS,
\end{align}
in which $\r = \x - \x_{0}$, with $\x_{0}$ as a reference point for computing torques,
and $\Sigma = -p\I + \T$.
Letting $\cM(t)$ be the torque exerted by the fluid on the body and noticing that 
$\oint_{\Sbt}  \rcross [ \ndot \Sigma] \dS = -\cM(t)$, we have

\begin{equation}
	\cM(t) = -\frac{d}{dt} \int_{\Vmt} (\rcross \rho \u) \dV + \oint_{\Scvt}  \rcross [ \ndot \Sigma ] \dS. 
\end{equation}

\noindent Using the RTT, the integral of an arbitrary quantity $\Phi$ over the material volume 
$\Vmt$ can be related to integral over arbitrary volume $\Vt$ with surface velocity moving with $\uS$ as
\begin{equation}
 \frac{d}{dt} \int_{\Vmt} \Phi \dV =\frac{d}{dt} \int_{\Vt} \Phi \dV 
 + \oint_{\partial \Vmt = \partial \Vt} \ndot (\u-\uS) \Phi \dS. \label{eq_vm_to_v}
 \end{equation} 

\noindent Using Eq.~\eqref{eq_vm_to_v} with $\Phi = \rcross \rho \u$, the expression 
for torque becomes

\begin{align}
	\cM(t) &= -\frac{d}{dt} \int_{\Vt} (\rcross \rho \u) \dV 
	 + \oint_{\Scvt} [\rcross (- p \; \n +  \ndot \T) - \ndot (\u-\uS) (\rcross \rho \u)  ] \dS \nonumber \\
	 & - \oint_{\Sbt} \ndot (\u - \uS)(\rcross \rho \u) \dS. \label{eq:Torque}
\end{align}
Finally, by manipulating the term derivative term in 
Eq.~\eqref{eq:Torque} using the RTT we get an expression for torque on an immersed 
body as 

\begin{equation}
\cM(t) = -\int_{\Vcvt} \rho \rcross \frac{\partial \u}{\partial t} \dV 
+ \frac{d}{dt} \int_{\Vbt} \rho (\rcross \u) \dV 
+ \oint_{\Scvt}  [\rcross (- p \; \n + \ndot \T) - (\ndot \u)\rho (\rcross \u) \;]  \dS.
\end{equation}


\section{Numerical discretization} \label{app_numerical_integration}

Here we describe the discrete evaluation of Eqs.~\eqref{eq_OurForce} 
and~\eqref{eq_OurTorque} to obtain the net hydrodynamic force and torque
on an immersed body.  For notational simplicity, we present the discretized
equations in two spatial dimensions. An extension to three spatial dimensions 
is straightforward. A discrete grid covers the physical domain 
$\Omega$ with mesh spacing $\dx$ and $\dy$ in each direction.
The position of each grid cell center is given by $\x_{i,j} = \left(x_{i,j},y_{i,j}\right)$. 
For a given cell center, $\x_{i-\half,j}$ denotes the physical location of the cell face 
that is half a grid space away from $\x_{i,j}$ in the negative $x$-direction, i.e. 
$\x_{i-\half,j} = \left(x_{i,j} - \frac{\dx}{2},y_{i,j}\right)$. Similarly $\x_{i,j-\half}$ denotes the
physical location of the cell face that is half a grid cell away from  $\x_{i,j}$ 
in the negative $y$-direction, i.e. $\x_{i,j-\half} = \left(x_{i,j},y_{i,j}- \frac{\dy}{2}\right)$.
The discrete approximations described here are also valid when adaptive mesh
refinement is used, although volume weights need to be appropriately
modified for different velocity components. 

Let $t^n$ be the time at time step $n$. After stepping forward from time $t^n$
to $t^{n+1} = t^n + \dt$, a pressure solution is obtained at cell centers
$p_{i,j}^{n+1} = p\left(\x_{i,j},t^{n+1}\right)$, while velocity components are obtained
at cell faces: $u_{i-\half,j}^{n+1} = u\left(\x_{i-\half,j}, t^{n+1}\right)$ and 
$v_{i,j-\half}^{n+1} = v\left(\x_{i,j-\half}, t^{n+1}\right)$. These are the only Eulerian
quantities needed to evaluate the discrete approximations to Eqs.~\eqref{eq_OurForce} 
and~\eqref{eq_OurTorque}. Only rectangular control volumes are considered in the 
present work. A CV is described by its lower left and upper right corners: letting $(\xL,\yL)$ and $(\xU,\yU)$
denote the lower and upper corners respectively, the control volume is defined to be 
the Cartesian product of intervals $\Vcvt = \left\{\x \in \Omega \mid \x \in [\xL,\xU] \times [\yL,\yU] \right\}$. 
Moreover, $\Scvt$ is forced to remain on grid lines and it not allowed to cross into the interior
of grid cells. This greatly simplifies the required numerical approximations. 
Refer to Fig.~\ref{fig_CV_discrete} for a sketch of the control volume configuration 
over a staggered mesh discretization. Let the control volume and surface at a time 
instance $t^{n+1}$ be denoted by $\Vcv^{n+1} = \Vcv(t^{n+1})$ 
and $\Scv^{n+1} = \Scv(t^{n+1})$, respectively.

\subsection{Discrete approximation to surface integrals}
The control surface $\Scv^{n+1}$ is composed of four segments (eight
faces in $3$D) denoted by $\cB$, $\cL$, $\cT$, and $\cR$ in Fig.~\ref{fig_CV_discrete}. 
Consequently, computing surface normals on each of these segments is simple, 
e.g for the bottom segment $\cB^{n+1}$, $\n = -\e_y$. The discretized surface 
integral of a quantity $\vec{\Phi}$ over $\Scv^{n+1}$ is simply the sum over these 
four segments

\begin{equation}
\label{eq_discrete_oint_CV}
\oint_{\Scv^{n+1}} \ndot \vec{\Phi} \dS =
\oint_{\cR^{n+1}} \e_x \cdot \vec{\Phi} \dS
- \oint_{\cL^{n+1}} \e_x \cdot \vec{\Phi} \dS
+ \oint_{\cT^{n+1}} \e_y \cdot \vec{\Phi} \dS
- \oint_{\cB^{n+1}} \e_y \cdot \vec{\Phi} \dS.
\end{equation}

Moreover, it is sufficient to show the discrete approximation to the surface
integral over a single segment since the contribution from the other three
segments are computed analogously. Over the bottom surface $\cB$, the
discretization of each term is given by

\begin{align}
\oint_{\cB} \ndot \left(-p \I\right) \dS  & = \oint_{\cB} - \e_y \left(-p\right)\dS  \nonumber \\
		  & \approx \sum_{(i,j-\half)\in \cB}  -\e_y \frac{-(p_{i,j} + p_{i,j-1})}{2} \dx,\label{eq_discrete_oint_B_p}  \\
\oint_{\cB} \ndot \left(- \u\rho \u\right) \dS   & = -\rho \oint_{\cB} -v \left(u \e_x + v \e_y \right) \dS  \nonumber \\  
								& \approx -\rho \sum_{(i,j-\half)\in \cB} -v_{i,j-\half}\left[\frac{u_{i-\half,j}+u_{i+\half,j}+u_{i-\half,j-1}+u_{i+\half,j-1}}{4}\e_x+v_{i,j-\half}\e_y\right] \dx  \label{eq_discrete_oint_B_u}, \\
 \oint_{\cB} \ndot \mu\left(\grad \u + \grad \u^T\right) \dS  & = \mu \oint_{\cB}  -\left[\left(\D{u}{y} + \D{v}{x} \right)\e_x + 2 \D{v}{y}\e_y \right] \dS  \nonumber \\
				& \approx \mu \sum_{(i,j-\half)\in \cB}  - \left[\left(\frac{u_{i+\half,j} - u_{i+\half,j-1} + u_{i-\half,j} - u_{i-\half,j-1}}{2\dy} + \frac{v_{i+1,j-\half} - v_{i-1,j-\half} }{2 \dx}\right)\e_x \right] \dx \nonumber \\
& + \mu \sum_{(i,j-\half)\in \cB} - \left[ 2 \frac{v_{i,j+\half} - v_{i,j-\3half}}{2 \dy}\e_y \right]\dx\label{eq_discrete_oint_B_visc}.
\end{align}

\begin{figure}[H]
  \centering
    \includegraphics[scale = 0.4]{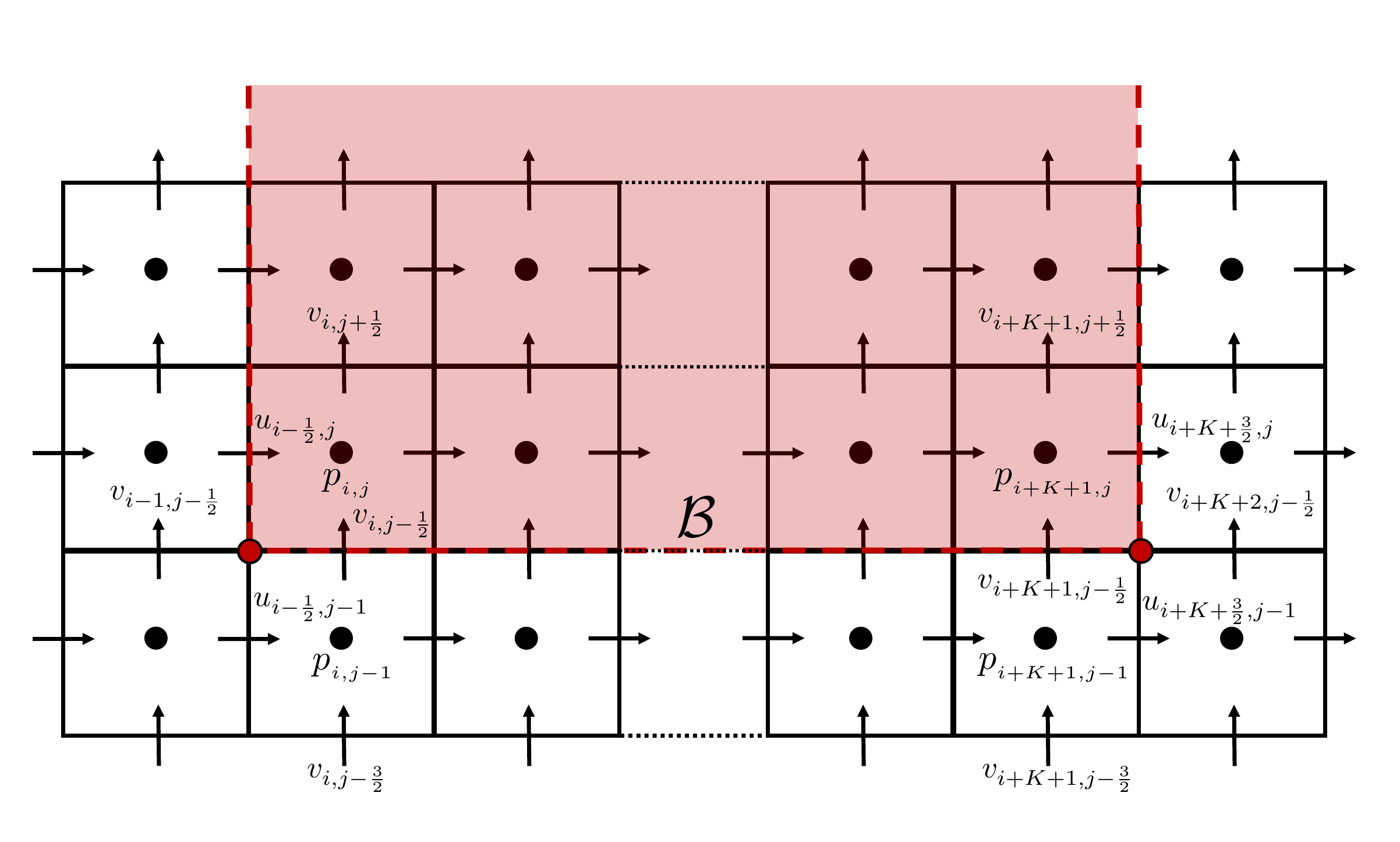}
  \caption{
  The velocity and pressure locations on the staggered--grid
  discretization required to evaluate the surface integral terms
  in the hydrodynamic force calculation on the bottom face of $\Scv$.}
  \label{fig_CV_bottom_face}
\end{figure}

\noindent Fig.~\ref{fig_CV_bottom_face} shows a schematic of the pressure
and velocity values required to evaluate Eqs.~\eqref{eq_discrete_oint_B_p},
\eqref{eq_discrete_oint_B_u}, and~\eqref{eq_discrete_oint_B_visc}.
Evaluating the surface integral on $\cB$ in the torque calculation about a point
$\x_0$ is done by computing $\r_{i,j-\half} = \x_{i,j-\half} - \x_0$ and evaluating
the cross product between $\r_{i,j-\half}$ and the integrand.

\subsection{Change in control volume momentum}
A discretization of the the time derivative term in Eq.~\eqref{eq_OurForce}
for the integral over $\Vcvt$ at time step $n+1$ is given by

\begin{equation}
\label{eq_temporal_discrete_CV_momentum}
-\int_{\Vcv^{n+1}} \rho \D{\u}{t} \dV  \approx 
-\int_{\Vcv^{n+1}} \rho \frac{\u^{n+1} - \u^{n}}{\dt} \dV
= -\frac{1}{\dt}\int_{\Vcv^{n+1}} \rho \u^{n+1} \dV + \frac{1}{\dt}\int_{\Vcv^{n+1}} \rho \u^{n} \dV,
\end{equation}

\noindent where a discrete approximation to the total linear momentum within $\Vcv^{n+1}$ can
be written as
\begin{equation}
\label{eq_discrete_CV_momentum}
\int_{\Vcv^{n+1}} \rho \u \dV \approx 
 \e_x \sum_{(i-\half,j) \in \Vcv^{n+1}} \rho u_{i-\half,j} \delV_{i-\half,j}
+ \e_y \sum_{(i,j-\half) \in \Vcv^{n+1}} \rho v_{i,j-\half} \delV_{i,j-\half}.
\end{equation}

\noindent Here, $\delV = \frac{\dx\dy}{2}$ when either $(i-\half,j) \in \Scv^{n+1}$ or
$(i,j-\half) \in \Scv^{n+1}$, and $\delV = \dx\dy$ otherwise, to ensure that 
$\sum_{(i-\half,j)} \delV_{i-\half,j} = \sum_{(i,j-\half)} \delV_{i,j-\half} = \left|\Vcv^{n+1}\right|$, the
volume of the CV. 

In the original hydrodynamic force formula Eq.~\eqref{eq_NocaForce} 
introduced by Noca \cite{Noca97}, a discretization of the momentum
term is given by

\begin{equation}
\label{eq_temporal_discrete_CV_momentum_Noca}
-\d{}{t}\int_{\Vcv^{n+1}} \rho \u \dV  \approx 
 -\frac{1}{\dt}\int_{\Vcv^{n+1}} \rho \u^{n+1} \dV + \frac{1}{\dt}\int_{\Vcv^{n}} \rho \u^{n} \dV.
\end{equation}

\noindent Notice that Eqs.~\eqref{eq_temporal_discrete_CV_momentum} and
\eqref{eq_temporal_discrete_CV_momentum_Noca} are nearly identical, although
the former only requires an evaluation over a single CV, while the latter requires an
evaluation over two time-lagged CVs.

The analogous term in the torque calculation Eq.~\eqref{eq_OurTorque} is discretized
differently. Each $u$ velocity location $\x_{i-\half,j}$ is looped over and an
approximation to $v$ is computed. Then $\r_{i-\half,j} = \x_{i-\half} - \x_0$ is computed
and used in the cross product. Mathematically, this is realized as

\begin{equation}
\label{eq_temporal_discrete_CV_angular_momentum}
-\int_{\Vcv^{n+1}} \rho \r \wedge \D{\u}{t} \dV \approx
-\frac{1}{\dt}\int_{\Vcv^{n+1}} \rho \r \wedge \u^{n+1} \dV + \frac{1}{\dt}\int_{\Vcv^{n+1}} \rho \r \wedge \u^{n} \dV,
\end{equation}  

\noindent in which 
\begin{equation}
\label{eq_discrete_CV_angular_momentum}
\int_{\Vcv^{n+1}} \rho \r \wedge \u \dV \approx 
\sum_{(i-\half,j) \in \Vcv^{n+1}} \r_{i-\half,j} \wedge \left[u_{i-\half,j} \e_x 
+ \frac{v_{i-1,j-\half} + v_{i,j-\half} + v_{i-1,j+\half} + v_{i,j+\half}}{4}\e_y\right] \delV_{i-\half,j}.
\end{equation}

\subsection{Change in body momentum}
The final term that needs to be discretely approximated is the change in momentum
of the immersed body. This is presented as an integral of the body's velocity over the 
region $\Vbt \subset \Omega$ in an Eulerian reference frame. Since the body's position and velocity 
are described in a Lagrangian reference frame over a region $I_b \subset U$, it is generally much easier to evaluate
the Lagrangian form of this integral instead. Using the definition of $\delta(\x)$, it can be shown
that the momentum over these two different reference frames are equivalent:

\begin{equation}
 \int_{\Vb} \rho \u(\x,t) \Dx= 
 \int_{I_b} \rho \U(\s,t) \Ds.
\end{equation}

\noindent Letting $\cG^n$ denote the collection of discrete IB points corresponding to the region $I_b$ 
at time step $n$, the object's momentum is obtained by
\begin{equation}
\Pb^n = \sum_{(l,m) \in \cG^{n}} \rho \U_{l,m}^{n} \Dels_{l,m},
\end{equation}
\noindent in which $\U_{l,m}^n$ denotes the velocity of IB node $(l,m)$ at time step $n$,
and $\Dels_{l,m}$ denotes the discrete volume occupied by the node. 
The change in momentum required for the evaluation of hydrodynamic forces is then
given by
  
\begin{equation}
\d{}{t} \int_{I_b} \rho \U(\s,t) \Ds \approx 
\frac{\Pb^{n+1} - \Pb^n}{\dt}.
\end{equation}

\noindent The change in the body's angular momentum for the torque calculation is done similarly.


\section*{Bibliography}
\begin{flushleft}
 \bibliography{CV_Drag}
\end{flushleft}

\end{document}